%% file: main.tex
\documentclass{article}
\usepackage[twoside,
paperwidth=210mm,
paperheight=297mm,
textheight=622pt,
textwidth=468pt,
centering,
headheight=50pt,
headsep=12pt,
footskip=18pt,
footnotesep=24pt plus 2pt minus 12pt,
columnsep=2pc]{geometry}

\usepackage[auth-sc]{authblk}

\usepackage{mathtools}

\allowdisplaybreaks

\usepackage{amssymb}
\usepackage[mathscr]{eucal}

\usepackage{mismath}

\usepackage{hyperref}

\usepackage{bm} 

\usepackage{graphicx}
\graphicspath{{figs/}}

\usepackage{pgf}
\usepackage{pgfplots}
\pgfplotsset{compat=1.14}

\usepackage{tikz}
\usetikzlibrary{arrows.meta}
\tikzset{%
  >={Latex[width=2mm,length=2mm]},
            base/.style = {rectangle, rounded corners, draw=black,
                           minimum width=2cm, minimum height=0.75cm,
                           text centered, font=\sffamily},
  activityStarts/.style = {base, fill=blue!30},
       startstop/.style = {base, fill=red!30},
    activityRuns/.style = {base, fill=green!30},
         process/.style = {base, minimum width=2.5cm, fill=orange!15},
}

\usetikzlibrary{math}
\usetikzlibrary{calc}

\usepackage{subcaption}

\usepackage{multirow}

\usepackage{makecell}
\setcellgapes{2pt}

\usepackage{booktabs}

\usepackage{siunitx}
\usepackage{algorithm}
\usepackage{algpseudocode}
\algrenewcommand\algorithmiccomment[2][\normalsize]{{#1\hfill #2}}

\usepackage[capitalise]{cleveref}

\usepackage{easyReview}

\newcommand{\Real}{\mathbb{R}}

\newcommand{\Integer}{\mathbb{Z}}

\newcommand{\SmallHalf}{{\scriptscriptstyle \tfrac{1}{2}}}

\DeclareMathOperator{\Div}{div}

\usepackage{xcolor}
\definecolor{MyColor1}{HTML}{BFCCB5}
\definecolor{MyColor2}{HTML}{7C96AB}
\definecolor{MyColor3}{HTML}{B7B7B7}
\definecolor{MyColor4}{HTML}{EDC6B1}

\DeclarePairedDelimiter{\RoundBrackets}{(}{)}
\DeclarePairedDelimiter{\CurlyBrackets}{\{}{\}}

\usepackage{stmaryrd}

\DeclareSymbolFont{sfletters}{OML}{cmbrm}{m}{it}

\DeclareMathSymbol{\salpha}{\mathord}{sfletters}{"0B}
\DeclareMathSymbol{\sbeta}{\mathord}{sfletters}{"0C}
\DeclareMathSymbol{\sgamma}{\mathord}{sfletters}{"0D}
\DeclareMathSymbol{\sdelta}{\mathord}{sfletters}{"0E}
\DeclareMathSymbol{\sepsilon}{\mathord}{sfletters}{"0F}
\DeclareMathSymbol{\szeta}{\mathord}{sfletters}{"10}
\DeclareMathSymbol{\seta}{\mathord}{sfletters}{"11}
\DeclareMathSymbol{\stheta}{\mathord}{sfletters}{"12}
\DeclareMathSymbol{\siota}{\mathord}{sfletters}{"13}
\DeclareMathSymbol{\skappa}{\mathord}{sfletters}{"14}
\DeclareMathSymbol{\slambda}{\mathord}{sfletters}{"15}
\DeclareMathSymbol{\smu}{\mathord}{sfletters}{"16}
\DeclareMathSymbol{\snu}{\mathord}{sfletters}{"17}
\DeclareMathSymbol{\sxi}{\mathord}{sfletters}{"18}
\DeclareMathSymbol{\spi}{\mathord}{sfletters}{"19}
\DeclareMathSymbol{\srho}{\mathord}{sfletters}{"1A}
\DeclareMathSymbol{\ssigma}{\mathord}{sfletters}{"1B}
\DeclareMathSymbol{\stau}{\mathord}{sfletters}{"1C}
\DeclareMathSymbol{\supsilon}{\mathord}{sfletters}{"1D}
\DeclareMathSymbol{\sphi}{\mathord}{sfletters}{"1E}
\DeclareMathSymbol{\schi}{\mathord}{sfletters}{"1F}
\DeclareMathSymbol{\spsi}{\mathord}{sfletters}{"20}
\DeclareMathSymbol{\somega}{\mathord}{sfletters}{"21}
\DeclareMathSymbol{\svarepsilon}{\mathord}{sfletters}{"22}
\DeclareMathSymbol{\svartheta}{\mathord}{sfletters}{"23}
\DeclareMathSymbol{\svarpi}{\mathord}{sfletters}{"24}
\DeclareMathSymbol{\svarrho}{\mathord}{sfletters}{"25}
\DeclareMathSymbol{\svarsigma}{\mathord}{sfletters}{"26}
\DeclareMathSymbol{\svarphi}{\mathord}{sfletters}{"27}
\DeclareMathSymbol{\sDelta}{\mathord}{sfletters}{"01}
\DeclareMathSymbol{\sTheta}{\mathord}{sfletters}{"02}

\usepackage{amsthm}

\theoremstyle{definition}

\crefname{assumption}{assumption}{assumptions}
\Crefname{assumption}{Assumption}{Assumptions}

\crefname{problem}{problem}{problems}
\Crefname{problem}{Problem}{Problems}

\theoremstyle{remark}

\usepackage{lineno}

\title{A robust solver for large-scale heat transfer topology optimization}
\author[1]{Yingjie Zhou}
\author[1]{Changqing Ye}
\author[1]{Yucheng Liu}
\author[2]{Shubin Fu}
\author[1]{Eric~T.~Chung}
\affil[1]{Department of Mathematics, The Chinese University of Hong Kong, Shatin, Hong~Kong~SAR, China.}
\affil[2]{Eastern Institute for Advanced Study, Ningbo, China.}

\date{}
 
\begin{document}
\maketitle

\begin{abstract}
This paper presents a large-scale parallel solver, specifically designed to tackle the challenges of solving high-dimensional and high-contrast linear systems in heat transfer topology optimization. The solver incorporates an interpolation technique to accelerate convergence in high-resolution domains, along with a multiscale multigrid preconditioner to handle complex coefficient fields with significant contrast. All modules of the optimization solver are implemented on a high performance computing cluster by the PETSc numerical library. Through a series of numerical investigations, we demonstrate the effectiveness of our approach in enhancing convergence and robustness during the optimization process, particularly in high-contrast scenarios with resolutions up to $1024^3$. Our performance results indicate that the proposed preconditioner achieves over $2\times$ speedup against the default algebraic multigrid in PETSc for high-contrast cases.


\end{abstract}

\section{Introduction}

With the simultaneous pursuit of performance and portability, electronic devices are required to integrate more circuits within constrained spaces. Accordingly, the heat generated by these devices is becoming a critical concern, as it can lead to overheating and adversely affect both the performance and reliability of the devices. The ongoing increase in power density, coupled with limited surface area, poses significant challenges for heat dissipation, especially for conventional cooling systems such as heat pipe radiators and heat sinks. As a result, numerous studies have focused on the optimum design of the material distribution to enhance heat transfer from interior regions of the devices to their surface areas.

Bejan et al.\ \cite{bejan1997constructal} first introduced a fundamental problem concerning the optimal distribution of limited materials for effective heat dissipation to surfaces. They approached this challenge by optimizing the assembly of some predefined geometric shapes, which is conducive to large-scale industrial manufacturing. However, achieving the optimal design remains difficult \cite{boichot2009tree} due to excessive constraints. As 3D printing technology matures, more refined microstructures can be fabricated to enhance thermal conduction, thereby promoting the application of topology optimization (TO) in the context of heat conduction systems \cite{bendsoe2013topology}.

Many topology optimization methods have been developed for thermal conduction systems according to review articles \cite{dbouk2017review,fawaz2022topology}. Zhang and Liu \cite{zhang2008design}, Dirker and Meyer \cite{dirker2013topology} investigated the optimal design of heat conductive paths by distributing highly conductive materials based on the solid isotropic with material penalization (SIMP) method. Lohan et al.\ \ \cite{lohan2020study} discussed the effect of different practical objectives and constraints to the optimization results. Burger et al.\ \cite{burger2013three} and Dede \cite{dede2009multiphysics} extended the problem to the three-dimensional cases, employing various discretization schemes and optimization targets. 

There are two major obstacles in the process of topology optimization. The first one is the high computational cost, especially when solving the large-scale problem in the three-dimensional space. The computational burden primarily arises from the huge number of design variables and many researchers have devoted considerable efforts to mitigate this issue. Aage et al.\ \cite{aage2015topology, aage2013parallel} implemented the parallel computing methods to accelerate the computation using the Portable and Extendable Toolkit for Scientific Computing (PETSc) \cite{balay2019petsc}. Kim and Yoon \cite{kim2000multi} designed a multi-resolution multiscale method to save the computational cost. Huang et al.\ \cite{huang2023problem} propose a machine learning based enhanced substructure-based framework to extract unique features in high-dimensional space. Several methods have also been proposed to improve the convergence speed and computational time of topology optimization \cite{svanberg2007mma, kim2021generalized}. Following the same objective, a multigrid preconditioned conjugate gradients solver is developed to reduce the iterations of Krylov subspace method by leveraging the geometric properties of the problem  \cite{amir2014multigrid}. Amir et al.\ \cite{amir2010efficient, amir2011reducing, limkilde2018reducing} made progress in establishing an adequate stopping criterion for the state and adjoint equations solver in nested approach. Guo et al.\ propose a new computational framework instead of traditional node point-based solution \cite{guo2016explicit,zhang2017explicit}. The second challenge lies in the high-contrast in thermal conductivity between different materials, which complicates and destabilizes the solution of the linear systems arising from the heat transfer model \cite{aage2013parallel}. Researches on this issue remains limited. Lazarov firstly investigated the effectiveness of multiscale finite element methods in addressing such high-contrast matrices in the topology optimization \cite{lazarov2013topology}.
 
Multiscale techniques, including the multiscale finite element method (MsFEM) \cite{hou1997multiscale} and the generalized multiscale finite element method (GMsFEM) \cite{efendiev2013generalized,chung2023multiscale,chung2015mixed} are proven capable of handling extremely complicated systems arising from highly heterogeneous media. The central idea of these methods is to reduce the problem defined on fine grids by translating it to a coarser grid. The multiscale basis functions capture the heterogeneity of the original medium and are constructed by solving a set of carefully designed local problems within each coarse element. While these multiscale methods have been successfully applied to a variety of multiscale models, their effectiveness tends to diminish as the contrast in permeability and correlation length increases \cite{arbogast2015two}. Besides, in the context of topology optimization, it is crucial to develop an efficient solver that can directly process the original fine-scale data, rather than relying exclusively on multiscale methods. Moreover, when adjoint approaches \cite{giles2000introduction} are utilized to compute the gradients of the objective function, we normally need to solve multiple linear systems with different sources within each optimization iteration. A robust preconditioner to solve the original fine-scale problem is essential to ensure the efficiency of the optimization method.
Direct solver such as MUMPS \cite{amestoy2001fully} and SuperLU \cite{li2005overview}, although better suited for solving such problems with multiple sources, face challenges in large-scale applications due to low parallel efficiency and excessive memory consumption. Therefore, adopting preconditioned iterative solvers is preferred in topology optimization problems. 

To addressing the aforementioned obstacles, employing multiscale techniques to develop classical preconditioners is an effective strategy to accelerate iterative solvers \cite{fu2024efficient, fu2024adaptive}. By engineering multiscale coarse spaces, we can effectively capture the high-contrast features inherent in the problem, providing a robust solution to the second challenge \cite{ye2024robust}. Besides, applying parallel computing technique is possible since generalized eigenvalue problem for building the coarse space is solved in a communication-free manner  \cite{fu2024efficient}. This capability significantly alleviates the computational burden associated with high-dimensional problems, thereby addressing the first issue of high computational cost.

In this study, we implement a parallel solver for the topology optimization of heat conduction, focusing on high-resolution and high-contrast simulations. To ensure rapid convergence in high-resolution cases, we first address the problem in corresponding low-resolution scenarios and gradually interpolate the results into the original domains to achieve satisfactory outcomes. Additionally, we employ a robust multigrid preconditioner based on solving spectral problems within a nested multiscale space. To enhance the efficiency of this solver, all modules of the topology optimization process, including the application of the Method of Moving Asymptotes \cite{svanberg1987method}, and the construction of our preconditioner, are parallelized by PETSc. We use PCMG, an interface built within PETSc to manage smoothers and interpolation operators for the multigrid technique, resulting in significant improvements in computational efficiency compared to our previous work \cite{Ye2024}. Numerous experiments are conducted to elucidate the efficiency of our solver for many different conditions. In particular, we examine the robustness of our preconditioner against contrast ratios of media, comparing its performance with the default algebraic multigrid within PETSc.

The rest of this paper is organized as follows.
In \cref{sec:preliminaries}, we introduce the model, discretization and optimization technique.
\Cref{sec:methods} constitutes the core of this paper, where we develop a robust multigrid preconditioner based on solving spectral problems in a nested multiscale space.
Numerical experiments conducted with various parameters are reported in \cref{sec:experiments}. Finally, \cref{sec:conclusions} provides concluding remarks.

\section{Preliminaries}\label{sec:preliminaries}
\subsection{Steady state heat transfer topology optimization}
\begin{figure}[!ht]
  \def\OverLen{1.0}
  \def\LENGTH{3.0}
  \def\Opacity{0.25}
  \centering
  \begin{tikzpicture}
    \coordinate (CoodOrigin) at ({-\LENGTH / 2}, {\LENGTH / 2}, \LENGTH);
    \draw[->,gray,thick] (CoodOrigin) -- ++(\OverLen, 0, 0) node[below] {$x$};
    \draw[->,gray,thick] (CoodOrigin) -- ++(0, \OverLen, 0) node[left] {$y$};
    \draw[->,gray,thick] (CoodOrigin) -- ++(0, 0, \OverLen) node[below] {$z$};



    \draw[thick,dash dot] (0, 0, 0) -- ++(\LENGTH, 0, 0);
    \draw[thick] (0, \LENGTH, 0) -- ++(\LENGTH, 0, 0);
    \draw[thick,dash dot] (0, 0, 0) -- ++(0, \LENGTH, 0);
    \draw[thick] (\LENGTH, 0, 0) -- ++(0, \LENGTH, 0);
    \draw[thick] (0, 0, \LENGTH) -- ++(\LENGTH, 0, 0);
    \draw[thick] (0, \LENGTH, \LENGTH) -- ++(\LENGTH, 0, 0);
    \draw[thick] (0, 0, \LENGTH) -- ++(0, \LENGTH, 0);
    \draw[thick] (\LENGTH, 0, \LENGTH) -- ++(0, \LENGTH, 0);
    \draw[thick] (\LENGTH, 0, 0) -- ++(0, 0, \LENGTH);
    \draw[thick] (\LENGTH, \LENGTH, 0) -- ++(0, 0, \LENGTH);
    \draw[thick] (0, \LENGTH, 0) -- ++(0, 0, \LENGTH);
    \draw[thick,dash dot] (0, 0, 0) -- ++(0, 0, \LENGTH);
    \draw[<->, thick] (0, \LENGTH, -0.2*\LENGTH) -- ++ (\LENGTH, 0, 0) node[above, midway, font = \small] {$L_x$};
    \draw[thick] (0, \LENGTH, -0.05*\LENGTH) -- ++ (0, 0, -0.2*\LENGTH);
    \draw[thick] (\LENGTH, \LENGTH, -0.05*\LENGTH) -- ++ (0, 0, -0.2*\LENGTH);

    \draw[<->, thick] (-0.1*\LENGTH, \LENGTH, 0) -- ++ (0, 0, \LENGTH) node[above, midway, font = \small] {$L_z$};
    \draw[thick] (-0.02*\LENGTH, \LENGTH, \LENGTH) -- ++ (-0.12*\LENGTH, 0, 0);
    \draw[thick] (-0.02*\LENGTH, \LENGTH, 0) -- ++ (-0.12*\LENGTH, 0, 0);

    \draw[<->, thick] (1.1*\LENGTH, \LENGTH, 0) -- ++ (0, -\LENGTH, 0) node[right, midway, font = \small] {$L_y$};
    \draw[thick] (1.02*\LENGTH, 0, 0) -- ++ (0.13*\LENGTH, 0, 0);
    \draw[thick] (1.02*\LENGTH, \LENGTH, 0) -- ++ (0.13*\LENGTH, 0, 0);

    \draw[fill=gray, opacity=0.8] (0.3*\LENGTH,0,0.3*\LENGTH) -- (0.7*\LENGTH,0,0.3*\LENGTH) -- (0.7*\LENGTH,0,0.7*\LENGTH) -- (0.3*\LENGTH,0,0.7*\LENGTH) -- cycle node[below right] {$\Gamma_\mathup{D}$};

  \end{tikzpicture}
  \caption{The design domain of the model problem.} \label{fig:3dproblem}
\end{figure}
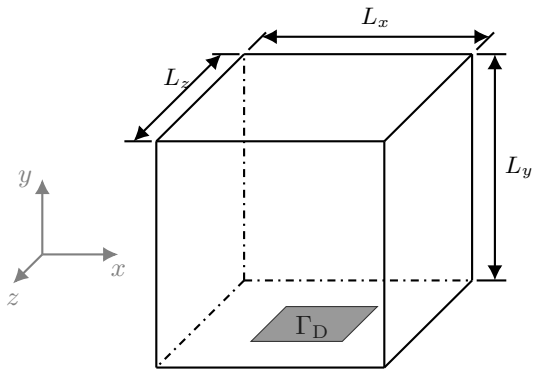

Our optimization problem is defined within a cubic domain $\Omega \subset \mathbb{R}^3$ shown in \cref{fig:3dproblem}, with dimensions $L_x\times L_y \times L_z$. The boundary $\Gamma = \overline{\Gamma_\mathup{D}\cup \Gamma_\mathup{N}}$, $\Gamma_\mathup{D}\cap \Gamma_\mathup{N}=\emptyset$ is partitioned to a Dirichlet (D) and a Neumann (N) part. There are two materials distributed in the domain with different thermal properties. The first one exhibits low thermal conductivity and high volumetric heat generating rate, while the second one features high thermal conductivity and low volumetric heat generating rate.

In order to numerically determine the optimized material distribution, we utilize a structured mesh existing in the domain $\Omega$ for discretization, with a resolution of $M=N_x \times N_y \times N_z$.
Each element has physical dimensions of $h_x \times h_y \times h_z$, where $h_{\diamond} = L_{\diamond} / N_{\diamond}$ for $\diamond \in {x, y, z}$.
As a convention, we denote by $\mathcal{K}_h$ the set of all volume elements, $\mathcal{F}_h$ the set of faces for all elements, $\mathcal{F}_h^\circ$ the subset of $\mathcal{F}_h$ by excluding faces on $\Gamma$. 
The volume element indexed by $(i, j, k) \in \mathcal{I}_h\coloneqq \CurlyBrackets{(i^*,j^*,k^*)\in \Integer^3:0\leq i^*< N_x,\ 0\leq j^* < N_y,\ 0\leq k^*< N_z}$ can be denoted as $K_{i,j,k} \in \mathcal{K}_h$.
Similarly, any face $F \in \mathcal{F}_h$ can be represented as $F_{i\pm \SmallHalf,j,k}$, $F_{i,j\pm \SmallHalf,k}$, or $F_{i,j,k\pm \SmallHalf}$.
Consequently, we treat that $\kappa_{i,j,k}$ within the element $K_{i,j,k}$ remains constant.

In the context of topology optimization, we introduce design variable $x_{i,j,k}\in[0,1]$ in each element to determine the material distribution. The thermal conductivity and heat generating rate in each cell are defined as:
\begin{equation} \label{eq:material}
  \begin{aligned}
    \kappa(x_{i,j,k}) &= \kappa_L + x_{i,j,k}^p (\kappa_H-\kappa_L), \\ 
    f_i &= f_0(1-x_{i,j,k}^p),
  \end{aligned}
\end{equation}
where $\kappa_H$ and $\kappa_L$ are the thermal conductivity for high-conductive material and low-conductive material, $f_0$ is the heat generating rate for low-conductive material and $p\in\mathbb{Z}^+$ is the penalisation factor based on the solid isotropic with material penalisation (SIMP) approach \cite{bendsoe2013topology}. In practice, as the value of 
$p$ increases, the speed at which the design variable $x_{i,j,k}$ converges to 0 or 1 also improves.

Motivated by  \cite{burger2013three,dirker2013topology}, we consider here a generic constrained optimization problem,

\begin{equation}\label{eq:opt}
\min_{\mathbf{x},T} \quad c(\mathbf{x},T)
\end{equation}
subject to 
\[
  \left\{
  \begin{aligned}
      \sum_{i,j,k} x_{i,j,k} / M & \leq V^*, &\forall (i, j, k)&\in \mathcal{I}_h\\
      0\leq x_{i,j,k} & \leq 1, &\forall (i, j, k)& \in \mathcal{I}_h \\
  \end{aligned}
  \right.
\]
Here $\mathbf{x}$ defined by $\{x_{i,j,k}\}_{(i,j,k)\in\mathcal{I}_h}$ is the vector of design variables across the whole domain with the length of $M$. Noted that $V^*$ is the volume constraint on $\mathbf{x}$ that the high conductive material cannot exceed a predefined proportion of all materials. Our purpose is to optimize the material distribution determined by design variable $x_{i,j,k}$ to minimize the objective function $c(\mathbf{x},T)$.


\subsection{Discretization scheme}
To find out the distribution of temperature $T$ in \cref{eq:opt}, we adopt the finite volume method (FVM) since it guarantees element-wise conservation of the fundamental physical quantities. We introduce the flux field $\bm{v}=-\kappa\nabla T$ to get the following first order equation system derived from steady state heat conduction equation:
\begin{equation}\label{eq:mixed system}
  \left\{
  \begin{aligned}
      \kappa^{-1} \bm{v}+\nabla T &= 0  &\text{in}\  &\Omega \\
      \Div \bm{v} &= f  &\text{in}\  &\Omega \\
      \kappa\nabla T \cdot \bm{n}  &= 0 &\text{on}\  &\Gamma_\mathup{N} \\
      T                                &= T_\mathup{D}  &\text{on}\  &\Gamma_\mathup{D}
  \end{aligned}
  \right.
\end{equation}
where $\kappa(\mathbf{x})$ is the heat conduction coefficient, $f$ is the volumetric heat source, $T$ is the temperature field across the whole domain and $T_\mathup{D}$ is the temperature on the Dirichlet boundary. 

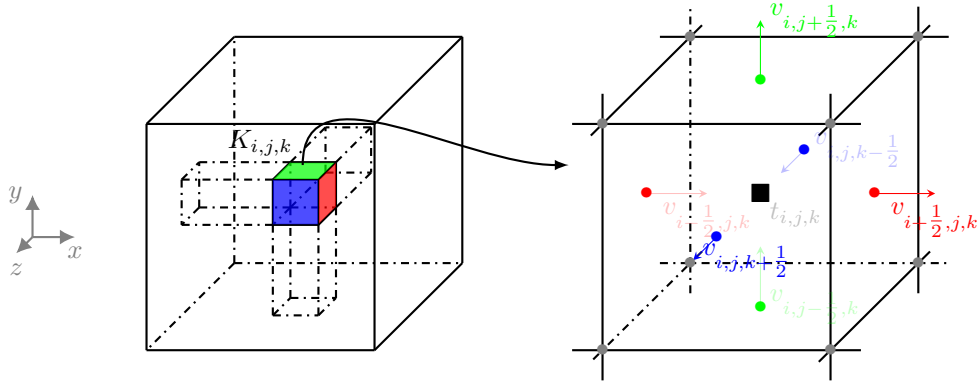
\begin{figure}[!ht]
  \def\OverLen{0.8}
  \def\LENGTH{3.0}
  \def\Opacity{0.25}
  \centering
  \begin{tikzpicture}
    \coordinate (CoodOrigin) at ({-\LENGTH / 2}, {\LENGTH / 2}, \LENGTH);
    \draw[->,gray,thick] (CoodOrigin) -- ++(0.7*\OverLen, 0, 0) node[below] {$x$};
    \draw[->,gray,thick] (CoodOrigin) -- ++(0, 0.7*\OverLen, 0) node[left] {$y$};
    \draw[->,gray,thick] (CoodOrigin) -- ++(0, 0, 0.7*\OverLen) node[below] {$z$};

    \coordinate (izjzkz) at (3*\LENGTH, 0, 0);
    \coordinate (iojzkz) at (2*\LENGTH, 0, 0);
    \coordinate (izjokz) at (3*\LENGTH, \LENGTH, 0);
    \coordinate (iojokz) at (2*\LENGTH, \LENGTH, 0);
    \coordinate (izjzko) at (3*\LENGTH, 0, \LENGTH);
    \coordinate (iojzko) at (2*\LENGTH, 0, \LENGTH);
    \coordinate (izjoko) at (3*\LENGTH, \LENGTH, \LENGTH);
    \coordinate (iojoko) at (2*\LENGTH, \LENGTH, \LENGTH);


    \draw[thick,dash dot] (0, 0, 0) -- ++(\LENGTH, 0, 0);
    \draw[thick] (0, \LENGTH, 0) -- ++(\LENGTH, 0, 0);
    \draw[thick,dash dot] (0, 0, 0) -- ++(0, {\LENGTH}, 0);
    \draw[thick] (\LENGTH, 0, 0) -- ++(0, {\LENGTH}, 0);
    \draw[thick] (0, 0, \LENGTH) -- ++(\LENGTH, 0, 0);
    \draw[thick] (0, \LENGTH, \LENGTH) -- ++(\LENGTH, 0, 0);
    \draw[thick] (0, 0, \LENGTH) -- ++(0, {\LENGTH}, 0);
    \draw[thick] (\LENGTH, 0, \LENGTH) -- ++(0, {\LENGTH}, 0);
    \draw[thick] (\LENGTH, 0, 0) -- ++(0, 0, \LENGTH);
    \draw[thick] (\LENGTH, \LENGTH, 0) -- ++(0, 0, \LENGTH);
    \draw[thick] (0, \LENGTH, 0) -- ++(0, 0, \LENGTH);
    \draw[thick,dash dot] (0, 0, 0) -- ++(0, 0, \LENGTH);

    \draw[thick,dash dot] (2*\LENGTH-0.5*\OverLen, 0, 0) -- ++({\LENGTH+\OverLen}, 0, 0);
    \draw[thick] (2*\LENGTH-0.5*\OverLen, \LENGTH, 0) -- ++({\LENGTH+\OverLen}, 0, 0);
    \draw[thick,dash dot] (2*\LENGTH, {-0.5*\OverLen}, 0) -- ++(0, {{\LENGTH+\OverLen}}, 0);
    \draw[thick] (3*\LENGTH, {-0.5*\OverLen}, 0) -- ++(0, {{\LENGTH+\OverLen}}, 0);
    \draw[thick] (2*\LENGTH-0.5*\OverLen, 0, \LENGTH) -- ++({\LENGTH+\OverLen}, 0, 0);
    \draw[thick] (2*\LENGTH-0.5*\OverLen, \LENGTH, \LENGTH) -- ++({\LENGTH+\OverLen}, 0, 0);
    \draw[thick] (2*\LENGTH, {-0.5*\OverLen}, \LENGTH) -- ++(0, {{\LENGTH+\OverLen}}, 0);
    \draw[thick] (3*\LENGTH, {-0.5*\OverLen}, \LENGTH) -- ++(0, {{\LENGTH+\OverLen}}, 0);
    \draw[thick] (3*\LENGTH, 0, {-0.5*\OverLen}) -- ++(0, 0, {\LENGTH+\OverLen});
    \draw[thick] (2*\LENGTH, \LENGTH, {-0.5*\OverLen}) -- ++(0, 0, {\LENGTH+\OverLen});
    \draw[thick] (3*\LENGTH, \LENGTH, {-0.5*\OverLen}) -- ++(0, 0, {\LENGTH+\OverLen});
    \draw[thick,dash dot] (2*\LENGTH, 0, {-0.5*\OverLen}) -- ++(0, 0, {\LENGTH+\OverLen});

    \draw[thick,dash dot] (0.4*\LENGTH, 0.4*\LENGTH, 0.4*\LENGTH) -- ++(0.2*\LENGTH, 0, 0);
    \draw[thick,dash dot] (0.4*\LENGTH, 0.4*\LENGTH, 0.4*\LENGTH) -- ++(0, 0.2*\LENGTH, 0);
    \draw[thick,dash dot] (0.4*\LENGTH, 0.4*\LENGTH, 0.4*\LENGTH) -- ++(0, 0,0.2*\LENGTH);
    \draw[thick] (0.4*\LENGTH, 0.6*\LENGTH, 0.4*\LENGTH) -- ++(0.2*\LENGTH,0,0);
    \draw[thick] (0.4*\LENGTH, 0.6*\LENGTH, 0.4*\LENGTH) -- ++(0,0,0.2*\LENGTH);
    \draw[thick] (0.4*\LENGTH, 0.4*\LENGTH, 0.6*\LENGTH) -- ++(0.2*\LENGTH,0,0);
    \draw[thick] (0.4*\LENGTH, 0.4*\LENGTH, 0.6*\LENGTH) -- ++(0,0.2*\LENGTH,0);
    \draw[thick] (0.6*\LENGTH, 0.4*\LENGTH, 0.4*\LENGTH) -- ++(0,0.2*\LENGTH,0);
    \draw[thick] (0.6*\LENGTH, 0.4*\LENGTH, 0.6*\LENGTH) -- ++(0,0.2*\LENGTH,0);
    \draw[thick] (0.4*\LENGTH, 0.6*\LENGTH, 0.6*\LENGTH) -- ++(0.2*\LENGTH,0,0);
    \draw[thick] (0.6*\LENGTH, 0.4*\LENGTH, 0.4*\LENGTH) -- ++(0,0,0.2*\LENGTH);
    \draw[thick] (0.6*\LENGTH, 0.6*\LENGTH, 0.4*\LENGTH) -- ++(0,0,0.2*\LENGTH);

    \draw[thick,dash dot] (0.4*\LENGTH, 0, 0.4*\LENGTH) -- ++(0, 0.4*\LENGTH, 0);
    \draw[thick,dash dot] (0.4*\LENGTH, 0, 0.6*\LENGTH) -- ++(0, 0.4*\LENGTH, 0);
    \draw[thick,dash dot] (0.6*\LENGTH, 0, 0.4*\LENGTH) -- ++(0, 0.4*\LENGTH, 0);
    \draw[thick,dash dot] (0.6*\LENGTH, 0, 0.6*\LENGTH) -- ++(0, 0.4*\LENGTH, 0);
    \draw[thick,dash dot] (0.4*\LENGTH, 0, 0.4*\LENGTH) -- ++(0.2*\LENGTH, 0, 0);
    \draw[thick,dash dot] (0.4*\LENGTH, 0, 0.6*\LENGTH) -- ++(0, 0, -0.2*\LENGTH);
    \draw[thick,dash dot] (0.4*\LENGTH, 0, 0.6*\LENGTH) -- ++(0.2*\LENGTH, 0, 0);
    \draw[thick,dash dot] (0.6*\LENGTH, 0, 0.6*\LENGTH) -- ++(0, 0, -0.2*\LENGTH);

    \draw[thick,dash dot] (0, 0.4*\LENGTH, 0.4*\LENGTH) -- ++(0.4*\LENGTH, 0, 0);
    \draw[thick,dash dot] (0, 0.4*\LENGTH, 0.6*\LENGTH) -- ++(0.4*\LENGTH, 0, 0);
    \draw[thick,dash dot] (0, 0.6*\LENGTH, 0.4*\LENGTH) -- ++(0.4*\LENGTH, 0, 0);
    \draw[thick,dash dot] (0, 0.6*\LENGTH, 0.6*\LENGTH) -- ++(0.4*\LENGTH, 0, 0);
    \draw[thick,dash dot] (0, 0.4*\LENGTH, 0.4*\LENGTH) -- ++(0, 0.2*\LENGTH, 0);
    \draw[thick,dash dot] (0, 0.4*\LENGTH, 0.4*\LENGTH) -- ++(0, 0,  0.2*\LENGTH);
    \draw[thick,dash dot] (0, 0.6*\LENGTH, 0.6*\LENGTH) -- ++(0, -0.2*\LENGTH, 0);
    \draw[thick,dash dot] (0, 0.6*\LENGTH, 0.6*\LENGTH) -- ++(0, 0, -0.2*\LENGTH);

    \draw[thick,dash dot] (0.4*\LENGTH, 0.4*\LENGTH, 0) -- ++(0, 0, 0.4*\LENGTH);
    \draw[thick,dash dot] (0.4*\LENGTH, 0.6*\LENGTH, 0) -- ++(0, 0, 0.4*\LENGTH);
    \draw[thick,dash dot] (0.6*\LENGTH, 0.4*\LENGTH, 0) -- ++(0, 0, 0.4*\LENGTH);
    \draw[thick,dash dot] (0.6*\LENGTH, 0.6*\LENGTH, 0) -- ++(0, 0, 0.4*\LENGTH);
    \draw[thick,dash dot] (0.4*\LENGTH, 0.4*\LENGTH, 0) -- ++(0.2*\LENGTH, 0, 0);
    \draw[thick,dash dot] (0.4*\LENGTH, 0.4*\LENGTH, 0) -- ++(0, 0.2*\LENGTH, 0);
    \draw[thick,dash dot] (0.6*\LENGTH, 0.6*\LENGTH, 0) -- ++(0, -0.2*\LENGTH, 0);
    \draw[thick,dash dot] (0.6*\LENGTH, 0.6*\LENGTH, 0) -- ++(-0.2*\LENGTH, 0, 0);

    \draw[fill=green, opacity=0.7] (0.4*\LENGTH,0.6*\LENGTH,0.4*\LENGTH) -- (0.6*\LENGTH,0.6*\LENGTH,0.4*\LENGTH) -- (0.6*\LENGTH,0.6*\LENGTH,0.6*\LENGTH) -- (0.4*\LENGTH,0.6*\LENGTH,0.6*\LENGTH);
    \draw[fill=red, opacity=0.7] (0.6*\LENGTH,0.4*\LENGTH,0.4*\LENGTH) -- (0.6*\LENGTH,0.4*\LENGTH,0.6*\LENGTH) -- (0.6*\LENGTH,0.6*\LENGTH,0.6*\LENGTH) -- (0.6*\LENGTH,0.6*\LENGTH,0.4*\LENGTH);
    \draw[fill=blue, opacity=0.7] (0.4*\LENGTH,0.4*\LENGTH,0.6*\LENGTH) -- (0.6*\LENGTH,0.4*\LENGTH,0.6*\LENGTH) -- (0.6*\LENGTH,0.6*\LENGTH,0.6*\LENGTH) -- (0.4*\LENGTH,0.6*\LENGTH,0.6*\LENGTH);

    \node[red] at (2*\LENGTH, 0.5*\LENGTH, 0.5*\LENGTH) {$\bullet$};
    \draw[red,-stealth, opacity=\Opacity] (2*\LENGTH, 0.5*\LENGTH, 0.5*\LENGTH) -- ++(\OverLen, 0, 0) node[below]  {$v_{i-\SmallHalf,j,k}$};
    \node[red] at (3*\LENGTH, 0.5*\LENGTH, 0.5*\LENGTH) {$\bullet$};
    \draw[red,-stealth] (3*\LENGTH, 0.5*\LENGTH, 0.5*\LENGTH) -- ++(\OverLen, 0, 0) node[below]  {$v_{i+\SmallHalf,j,k}$};

    \node[green] at (2.5*\LENGTH, 0, 0.5*\LENGTH) {$\bullet$};
    \node[green, right, opacity=\Opacity] at (2.5*\LENGTH, 0, 0.5*\LENGTH) {$v_{i, j-\SmallHalf, k}$};
    \draw[green, -stealth, opacity=\Opacity] (2.5*\LENGTH, 0, 0.5*\LENGTH) -- ++(0, \OverLen, 0);
    \node[green] at (2.5*\LENGTH, \LENGTH, 0.5*\LENGTH) {$\bullet$};
    \draw[green,-stealth] (2.5*\LENGTH, \LENGTH, 0.5*\LENGTH) -- ++(0, \OverLen, 0) node[right]  {$v_{i, j+\SmallHalf, k}$};

    \node[blue] at (2.5*\LENGTH, 0.5*\LENGTH, 0) {$\bullet$};
    \node[blue, right, opacity=\Opacity] at (2.5*\LENGTH, 0.5*\LENGTH, 0) {$v_{i, j, k-\SmallHalf}$};
    \draw[blue, -stealth, opacity=\Opacity5] (2.5*\LENGTH, 0.5*\LENGTH, 0) -- ++(0, 0, \OverLen);
    \node[blue] at (2.5*\LENGTH, 0.5*\LENGTH, \LENGTH) {$\bullet$};
    \draw[blue,-stealth] (2.5*\LENGTH, 0.5*\LENGTH, \LENGTH) -- ++(0, 0, \OverLen)node[right]  {$v_{i, j, k+\SmallHalf}$};

    \node at ({2.5*\LENGTH}, {0.5*\LENGTH}, {0.5*\LENGTH}) {$\blacksquare$};
    \node[opacity=\Opacity, below right] at ({2.5*\LENGTH}, {0.5*\LENGTH}, {0.5*\LENGTH}) {$t_{i,j,k}$};

    \draw[-latex,thick](0.9,1.3)node[above left]{$K_{i,j,k}$}
    to[out=90,in=180] (4.4,1.3);

    \node[gray] at (izjzkz) {$\bullet$};
    \node[gray] at (iojzkz) {$\bullet$};
    \node[gray] at (izjokz) {$\bullet$};
    \node[gray] at (iojokz) {$\bullet$};
    \node[gray] at (izjzko) {$\bullet$};
    \node[gray] at (iojzko) {$\bullet$};
    \node[gray] at (izjoko) {$\bullet$};
    \node[gray] at (iojoko) {$\bullet$};

  \end{tikzpicture}
  \caption{An illustration of $\bm{V}_h$ and $W_h$ for element $K_{i,j,k}$.} \label{fig:dof}
\end{figure}

To solve the above equation, the two-point fluxapproximation (TPFA) scheme \cite{barth2018finite} is introduced. The discretization of $T$, denoted by $t_h$, could simply be represented as $t_{i,j,k}$; the discretization of $\bm{v}$, denoted by $\bm{v}_h$ could be expressed as three 3D arrays $v_{i\pm\SmallHalf,j,k}$, $v_{i,j\pm\SmallHalf,k}$ and $v_{i,j,k\pm\SmallHalf}$ defined in the faces of all elements. \cref{fig:dof} provides an elaborate illustration of $\bm{V}_h$ and $W_h$.

Velocity elimination technique \cite{arbogast1997mixed,chen2020generalized,russell1983finite} employs the trapezoidal quadrature rule to remove the flux field in the system:
\begin{equation}\label{eq:flux}
  v_F= -\kappa_F t_F^\delta,\quad \forall F\in \mathcal{F}_h^\circ
\end{equation}
Here $v_F, \kappa_F$ and $t_F^\delta$ have different meanings depending on $F$. For instance, if $F$ is an internal face indexed by $F_{i-\SmallHalf,j,k}$, then $v_F = v_{i-\SmallHalf,j,k}$, $ \kappa_F=2/\RoundBrackets*{1/\kappa_{i-1,j,k}^x+1/\kappa_{i,j,k}^x}$ and $t_F^\delta=\RoundBrackets*{t_{i,j,k}-t_{i-1,j,k}}/h^x$. An advantage of using TPFA lies in that the physical quantity $\{t_{i,j,k}\}$ and the design variable $\{x_{i,j,k}\}$ is indexed in a same manner, thus saving the effort to handle the flux field $\bm{v}$. A linear system for $\mathsf{t}\coloneqq\{t_{i,j,k},(i, j, k)\in \mathcal{I}_h\}$ is then obtained using a standard assembly procedure:
\begin{equation}\label{eq:new linear system}
\mathsf{A}\mathsf{t}=\mathsf{f},
\end{equation}

\subsection{Optimization method and sensitivity analysis}
Many objective functions $c(\mathbf{x}, T)$ in \cref{eq:opt} can be utilized in the TO problem, such as the compliance, the maximum temperature and the mean temperature across the domain. It has been shown that the mean temperature serves as a good indicator of the overall thermal performance of the system \cite{lohan2020study}, while also maintaining relatively low computational costs of its gradients. So we choose 
\begin{equation}\label{eq:cost function}
  c(\mathbf{x}, \mathsf{t})=\bar{T}\approx\frac{1}{M}\sum_{(i,j,k)\in \mathcal{I}_h}T_{i,j,k}
\end{equation}
as our cost function. Notice that \cref{eq:cost function} is a discretized version of $c(\mathbf{x},T)$.

The method of moving asymptotes (MMA) algorithm \cite{svanberg1987method,svanberg2007mma} is introduced and implemented in parallel by PETSc to solve the optimization problem. The sensitivity of the cost function is required by the MMA algorithm. It's computationally expensive to get the gradients of the cost function directly, especially with a large number of design variables $M$. So we use the adjoint approaches \cite{strang2007computational} to get the gradients more efficiently. From $\mathsf{A}(\mathbf{x})\mathsf{t}=\mathsf{f}(\mathbf{x})$, we can derive that
\[
\frac{\partial \mathsf{A}}{\partial x_{i,j,k}}\mathsf{t}+\mathsf{A}\frac{\partial \mathsf{t}}{\partial x_{i,j,k}}=\frac{\partial \mathsf{f}}{\partial x_{i,j,k}}, \quad \forall (i,j,k)\in\mathcal{I}_h
\]
and in the discretized manner, 
\[
    \begin{aligned}
      \frac{\partial c(\mathbf{x}, \mathsf{t})}{\partial x_{i,j,k}} &= \frac{\partial c}{\partial x_{i,j,k}}+\frac{\partial c}{\partial \mathsf{t}}\frac{\partial \mathsf{t}}{\partial x_{i,j,k}} \\
      &=\frac{\partial c}{\partial \mathsf{t}}\cdot(\mathsf{A}^{-1})(\frac{\partial \mathsf{f}}{\partial x_{i,j,k}}-\frac{\partial \mathsf{A}}{\partial x_{i,j,k}}\mathsf{t}) \\
      &=\mathsf{u}^T(\frac{\partial \mathsf{f}}{\partial x_{i,j,k}}-\frac{\partial \mathsf{A}}{\partial x_{i,j,k}}\mathsf{t}), \quad \forall (i,j,k)\in\mathcal{I}_h
    \end{aligned}
\]
The equation can be used to get the sensitivity of each element, thus it forms a matrix of sensitivities with the same distribution of $\mathbf{x}$ or $\mathsf{t}$. The adjoint vector $\mathsf{u}$ can be calculated by:
\begin{equation}\label{eq:adjoint}
\mathsf{A}\mathsf{u} = \frac{\partial c}{\partial \mathsf{t}}
\end{equation}
since $\mathsf{A}$ is a symmetric matrix.

In each iteration of MMA, two linear systems \cref{eq:new linear system,eq:adjoint} need to be solved, which makes a robust and efficient preconditioner crucial for the optimization process. The next part of this paper is to develop a spectral multigrid preconditioner to accelerate the computation of these two linear systems.

\section{Methods}\label{sec:methods}
We begin by discussing the decomposition of the domain $\Omega$ to outline the construction of coarse spaces.
We assume that $\mathcal{Q}$ layers of hierarchical grids that are denoted by $\mathcal{K}^{(l)}$ with $l\in \CurlyBrackets{0,\dots,\mathcal{Q}-1}$ exists in $\Omega$.
Here hierarchical property here indicates that each $\mathcal{K}^{(i-1)}$ is a refined version of $\mathcal{K}^{(i)}$ for $i=1,\dots,\mathcal{Q}-1$.
It should be noted that the finest mesh corresponds to $\mathcal{K}_h=\mathcal{K}^{(0)}$, and the resolution of the model is directly associated with $\mathcal{K}_h$.

\begin{figure}[!ht]
  \def\DualCoarseLen{1.35}
  \def\CoarseLen{0.45}
  \def\Len{0.15}
  \centering
  \begin{tikzpicture}[scale=1,every node/.style={minimum size=1cm},on grid]
		
    \begin{scope}[
    	yshift=0,every node/.append style={
    	    yslant=0.5,xslant=-1},yslant=0.5,xslant=-1
    	             ]
      \fill[white,fill opacity=.9] (0,0) rectangle (3*\DualCoarseLen,3*\DualCoarseLen);
      \fill[MyColor1,fill opacity=.9] (0,0) rectangle (3*\DualCoarseLen,3*\DualCoarseLen);
      \fill[MyColor2,fill opacity=.9] (2*\DualCoarseLen,0) rectangle (3*\DualCoarseLen,\DualCoarseLen);
      \draw[step=\DualCoarseLen, black] (0,0) grid (3*\DualCoarseLen,3*\DualCoarseLen);
      \draw[gray,very thick] (0,0) rectangle (3*\DualCoarseLen,3*\DualCoarseLen);
    \end{scope}
    	
    \begin{scope}[
    	yshift=60,every node/.append style={
    	yslant=0.5,xslant=-1},yslant=0.5,xslant=-1
    	             ]
    	\fill[white,fill opacity=.9] (0,0) rectangle (3*\DualCoarseLen,3*\DualCoarseLen);
      \fill[MyColor4,fill opacity=0.6] (0,0) rectangle (3*\DualCoarseLen,3*\DualCoarseLen);
      \fill[MyColor1,fill opacity=0.6] (\DualCoarseLen, \DualCoarseLen) rectangle (2*\DualCoarseLen,2*\DualCoarseLen);
      \fill[MyColor2,fill opacity=0.6] (\DualCoarseLen+2*\CoarseLen, \DualCoarseLen) rectangle (2*\DualCoarseLen,\DualCoarseLen+\CoarseLen);

    	\draw[step=\DualCoarseLen, black] (0,0) grid (3*\DualCoarseLen,3*\DualCoarseLen);
      \draw[step=\CoarseLen, black] (\DualCoarseLen,\DualCoarseLen) grid (2*\DualCoarseLen,2*\DualCoarseLen);
    	\draw[black,very thick] (0,0) rectangle (3*\DualCoarseLen,3*\DualCoarseLen);
      \draw[gray,very thick] (\DualCoarseLen, \DualCoarseLen) rectangle (2*\DualCoarseLen,2*\DualCoarseLen);
    	\draw[black,dashed] (0,0) rectangle (3*\DualCoarseLen,3*\DualCoarseLen);

      \fill[MyColor3,fill opacity=0.9] (-\CoarseLen, 0) rectangle (0,3*\DualCoarseLen);
      \draw[step=\CoarseLen, black] (-\CoarseLen-\Len, -\CoarseLen-\Len) grid (0,3*\DualCoarseLen+\CoarseLen+\Len);

      \fill[MyColor3,fill opacity=0.9] (3*\DualCoarseLen, 0) rectangle (3*\DualCoarseLen+\CoarseLen,3*\DualCoarseLen);
      \draw[step=\CoarseLen, black] (3*\DualCoarseLen, -\CoarseLen-\Len) grid (3*\DualCoarseLen+\CoarseLen+\Len,3*\DualCoarseLen+\CoarseLen+\Len);

      \fill[MyColor3,fill opacity=0.9] (0,-\CoarseLen) rectangle (3*\DualCoarseLen, 0);
      \draw[step=\CoarseLen, black] (0, -\CoarseLen-\Len) grid (3*\DualCoarseLen, 0);

      \fill[MyColor3,fill opacity=0.9] (0, 3*\DualCoarseLen) rectangle (3*\DualCoarseLen,3*\DualCoarseLen+\CoarseLen);
      \draw[step=\CoarseLen, black] (0, 3*\DualCoarseLen) grid (3*\DualCoarseLen,3*\DualCoarseLen+\CoarseLen+\Len);

    \end{scope}
    	
    \begin{scope}[
    	yshift=120,every node/.append style={
    	    yslant=0.5,xslant=-1},yslant=0.5,xslant=-1
    	  ]
        \fill[white,fill opacity=.9] (0,0) rectangle (3*\DualCoarseLen,3*\DualCoarseLen);
        \fill[MyColor4,fill opacity=0.9] (\DualCoarseLen, \DualCoarseLen) rectangle (2*\DualCoarseLen,2*\DualCoarseLen);
        \fill[MyColor1,fill opacity=0.9] (\DualCoarseLen+\CoarseLen, \DualCoarseLen+\CoarseLen) rectangle (\DualCoarseLen+2*\CoarseLen, \DualCoarseLen+2*\CoarseLen);
        \fill[MyColor2,fill opacity=0.9] (\DualCoarseLen+\CoarseLen+2*\Len, \DualCoarseLen+\CoarseLen) rectangle (\DualCoarseLen+2*\CoarseLen, \DualCoarseLen+\CoarseLen+\Len);

        \draw[step=\DualCoarseLen, black] (0,0) grid (3*\DualCoarseLen,3*\DualCoarseLen);
        \draw[step=\CoarseLen, black] (\DualCoarseLen,\DualCoarseLen) grid (2*\DualCoarseLen,2*\DualCoarseLen);
        \draw[black,very thick] (\DualCoarseLen,\DualCoarseLen) rectangle (2*\DualCoarseLen,2*\DualCoarseLen);

        \draw[step=\Len, black] (\DualCoarseLen-1.01*\Len,\DualCoarseLen-1.01*\Len) grid (\DualCoarseLen,2*\DualCoarseLen+1.01*\Len);
        \fill[MyColor3,fill opacity=0.6] (\DualCoarseLen-\Len,\DualCoarseLen) rectangle (\DualCoarseLen,2*\DualCoarseLen);

        \draw[step=\Len, black] (\DualCoarseLen,\DualCoarseLen-1.01*\Len) grid (2*\DualCoarseLen,\DualCoarseLen);
        \fill[MyColor3,fill opacity=0.6] (\DualCoarseLen,\DualCoarseLen-1.01*\Len) rectangle (2*\DualCoarseLen,\DualCoarseLen);

        \draw[step=\Len, black] (\DualCoarseLen,2*\DualCoarseLen) grid (2*\DualCoarseLen,2*\DualCoarseLen+\Len);
        \fill[MyColor3,fill opacity=0.6] (\DualCoarseLen,2*\DualCoarseLen) rectangle (2*\DualCoarseLen,2*\DualCoarseLen+\Len);

        \draw[step=\Len, black] (2*\DualCoarseLen,\DualCoarseLen-1.01*\Len) grid (2*\DualCoarseLen+\Len,2*\DualCoarseLen+\Len);
        \fill[MyColor3,fill opacity=0.6] (2*\DualCoarseLen,\DualCoarseLen) rectangle (2*\DualCoarseLen+\Len,2*\DualCoarseLen);

        \draw[step=\Len, black] (\DualCoarseLen+\CoarseLen, \DualCoarseLen+\CoarseLen) grid (\DualCoarseLen+2*\CoarseLen, \DualCoarseLen+2*\CoarseLen);
        \draw[black,dashed] (0,0) rectangle (3*\DualCoarseLen,3*\DualCoarseLen);
    \end{scope}

    \draw[-latex,thick] (5.9,2.4) node[right]{$\tau$}
         to[out=180,in=60] (3,2);

    \draw[-latex,thick](5.9,4.4)node[right]{$K_c$}
        to[out=180,in=60] (0.2,4.2);

    \draw[-latex,thick](5.9,8.4)node[right]{$K_{cc}$}
        to[out=180,in=90] (0.3,6.5);

    \draw[-latex,thick](5.9,6.4)node[right]{$\text{Ghost\ layers}$}
        to[out=180,in=60] (1.2,6.4);

    \draw[-latex,thick](-0.8,6.2)to[out=180,in=90] (-3,4);
    \draw[-latex,thick](-0.8,4.1)to[out=180,in=90] (-3,1.9);
\end{tikzpicture}

  \caption{An illustration of hierarchical meshes, a fine element $\tau$, a coarse element $K_\mathup{c}$, a coarse-coarse element $K_\mathup{cc}$ and ghost layers. In the MPI setting, each MPI process is responsible of a coarse-coarse element in $\mathcal{K}_\mathup{cc}$ along with ghost layers and no inter-process communications are required between coarse elements.}
  \label{fig:grid}
\end{figure}
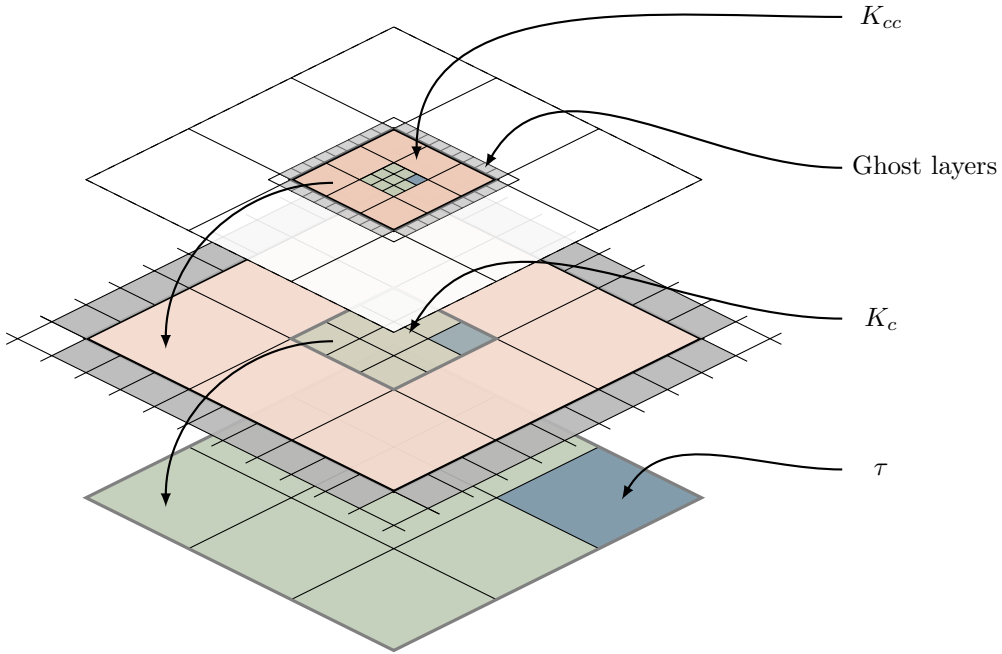

In this paper, we concentrate on a three-grid setting, specifically $\mathcal{Q}=3$.
For clarity, we denote by $\mathcal{K}_\mathup{c}=\mathcal{K}^{(1)}$ the coarse mesh and $\mathcal{K}_\mathup{cc}=\mathcal{K}^{(2)}$ the coarse-coarse mesh.
Additionally, this configuration is tailored for the MPI parallel computing architecture, each MPI process deals with a coarse-coarse element $K_{cc}$ in $\mathcal{K}_\mathup{cc}$. Moreover, each coarse-coarse element $K_{cc}$ is further divided into coarse elements $K_c$ in $\mathcal{K}_c$ and the coarse elements are further divided into fine elements in $\mathcal{K}_h$. Notably, each MPI process (or a coarse-coarse element $K_c$) also handles a ghost layer that have a width of one fine element, facilitating inter-process communications since the elements required by assembling \cref{eq:new linear system} will not exceed one layer.
\Cref{fig:grid} illustrates the architecture of hierarchical meshes, where a fine element $\tau$ in $\mathcal{K}_h$, a coarse element $K_\mathup{c}$ in $\mathcal{K}_\mathup{c}$, a coarse-coarse element $K_\mathup{cc}$ in $\mathcal{K}_\mathup{cc}$ and ghost layers are visualized in different colors.
We denote by $m_\mathup{c}$ the total number of coarse elements, $m_\mathup{cc}$ the total number of coarse-coarse elements, and $M$ the total number of fine elements.

A typical three-grid iteration is demonstrated in \cref{alg:three grid}.
While three-grid or multigrid methods have been proven to converge under certain assumptions \cite{notay2007convergence,xu2022convergence}, we utilize them as acceleration techniques in preconditioned iterative solvers (ref.\  \cite{Saad2003}).
The matrices $\mathsf{S}$ and $\mathsf{S}_\mathup{c}$ serve as smoothers on the fine and coarse grids, respectively. 
The two restriction matrices, $\mathsf{R}_\mathup{c}$ and $\mathsf{R}_\mathup{cc}$, are crucial components of our preconditioner, while their transposes, $\mathsf{R}_\mathup{c}^\intercal$ and $\mathsf{R}_\mathup{cc}^\intercal$, are commonly referred to as prolongation matrices.
Moreover, we have the relations $\mathsf{A}_\mathup{c}=\mathsf{R}_\mathup{c}\mathsf{A}\mathsf{R}_\mathup{c}^\intercal$ and $\mathsf{A}_\mathup{cc}=\mathsf{R}_\mathup{cc}\mathsf{A}_\mathup{c}\mathsf{R}_\mathup{cc}^\intercal=(\mathsf{R}_\mathup{cc}\mathsf{R}_\mathup{c})\mathsf{A}(\mathsf{R}_\mathup{cc}\mathsf{R}_\mathup{c})^\intercal$.

\begin{algorithm}
  \caption{A typical iteration in the three-grid method. }\label{alg:three grid}
  \begin{algorithmic}[1]
    \Require{The operators---$\mathsf{A}$, $\mathsf{A}_\mathup{c}$, $\mathsf{S}^{-1}$, $\mathsf{S}^{-1}_\mathup{c}$, $\mathsf{R}_\mathup{c}$ and $\mathsf{R}_\mathup{cc}$; the initial residual---$\mathsf{r}^{(0)}$}
    \State{Presmoothing: $\mathsf{r}^{(1)}\leftarrow \mathsf{S}^{-1}\mathsf{r}^{(0)}$}  \Comment[\footnotesize]{$\mathsf{S}$ is the smoother on the fine grid}
    \State{Restriction: $\mathsf{r}_\mathup{c}\leftarrow \mathsf{R}_\mathup{c}(\mathsf{r}^{(0)}-\mathsf{A}\mathsf{r}^{(1)})$}  \Comment[\footnotesize]{Restrict the residual onto the coarse space}
    \State{\quad Presmoothing: $\mathsf{r}_\mathup{c}^{(0)}\leftarrow \mathsf{S}_\mathup{c}^{-1}\mathsf{r}_\mathup{c}$}  \Comment[\footnotesize]{$\mathsf{S}_\mathup{c}$ is the smoother on the coarse grid}
    \State{\quad Restriction: $\mathsf{r}_\mathup{cc}\leftarrow \mathsf{R}_\mathup{cc}(\mathsf{r}_\mathup{c}-\mathsf{A}_\mathup{c}\mathsf{r}_\mathup{c}^{(0)})$} \Comment[\footnotesize]{Restrict the residual onto the coarse-coarse space}
    \State{\quad\quad Correction: $\mathsf{e}_\mathup{cc}\leftarrow \mathsf{A}_\mathup{cc}\mathsf{r}_\mathup{cc}$} \Comment[\footnotesize]{Involve a direct solver on the coarse-coarse grid}
    \State{\quad Prolongation: $\mathsf{r}_\mathup{c}^{(1)}\leftarrow \mathsf{r}_\mathup{c}^{(0)}+\mathsf{R}_\mathup{cc}^\intercal\mathsf{e}_\mathup{cc}$ \Comment[\footnotesize]{Project back into the coarse space}}
    \State{\quad Postsmoothing: $\mathsf{e}_\mathup{c}\leftarrow \mathsf{r}_\mathup{c}^{(1)}+\mathsf{S}_\mathsf{c}^{-\intercal}(\mathsf{r}_\mathup{c}-\mathsf{A}_\mathup{c}\mathsf{r}_\mathup{c}^{(1)})$} \Comment[\footnotesize]{Build the map $\mathsf{r}_\mathup{c}\mapsto \mathsf{e}_\mathup{c}$}
    \State{Prolongation: $\mathsf{r}^{(2)}\leftarrow \mathsf{r}^{(1)}+\mathsf{R}_\mathup{c}^\intercal \mathsf{e}_\mathup{c}$}\Comment[\footnotesize]{Project back into the fine space}
    \State{Postsmoothing: $\mathsf{y}\leftarrow \mathsf{r}^{(2)}+\mathsf{S}^{-\intercal}(\mathsf{r}^{(0)}-\mathsf{A}\mathsf{r}^{(2)})$} \Comment[\footnotesize]{The end of the iteration}
  \end{algorithmic}
\end{algorithm}

We then detail the constructions of $\mathsf{R}_\mathup{c}$ and $\mathsf{R}_\mathup{cc}$.
For each coarse element $K_\mathup{c}^i\in \mathcal{K}_\mathup{c}$, we define $W_h(K_\mathup{c}^i)$ by restricting $W_h$ to $K_\mathup{c}^i$ and solving a local spectral problem:
\begin{equation}\label{eq:spepb}
  h_xh_yh_z\sum_{F\in \mathcal{F}_h^\circ(K_\mathup{c}^i)}\kappa_F\Phi_F^\delta q_F^\delta  = \lambda \int_{K_\mathup{c}^i} \kappa \Phi_h q_h \di V, \ \forall q_h \in W_h(K_\mathup{c}^i),
\end{equation}
where $\mathcal{F}_h^\circ(K_\mathup{c}^i)$ represents the set of all internal faces in $K_\mathup{c}^i$, $\Phi_h$ denotes an eigenvector corresponding to the eigenvalue $\lambda$ and $\Phi_F^\delta$,$\ q_F^\delta$ are the jump of $\Phi_h$,$\ q_h$ through face $F$.
After solving the spectral problem \cref{eq:spepb}, we construct the local coarse space $W_\mathup{c}(K_\mathup{c}^i)$ as $\spa\CurlyBrackets{\Phi_{h}^{i,j}\mid j=0,\dots,l_\mathup{c}^i-1}$, where $\CurlyBrackets{\Phi_{h}^{i,j}}_{j=0}^{l_\mathup{c}^{i}-1}$ are the eigenvectors associated with the smallest eigenvalues $L^i_\mathup{c}$.
By identifying $W_\mathup{c}(K_\mathup{c}^i)$ as a linear subspace of $W_h$, the global coarse space $W_\mathup{c}\subset W_h$ is built by $ W_\mathup{c}=W_\mathup{c}(K_\mathup{c}^1)\oplus \dots \oplus W_\mathup{c}(K_\mathup{c}^{m_\mathup{c}})$, and the dimension $n_\mathup{c}$ of $W_\mathup{c}$ is given by $n_\mathup{c}=l_\mathup{c}^1+\dots+l_\mathup{c}^{m_\mathup{c}}$.
Thus, we define $\mathsf{R}_\mathup{c}\in \Real^{n_\mathup{c}\times n}$ as the restriction matrix from the fine grid to the coarse grid, where each row corresponds to the algebraic representation of an eigenvector from \cref{eq:spepb} in $\Real^n$.

To further reduce the dimension of $W_\mathup{c}$, we solve another spectral problem on each coarse-coarse element $K_\mathup{cc}^i \in \mathcal{K}_\mathup{cc}$:
\begin{equation}\label{eq:spectral cc}
  h_xh_yh_z\sum_{F\in \mathcal{F}_h^\circ(K_\mathup{cc}^i)}\kappa_F\varPsi_F^\delta\varTheta_F^\delta = \lambda \int_{K_\mathup{cc}^i} \kappa \varPsi_\mathup{c} \varTheta_\mathup{c} \di V, \ \forall \varTheta_\mathup{c} \in W_\mathup{c}(K_\mathup{cc}^i),
\end{equation}
where $\mathcal{F}_h^\circ(K_\mathup{cc}^i)$ denotes the set of all internal edges in $K_\mathup{cc}^i$, and $W_\mathup{c}(K_\mathup{cc}^i)$ represents the restriction of $W_\mathup{c}$ on $K_\mathup{cc}^i$.
It should be noted that \cref{eq:spectral cc} is defined in the linear space $W_\mathup{c}(K_\mathup{cc}^i)$ rather than $W_h(K_\mathup{cc}^i)$, which offers two advantages:
\begin{itemize}
    \item Solving \cref{eq:spectral cc} results in a much smaller algebraic system and should be significantly easier to solve compared with the latter.
    \item The eigenvectors of \cref{eq:spectral cc} effectively compress the coarse space using the existed information, whereas the latter is unrelated to $W_\mathup{c}$.
\end{itemize}

Similarly, the local coarse-coarse space $W_\mathup{cc}(K_\mathup{cc}^i)$ can be constructed by spanning the set $\CurlyBrackets{\varPsi_\mathup{c}^{i,j}}_{j=0}^{l_\mathup{cc}^i-1}$, which consists of the eigenvectors corresponding to the smallest $l_\mathup{cc}^i$ eigenvalues.
Furthermore, by identifying $W_\mathup{cc}(K_\mathup{cc}^i)$ as a linear subspace of $W_\mathup{c}$, we build the coarse-coarse space $W_\mathup{cc}$ as $W_\mathup{cc}(K_\mathup{cc}^1)\oplus\dots\oplus W_\mathup{cc}(K_\mathup{cc}^{m_\mathup{cc}})$.
The dimension of $W_\mathup{cc}$, denoted by $n_\mathup{cc}$, is given by $l_\mathup{cc}^1+\dots+l_\mathup{cc}^{m_\mathup{cc}}$.
Consequently, we establish a nested inclusions: $W_\mathup{cc}\subset W_\mathup{c} \subset W_h$.
Returning to $\mathsf{R}_\mathup{cc}$, each row of $\mathsf{R}_\mathup{cc}\mathsf{R}_\mathup{c}$ is the algebraic representation eigenvector $\varPsi_c$ of \cref{eq:spectral cc} in $\Real^n$, and $\mathsf{R}_\mathup{cc}$ thus has dimensions of $n_\mathup{cc}\times n_\mathup{c}$.


A key observation from \cref{eq:spepb,eq:spectral cc} is that the right-hand operators in their algebraic forms are all diagonal matrices.
Furthermore, when the eigenvectors obtained from \cref{eq:spepb} are normalized w.r.t.\ the right-hand bilinear form and chosen as the basis for $W_\mathup{c}$, the right-hand operator in the algebraic form of \cref{eq:spectral cc} transforms into an identity matrix. This trick effectively simplifies a generalized eigenvalue problem into a standard one. Although constructing coarse spaces through spectral problems is often considered time-consuming (see \cite{xu2017algebraic}), the proposed spectral problems with a simple format of right-hand matrix can be anticipated to address this challenge.

Apart from the preconditioner discussed above, we employ an interpolation technique to enhance the convergence of MMA iterations. This approach arises from the observed difficulties the MMA algorithm encounters in high-resolution scenarios. To address this issue, we initially carry out the optimization process at a low-resolution level and get a rough result. Following this, we interpolate the results into the high-resolution cases. This two-step strategy significantly improves the performance of the MMA algorithm, allowing it to navigate the complex landscape of high-resolution problems more effectively. By leveraging the insights gained from the low-resolution optimization, we can achieve faster convergence and more stable iterations in the subsequent high-resolution phase. Furthermore, we employ an adaptive algorithm to manage the transformation between grids of varying sizes, which eliminates the need for inter-process communication. As a result, this technique can be applied and parallelized efficiently.

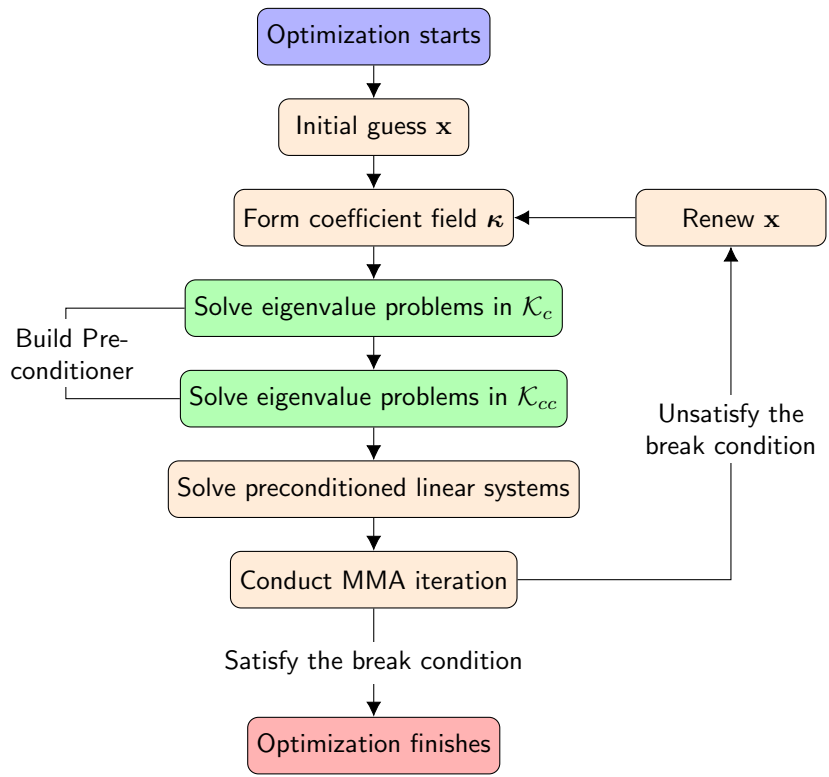
\begin{figure}[!ht]
  \centering
  \begin{tikzpicture}[node distance=1.2cm,
    every node/.style={fill=white, font=\sffamily}, align=center]
  \node (start)             [activityStarts]              {Optimization starts};
  \node (formDesign)     [process, below of=start]          {Initial guess $\mathbf{x}$};
  \node (formCoefficient)      [process, below of=formDesign]   {Form coefficient field $\bm{\kappa}$};
  \node (solveEigen1)     [activityRuns, below of=formCoefficient]   {Solve eigenvalue problems in $\mathcal{K}_c$};
  \node (solveEigen2)      [activityRuns, below of=solveEigen1]
                                                      {Solve eigenvalue problems in $\mathcal{K}_{cc}$};
  \node (solveLinear)      [process, below of=solveEigen2]
                                                                {Solve preconditioned linear systems};
  \node (mma)       [process, below of=solveLinear]
                                                                 {Conduct MMA iteration};
  \node (finish)    [startstop, below of=mma, yshift=-1cm] 
                                                              {Optimization finishes};
  \node (renew)    [process, right of=formCoefficient, xshift=3.5cm]
                                                              {Renew $\mathbf{x}$};
  \draw[->]             (start) -- (formDesign);
  \draw[->]     (formDesign) -- (formCoefficient);
  \draw[->]      (formCoefficient) -- (solveEigen1);
  \draw[->]     (solveEigen1) -- (solveEigen2);
  \draw[->]      (solveEigen2) -- (solveLinear);
  \draw[->]      (solveLinear) -- (mma);
  \draw[->]       (mma) -- node {Satisfy the break condition} (finish);
  \draw[->]    (renew) -- (formCoefficient);
  \draw[->]       (mma) -| node[yshift=2cm, text width=3cm]
                                   {Unsatisfy the break condition}
                                   (renew);

  \draw[-] (solveEigen2.west) -- ++(-1.5,0) -- ++(0,1.2)--                
     node[xshift=-0.7cm,yshift=-0.6cm, text width=2.5cm]
     {Build Preconditioner}(solveEigen1.west);
  \end{tikzpicture}

  \caption{A flow chart of the optimization process.}
  \label{fig:flowchart}
\end{figure}

Finally, a flow chart is presented in \cref{fig:flowchart} to encapsulate the complete process of our solver for topology optimization. It is important to highlight that all steps illustrated in this figure are efficiently parallelized, demonstrating the solver's capability to leverage computational resources effectively. This parallelization relies on that each spectral problem \cref{eq:spepb,eq:spectral cc} are all independent and can be solved at the same time.

\section{Numerical Experiments}\label{sec:experiments}
In this section, we present numerical experiments to demonstrate the effectiveness of the proposed solver in heat transfer topology optimization. We denote our preconditioner as Multiscale Multigrid (MMG) in the following experiments. The unit cube $(0,1)^3$ is set as the design domain $\Omega$ and is further discretized into $M = N_x\times N_y\times N_z$ elements. The topology optimization through MMG is implemented in PETSc \cite{balay2019petsc}. We leverage DMDA, a module provided by PETSc, for inter-process communications and data management. As described in \cref{fig:grid}, each MPI process is responsible for a coarse-coarse element in $\mathcal{K}_\mathup{cc}$, so we have $n_{cc}$ MPI processes in total.\  Moreover, each coarse-coarse element is further divided into $\mathtt{sd}$ parts in all three dimensions, which defines the coarse grid $\mathcal{K}_c$ and determines the relation $n_c = \mathtt{sd}^3\times n_{cc}$. To make the notations more self-explanatory, we take $\mathtt{DoF}$ for the number of fine elements (equal $M$), $\mathtt{proc}$ for the number of MPI processes (equals $n_{cc}$) in the following experiments.

Unlike previous work \cite{Ye2024}, we utilized PCMG, the geometric multigrid interface provided by PETSc, to encapsulate all components in multigrid preconditioner in this paper, resulting in improved computational performance. Another advantage of using PCMG is that we only need to construct $\mathsf{R}_\mathup{c}$ and $\mathsf{R}_\mathup{cc}$ for our preconditioner, a process we refer to as the setup phase. For smoothers $\mathsf{S}$ on $\mathsf{A}$ and $\mathsf{S}_\mathup{c}$ on $\mathsf{A}_\mathup{c}$, we apply block-Jacobi iterations supported by the built-in Krylov subspace methods in PETSc. The full LU factorization of $\mathsf{A}_{cc}$ is facilitated by an external package SuperLU\!\textunderscore DIST \cite{li2003superlu_dist}. We primarily compare our preconditioner with the algebraic multigrid (AMG) preconditioner supported by PETSc with the default parameters, since many other preconditioners like additive Schwarz method (ASM) are unavailable to converge within 10000 iterations in a iterative solver for linear systems when the coefficients field becomes complicated. 

Solving eigenvalue problems is a crucial ingredient in our method, and we employ SLEPc \cite{hernandez2005slepc}, a companion eigensolver for PETSc for implementation. To determine the dimensions of local coarse and coarse-coarse spaces, we use predefined $L_\mathup{c}$ and $L_\mathup{cc}$ eigenvectors for all coarse elements $K_c$ and coarse-coarse elements $K_{cc}$.

Another noteworthy point is that each MMA iteration involves solving two equations \cref{eq:new linear system,eq:adjoint} with the same matrix $\mathsf{A}$ but different right-hand sides. This allows us to reuse the preconditioner for both equations, which is particularly advantageous since our preconditioner is highly efficient for solving linear systems, even though it does not offer significant benefits during the setup phase. For simplicity, we denote \cref{eq:new linear system} as $\mathtt{LS1}$ and \cref{eq:adjoint} as $\mathtt{LS2}$ in the following experiments. So the computational time for solving the linear systems in the topology optimization can be roughly divided into three components: set-up, solve $\mathtt{LS1}$, solve $\mathtt{LS2}$.

A typical feature of the heat transfer topology optimization problem is the high-contrast of the coefficient field. Here we denote the contrast $\kappa_H /\kappa_L$ as $10^\mathtt{cr}$ and choose $\mathtt{cr}$ from $\{1,2,3,4,5,6\}$ to enable a more comprehensive comparison. For other parameters in the optimization problem, we set $T_\mathup{D}=100$ and $f_0=1$.

We choose the PETSc conjugate gradient (CG) method with default parameters to solve the fine system with a relative tolerance of $10^{-6}$ and iteration counts are indicated by $\mathtt{iter}$. The HPC cluster in which we run the program is interconnected by an Infiniband network, and each node is equipped with dual Intel\textsuperscript{\textregistered} Xeon\textsuperscript{\textregistered} Gold 6258R CPUs (56 cores in total) and 192GB of memory. 
Our codes are available on GitHub\footnote{\href{https://github.com/pentaery/TOP_with_MMP}{https://github.com/pentaery/TOP\_with\_MMP}}.

\subsection{Results for topology optimization using the interpolation technique}
\begin{figure}[!ht]
  \centering
  \input{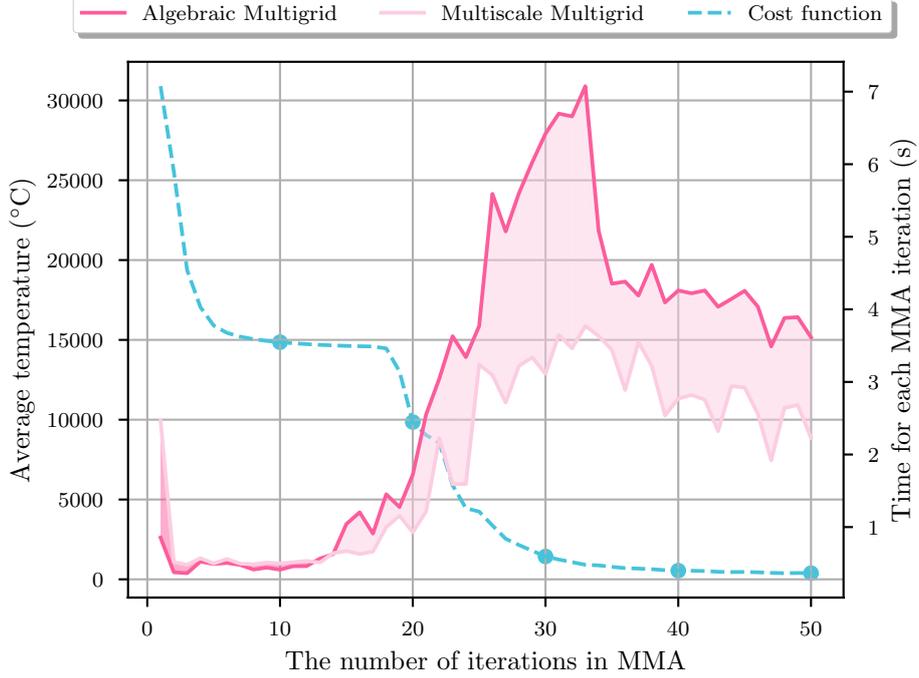}
  \caption{The changes of cost function and number of CG iterations until converged during 50 MMA iterations.} \label{fig:Iterations during MMA}
\end{figure}

In this subsection, we present the results of the steady state heat transfer topology optimization. Frist we give some insights on the implementation of MMA algorithms. A typical MMA iterations executed by 560 MPI processes is shown in \cref{fig:Iterations during MMA}. In this experiment, we set $\mathtt{DoF}$ as $128^3$, $\mathtt{cr}=6$ and $V^*=0.05$. The observed behavior of the cost function demonstrates two significant rapid drops before ultimately converging to an optimal solution. Additionally, the figure presents a comparative analysis of the computational time required for each MMA iteration, contrasting the performance of MMG and AMG. This comparison highlights the efficiency gains achieved through our method, as reflected in the reduced computational time per iteration. 

\begin{figure}[!ht]
  \centering
  \begin{subfigure}[b]{0.18\textwidth}
    \includegraphics[width=\textwidth]{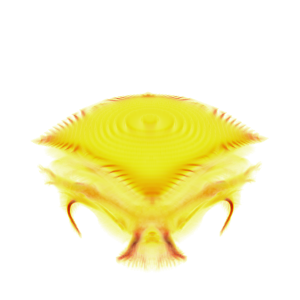}
    \caption{}\label{fig:eff convergence dof a}
  \end{subfigure}
  \begin{subfigure}[b]{0.18\textwidth}
    \includegraphics[width=\textwidth]{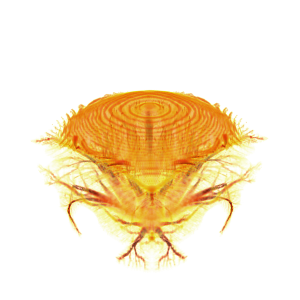}
    \caption{}\label{fig:eff convergence dof b}
  \end{subfigure}
  \begin{subfigure}[b]{0.18\textwidth}
    \includegraphics[width=\textwidth]{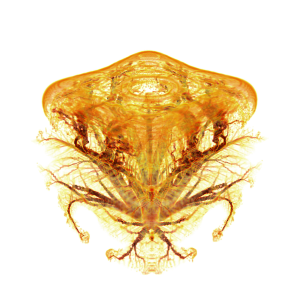}
    \caption{}\label{fig:eff convergence dof c}
  \end{subfigure}
  \begin{subfigure}[b]{0.18\textwidth}
    \includegraphics[width=\textwidth]{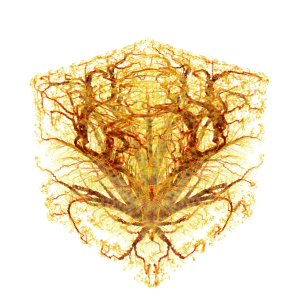}
    \caption{}\label{fig:eff convergence dof d}
  \end{subfigure}
  \begin{subfigure}[b]{0.18\textwidth}
    \includegraphics[width=\textwidth]{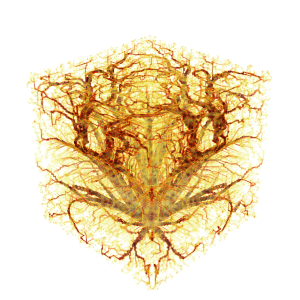}
    \caption{}\label{fig:eff convergence dof e}
  \end{subfigure}
  \caption{Five representative coefficient value fields selected from \cref{fig:Iterations during MMA} every 10 iterations, illustrating the increasing complexity of the problem.}\label{fig:change in TO}
\end{figure}

To illustrate certain details of the MMA iteration process, we pick 5 representative coefficient values every 10 MMA iterations (the blue point in \cref{fig:Iterations during MMA}) and plot the corresponding result in \cref{fig:change in TO}. We observe that our preconditioner exhibits performance comparable to that of AMG when the coefficient field is smooth. However, as the complexity increases, the proposed preconditioner significantly outperforms the AMG. Therefore, a robust preconditioner proves to be highly beneficial in solving multiple linear systems in topology optimization, particularly when there is a large contrast in the coefficient field.

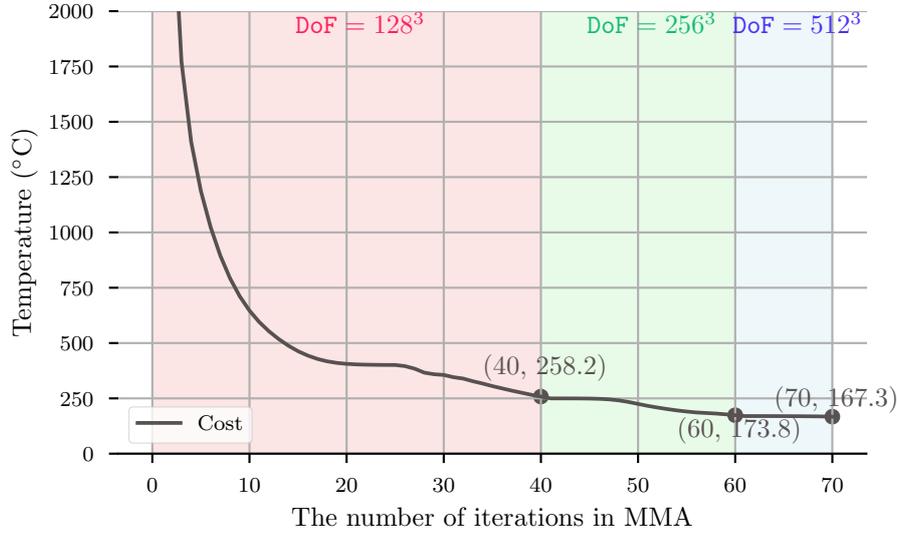
\begin{figure}[!ht]
  \centering
  \input{output.pgf}
  \caption{The changes of cost function in MMA as $\mathtt{DoF}$ transitions from $128^3$ to $512^3$ through two interpolation steps.} \label{fig:Interpolation}
\end{figure}

The convergence rate changes when the resolution level $\mathtt{DoF}$ increases. The optimization method struggles to converge if $\mathtt{DoF}$ is $256^3$ or $512^3$, prompting us to get a initial guess for the next resolution level. To be exact, if we implement MMA for $\mathtt{DoF}=256^3$, we first solve the optimization problem with $\mathtt{DoF}=128^3$ and then interpolate the resulting value to get the initial guess for $\mathtt{DoF}=256^3$. The same procedure is applied to $\mathtt{DoF}=512^3$. As illustrated in \cref{fig:Interpolation}, the cost function changes as $\mathtt{DoF}$ transitions from $128^3$ to $512^3$ through two interpolation steps, with $\mathtt{cr}=6$ and $V^*$ is set to 0.05. After conducting multiple experiments, we have found that this interpolation technique significantly enhances the convergence rate of MMA iterations while reducing computational time. For instance, to  to achieve satisfactory results for $\mathtt{DoF}=512^3$, it's suffices to
execute 40 MMA iterations with $\mathtt{DoF}=128^3$, 20 MMA iterations with $\mathtt{DoF}=256^3$ and finally 10 MMA iterations with $\mathtt{DoF}=512^3$. For a reference, if we conduct the same experiment with $\mathtt{DoF}=512^3$, the cost function decrease to 174.4 after 25 iterations and 129.4 after 70 iterations. However, the computational cost of conducting 25 iterations with $\mathtt{DoF}=512^3$ is much more than using the interpolation technique.

Consequently, we employ interpolation methods in the subsequent experiments when $\mathtt{DoF}$ exceeds $128^3$. In the following subsections, we will present our topology optimization results, highlighting several key parameters involved in the optimization process.

\subsubsection{Comparison between different resolution levels}
\begin{figure}[!ht]
  \centering
  \begin{subfigure}[b]{0.24\textwidth}
    \includegraphics[width=\textwidth]{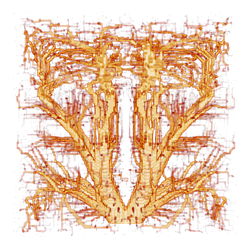}
    \caption{}
  \end{subfigure}
  \begin{subfigure}[b]{0.24\textwidth}
    \includegraphics[width=\textwidth]{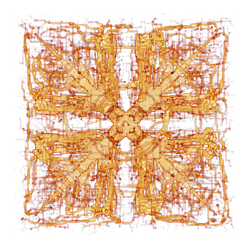}
    \caption{}
  \end{subfigure}
  \begin{subfigure}[b]{0.24\textwidth}
    \includegraphics[width=\textwidth]{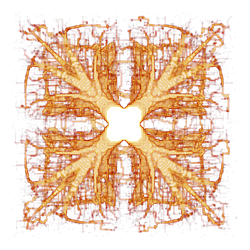}
    \caption{}
  \end{subfigure}
  \begin{subfigure}[b]{0.24\textwidth}
    \includegraphics[width=\textwidth]{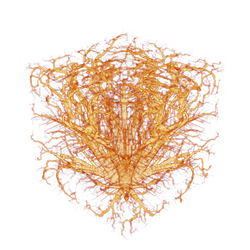}
    \caption{}
  \end{subfigure}

  \begin{subfigure}[b]{0.24\textwidth}
    \includegraphics[width=\textwidth]{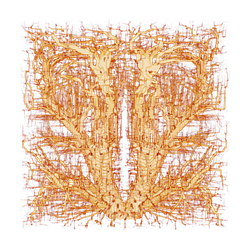}
    \caption{}
  \end{subfigure}
  \begin{subfigure}[b]{0.24\textwidth}
    \includegraphics[width=\textwidth]{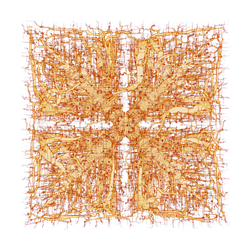}
    \caption{}
  \end{subfigure}
  \begin{subfigure}[b]{0.24\textwidth}
    \includegraphics[width=\textwidth]{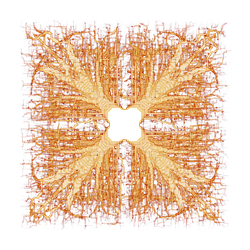}
    \caption{}
  \end{subfigure}
  \begin{subfigure}[b]{0.24\textwidth}
    \includegraphics[width=\textwidth]{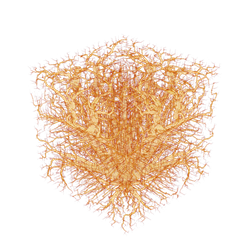}
    \caption{}
  \end{subfigure}

  \begin{subfigure}[b]{0.24\textwidth}
    \includegraphics[width=\textwidth]{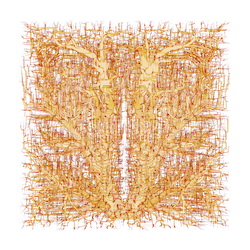}
    \caption{}
  \end{subfigure}
  \begin{subfigure}[b]{0.24\textwidth}
    \includegraphics[width=\textwidth]{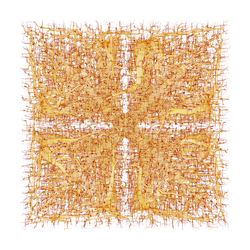}
    \caption{}
  \end{subfigure}
  \begin{subfigure}[b]{0.24\textwidth}
    \includegraphics[width=\textwidth]{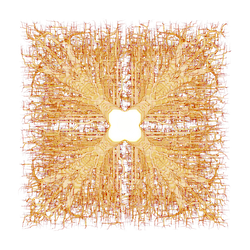}
    \caption{}
  \end{subfigure}
  \begin{subfigure}[b]{0.24\textwidth}
    \includegraphics[width=\textwidth]{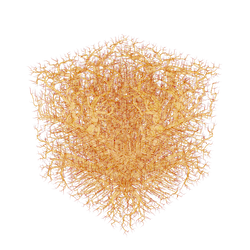}
    \caption{}
  \end{subfigure}
  \caption{Sample solution for $\mathtt{cr}=6$, $V^*=0.05$ with a Dirichlet boundary centered in the bottom surface. (a)-(d) are sequentially side, top, bottom and isometric views for $\mathtt{DoF}=128^3$. Similarly, (e)-(h) are 4 different views for $\mathtt{DoF}=256^3$ and (i)-(l) are 4 different views for $\mathtt{DoF}=512^3$.}\label{fig:differentdof}
\end{figure}

In this study, we investigated the influence of varying resolution levels $\mathtt{DoF}$ on the topological structure of our model. We fixed $\mathtt{cr}=6$ and established a constant volume ratio $V^*=0.05$. The resolution levels $\mathtt{DoF}$ were systematically varied from $128^3$ to $512^3$ and the results are shown in \cref{fig:differentdof}. This range was selected to comprehensively assess how increased resolution affect the model's spatial resolution and the performance of the resulting topological features. 

To ensure consistency across our experiments, the Dirichlet boundary $\Gamma_\mathup{D}$
was fixed at the center of the bottom surface, maintaining a length and width of 
0.1 $\times$ 0.1. This specific configuration allowed us to isolate the effects of resolution while controlling other potentially confounding factors.

The choice of resolution is critical in computational modeling, particularly in applications where fine details significantly impact the overall behavior and properties of the system. By comparing the results obtained at different resolution levels, we aimed to elucidate the relationship between resolution and the accuracy of topological representations.

It can be seen that most of the distribution layouts resemble natural trees, similar to the structures in \cite{burger2013three, dirker2013topology}. The results indicate that as $\mathtt{DoF}$ increases, the complexity and detail of the topological structures become more pronounced, revealing finer features that may be lost at lower resolutions. This suggests that higher resolutions can enhance our understanding of the underlying phenomena, thereby providing deeper insights into the model's behavior. The average temperature of these 3 results are sequentially 536.0$^\circ\mathrm{C}$, 141.3$^\circ\mathrm{C}$ and 137.9$^\circ\mathrm{C}$, which also means that the structure in a high-resolution would have a better heat performance.

\subsubsection{Comparison between different volume ratios}
Another crucial parameter in topology optimization is the volume ratio $V^*$. This parameter plays a significant role in determining the amount of high-conductive materials that can be distributed throughout the entire design domain. In our experimental setup, we fixed $\mathtt{DoF}=512^3$ and $\mathtt{cr}=4$. $V^*$ is then varied from 0.05 to 0.15 to assess its influence on the resulting topological structures. 
\begin{figure}[!ht]
  \centering
  \begin{subfigure}[b]{0.24\textwidth}
    \includegraphics[width=\textwidth]{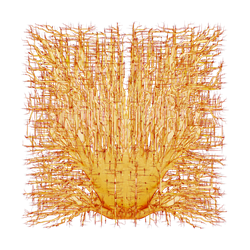}
    \caption{}
  \end{subfigure}
  \begin{subfigure}[b]{0.24\textwidth}
    \includegraphics[width=\textwidth]{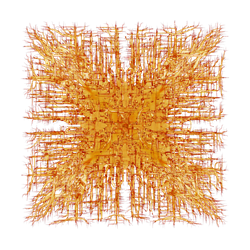}
    \caption{}
  \end{subfigure}
  \begin{subfigure}[b]{0.24\textwidth}
    \includegraphics[width=\textwidth]{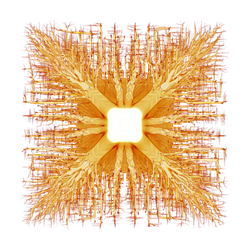}
    \caption{}
  \end{subfigure}
  \begin{subfigure}[b]{0.24\textwidth}
    \includegraphics[width=\textwidth]{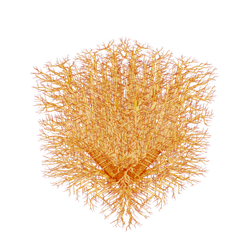}
    \caption{}
  \end{subfigure}

  \begin{subfigure}[b]{0.24\textwidth}
    \includegraphics[width=\textwidth]{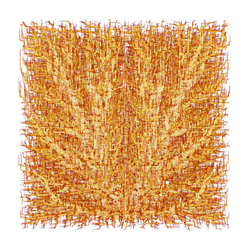}
    \caption{}
  \end{subfigure}
  \begin{subfigure}[b]{0.24\textwidth}
    \includegraphics[width=\textwidth]{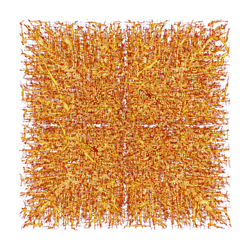}
    \caption{}
  \end{subfigure}
  \begin{subfigure}[b]{0.24\textwidth}
    \includegraphics[width=\textwidth]{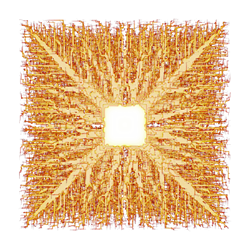}
    \caption{}
  \end{subfigure}
  \begin{subfigure}[b]{0.24\textwidth}
    \includegraphics[width=\textwidth]{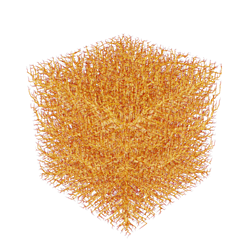}
    \caption{}
  \end{subfigure}

  \begin{subfigure}[b]{0.24\textwidth}
    \includegraphics[width=\textwidth]{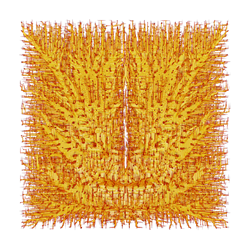}
    \caption{}
  \end{subfigure}
  \begin{subfigure}[b]{0.24\textwidth}
    \includegraphics[width=\textwidth]{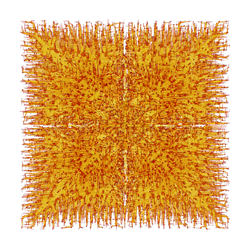}
    \caption{}
  \end{subfigure}
  \begin{subfigure}[b]{0.24\textwidth}
    \includegraphics[width=\textwidth]{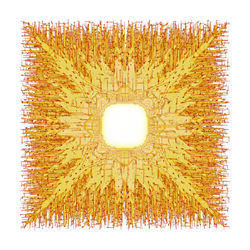}
    \caption{}
  \end{subfigure}
  \begin{subfigure}[b]{0.24\textwidth}
    \includegraphics[width=\textwidth]{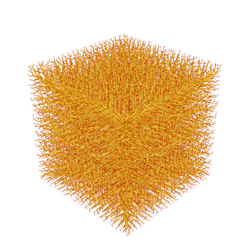}
    \caption{}
  \end{subfigure}

  \caption{Sample solution for $\mathtt{cr}=4$, $\mathtt{DoF}=512^3$ with a Dirichlet boundary centered in the bottom surface. (a)-(d) are sequentially side, top, bottom and isometric views for $V^*=0.05$. Similarly, (e)-(h) are 4 different views for $V^*=0.1$ and (i)-(l) are 4 different views for $V^*=0.15$.}\label{fig:differentv}
\end{figure}

The results in \cref{fig:differentv} show an increase in $V^*$ leads to a more intricate and compact topological structure. Specifically, we recorded the average temperatures for the three different configurations, which are 167.3$^\circ \mathrm{C}$, 111.5$^\circ \mathrm{C}$ and 108.5$^\circ \mathrm{C}$ respectively. These findings indicate a clear trend: as the volume ratio increases, the average temperature within the system decreases, suggesting that a higher allocation of high conductive materials results in more efficient thermal management.

\subsubsection{Comparison between different boundary conditions}

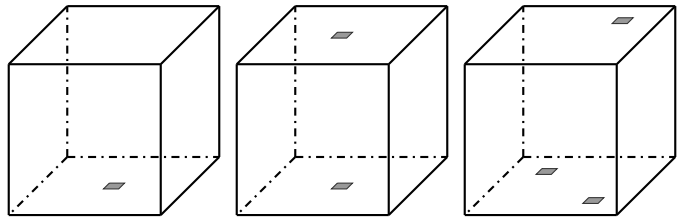
\begin{figure}[!ht]
  \def\OverLen{1.0}
  \def\LENGTH{2.0}
  \def\Opacity{0.25}
  \centering
  \begin{tikzpicture}
    \coordinate (CoodOrigin) at ({-\LENGTH / 2}, {\LENGTH / 2}, \LENGTH);

    \draw[thick,dash dot] (0, 0, 0) -- ++(\LENGTH, 0, 0);
    \draw[thick] (0, \LENGTH, 0) -- ++(\LENGTH, 0, 0);
    \draw[thick,dash dot] (0, 0, 0) -- ++(0, \LENGTH, 0);
    \draw[thick] (\LENGTH, 0, 0) -- ++(0, \LENGTH, 0);
    \draw[thick] (0, 0, \LENGTH) -- ++(\LENGTH, 0, 0);
    \draw[thick] (0, \LENGTH, \LENGTH) -- ++(\LENGTH, 0, 0);
    \draw[thick] (0, 0, \LENGTH) -- ++(0, \LENGTH, 0);
    \draw[thick] (\LENGTH, 0, \LENGTH) -- ++(0, \LENGTH, 0);
    \draw[thick] (\LENGTH, 0, 0) -- ++(0, 0, \LENGTH);
    \draw[thick] (\LENGTH, \LENGTH, 0) -- ++(0, 0, \LENGTH);
    \draw[thick] (0, \LENGTH, 0) -- ++(0, 0, \LENGTH);
    \draw[thick,dash dot] (0, 0, 0) -- ++(0, 0, \LENGTH);

    \draw[fill=gray, opacity=0.8] (0.45*\LENGTH,0,0.45*\LENGTH) -- (0.55*\LENGTH,0,0.45*\LENGTH) -- (0.55*\LENGTH,0,0.55*\LENGTH) -- (0.45*\LENGTH,0,0.55*\LENGTH) -- (0.45*\LENGTH,0,0.45*\LENGTH);

    \draw[thick,dash dot] (1.5*\LENGTH, 0, 0) -- ++(\LENGTH, 0, 0);
    \draw[thick] (1.5*\LENGTH, \LENGTH, 0) -- ++(\LENGTH, 0, 0);
    \draw[thick,dash dot] (1.5*\LENGTH, 0, 0) -- ++(0, \LENGTH, 0);
    \draw[thick] (2.5*\LENGTH, 0, 0) -- ++(0, \LENGTH, 0);
    \draw[thick] (1.5*\LENGTH, 0, \LENGTH) -- ++(\LENGTH, 0, 0);
    \draw[thick] (1.5*\LENGTH, \LENGTH, \LENGTH) -- ++(\LENGTH, 0, 0);
    \draw[thick] (1.5*\LENGTH, 0, \LENGTH) -- ++(0, \LENGTH, 0);
    \draw[thick] (2.5*\LENGTH, 0, \LENGTH) -- ++(0, \LENGTH, 0);
    \draw[thick] (2.5*\LENGTH, 0, 0) -- ++(0, 0, \LENGTH);
    \draw[thick] (1.5*\LENGTH, \LENGTH, 0) -- ++(0, 0, \LENGTH);
    \draw[thick] (2.5*\LENGTH, \LENGTH, 0) -- ++(0, 0, \LENGTH);
    \draw[thick,dash dot] (1.5*\LENGTH, 0, 0) -- ++(0, 0, \LENGTH);

    \draw[fill=gray, opacity=0.8] (1.95*\LENGTH,0,0.45*\LENGTH) -- (2.05*\LENGTH,0,0.45*\LENGTH) -- (2.05*\LENGTH,0,0.55*\LENGTH) -- (1.95*\LENGTH,0,0.55*\LENGTH) -- (1.95*\LENGTH,0,0.45*\LENGTH);

    \draw[fill=gray, opacity=0.8] (1.95*\LENGTH,\LENGTH,0.45*\LENGTH) -- (2.05*\LENGTH,\LENGTH,0.45*\LENGTH) -- (2.05*\LENGTH,\LENGTH,0.55*\LENGTH) -- (1.95*\LENGTH,\LENGTH,0.55*\LENGTH) -- (1.95*\LENGTH,\LENGTH,0.45*\LENGTH);

    \draw[thick,dash dot] (3*\LENGTH, 0, 0) -- ++(\LENGTH, 0, 0);
    \draw[thick] (3*\LENGTH, \LENGTH, 0) -- ++(\LENGTH, 0, 0);
    \draw[thick,dash dot] (3*\LENGTH, 0, 0) -- ++(0, \LENGTH, 0);
    \draw[thick] (4*\LENGTH, 0, 0) -- ++(0, \LENGTH, 0);
    \draw[thick] (3*\LENGTH, 0, \LENGTH) -- ++(\LENGTH, 0, 0);
    \draw[thick] (3*\LENGTH, \LENGTH, \LENGTH) -- ++(\LENGTH, 0, 0);
    \draw[thick] (3*\LENGTH, 0, \LENGTH) -- ++(0, \LENGTH, 0);
    \draw[thick] (4*\LENGTH, 0, \LENGTH) -- ++(0, \LENGTH, 0);
    \draw[thick] (4*\LENGTH, 0, 0) -- ++(0, 0, \LENGTH);
    \draw[thick] (3*\LENGTH, \LENGTH, 0) -- ++(0, 0, \LENGTH);
    \draw[thick] (4*\LENGTH, \LENGTH, 0) -- ++(0, 0, \LENGTH);
    \draw[thick,dash dot] (3*\LENGTH, 0, 0) -- ++(0, 0, \LENGTH);

    \draw[fill=gray, opacity=0.8] (3.2*\LENGTH,0,0.2*\LENGTH) -- (3.2*\LENGTH,0,0.3*\LENGTH) -- (3.3*\LENGTH,0,0.3*\LENGTH) -- (3.3*\LENGTH,0,0.2*\LENGTH) -- (3.2*\LENGTH,0,0.2*\LENGTH);
    \draw[fill=gray, opacity=0.8] (3.7*\LENGTH,0,0.7*\LENGTH) -- (3.7*\LENGTH,0,0.8*\LENGTH) -- (3.8*\LENGTH,0,0.8*\LENGTH) -- (3.8*\LENGTH,0,0.7*\LENGTH) -- (3.7*\LENGTH,0,0.7*\LENGTH);
    \draw[fill=gray, opacity=0.8] (3.7*\LENGTH,\LENGTH,0.2*\LENGTH) -- (3.7*\LENGTH,\LENGTH,0.3*\LENGTH) -- (3.8*\LENGTH,\LENGTH,0.3*\LENGTH) -- (3.8*\LENGTH,\LENGTH,0.2*\LENGTH) -- (3.7*\LENGTH,\LENGTH,0.2*\LENGTH);

  \end{tikzpicture}
  \caption{Three different boundaries set for the experiment, the gray parts are subject to the Dirichlet boundary conditions while the others are subject to Neumann boundary conditions} \label{fig:boundary}
\end{figure}
Different boundary conditions lead to different topology structures. In this experiment, we fixed $\mathtt{DoF}=512^3$ and $\mathtt{cr}=4$. To ensure the diversity of results, we selected three types of Dirichlet boundaries as \cref{fig:boundary} illustrated.
Each piece of Dirichlet part maintains a length and width of 0.1 $\times$ 0.1. The first type of Dirichlet boundary condition is the one utilized in our previous experiments, which is positioned at the center of the bottom surface. This setup serves as a reference for our comparisons. The second type consists of two symmetrical Dirichlet boundaries located on the top and bottom surfaces. The third type features an asymmetrical configuration, where two Dirichlet boundaries are placed on the bottom surface and one on the top surface. This non-uniform distribution of boundary conditions is expected to introduce additional complexity and variation into the optimization results.
\begin{figure}[!ht]
  \centering
  \begin{subfigure}[b]{0.24\textwidth}
    \includegraphics[width=\textwidth]{cr4bot512v10_side.png}
    \caption{}
  \end{subfigure}
  \begin{subfigure}[b]{0.24\textwidth}
    \includegraphics[width=\textwidth]{cr4bot512v10_top.png}
    \caption{}
  \end{subfigure}
  \begin{subfigure}[b]{0.24\textwidth}
    \includegraphics[width=\textwidth]{cr4bot512v10_bot.png}
    \caption{}
  \end{subfigure}
  \begin{subfigure}[b]{0.24\textwidth}
    \includegraphics[width=\textwidth]{cr4bot512v10_isometric.png}
    \caption{}
  \end{subfigure}

  \begin{subfigure}[b]{0.24\textwidth}
    \includegraphics[width=\textwidth]{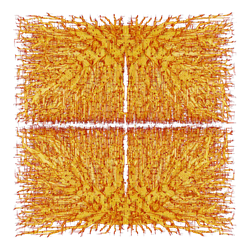}
    \caption{}
  \end{subfigure}
  \begin{subfigure}[b]{0.24\textwidth}
    \includegraphics[width=\textwidth]{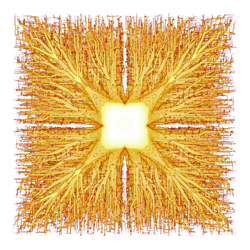}
    \caption{}
  \end{subfigure}
  \begin{subfigure}[b]{0.24\textwidth}
    \includegraphics[width=\textwidth]{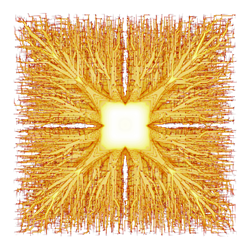}
    \caption{}
  \end{subfigure}
  \begin{subfigure}[b]{0.24\textwidth}
    \includegraphics[width=\textwidth]{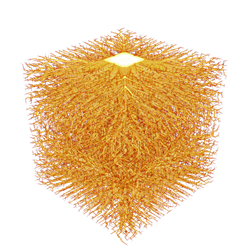}
    \caption{}
  \end{subfigure}

  \begin{subfigure}[b]{0.24\textwidth}
    \includegraphics[width=\textwidth]{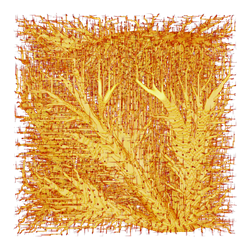}
    \caption{}
  \end{subfigure}
  \begin{subfigure}[b]{0.24\textwidth}
    \includegraphics[width=\textwidth]{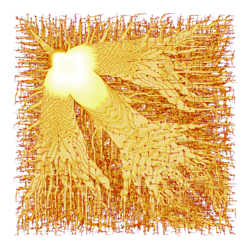}
    \caption{}
  \end{subfigure}
  \begin{subfigure}[b]{0.24\textwidth}
    \includegraphics[width=\textwidth]{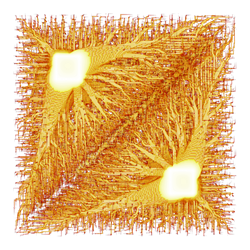}
    \caption{}
  \end{subfigure}
  \begin{subfigure}[b]{0.24\textwidth}
    \includegraphics[width=\textwidth]{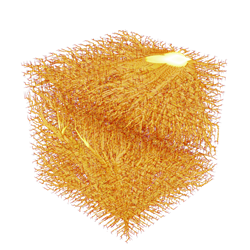}
    \caption{}
  \end{subfigure}
  \caption{Sample solution for $\mathtt{cr}=4$, $\mathtt{DoF}=512^3$ and $V^*=0.1$. (a)-(d) are sequentially side, top, bottom and isometric views for one Dirichlet boundary centered in the bottom surface. Similarly, (e)-(h) are 4 different views for the same Dirichlet boundaries centered in the bottom and top surfaces and (i)-(l) are 4 different views for 3 different Dirichlet boundaries.}\label{fig:differentbd}
\end{figure}

The results are presented in \cref{fig:differentbd}. We can observe that the typical tree-like topology optimization structure emanates from every Dirichlet boundary, intricately spreading into every corner of the domain. The average temperature for these three boundaries are 167.3$^\circ \mathrm{C}$, 105.6$^\circ \mathrm{C}$, 105.5$^\circ \mathrm{C}$ respectively.

\subsubsection{Comparison between different coefficient contrast}
\begin{figure}[!ht]
  \centering
  \begin{subfigure}[b]{0.24\textwidth}
    \includegraphics[width=\textwidth]{cr4bot512v5_side.png}
    \caption{}
  \end{subfigure}
  \begin{subfigure}[b]{0.24\textwidth}
    \includegraphics[width=\textwidth]{cr4bot512v5_top.png}
    \caption{}
  \end{subfigure}
  \begin{subfigure}[b]{0.24\textwidth}
    \includegraphics[width=\textwidth]{cr4bot512v5_bot.png}
    \caption{}
  \end{subfigure}
  \begin{subfigure}[b]{0.24\textwidth}
    \includegraphics[width=\textwidth]{cr4bot512v5_isometric.png}
    \caption{}
  \end{subfigure}

  \begin{subfigure}[b]{0.24\textwidth}
    \includegraphics[width=\textwidth]{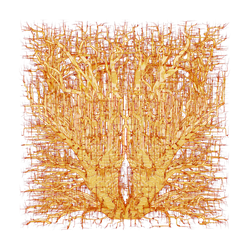}
    \caption{}
  \end{subfigure}
  \begin{subfigure}[b]{0.24\textwidth}
    \includegraphics[width=\textwidth]{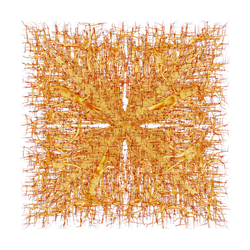}
    \caption{}
  \end{subfigure}
  \begin{subfigure}[b]{0.24\textwidth}
    \includegraphics[width=\textwidth]{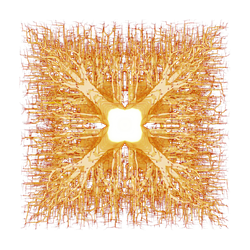}
    \caption{}
  \end{subfigure}
  \begin{subfigure}[b]{0.24\textwidth}
    \includegraphics[width=\textwidth]{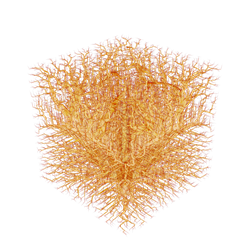}
    \caption{}
  \end{subfigure}

  \begin{subfigure}[b]{0.24\textwidth}
    \includegraphics[width=\textwidth]{cr6bot512v5_side.png}
    \caption{}
  \end{subfigure}
  \begin{subfigure}[b]{0.24\textwidth}
    \includegraphics[width=\textwidth]{cr6bot512v5_top.png}
    \caption{}
  \end{subfigure}
  \begin{subfigure}[b]{0.24\textwidth}
    \includegraphics[width=\textwidth]{cr6bot512v5_bot.png}
    \caption{}
  \end{subfigure}
  \begin{subfigure}[b]{0.24\textwidth}
    \includegraphics[width=\textwidth]{cr6bot512v5_isometric.png}
    \caption{}
  \end{subfigure}
  \caption{Sample solution for $V^*=0.05$, $\mathtt{DoF}=512^3$ with a Dirichlet boundary centered in the bottom surface. (a)-(d) are sequentially side, top, bottom and isometric views for $\mathtt{cr}=4$. Similarly, (e)-(h) are 4 different views for $\mathtt{cr}=5$ and (i)-(l) are 4 different views for $\mathtt{cr}=6$.}\label{fig:differentcr}
\end{figure}

The contrast of the coefficient field is another crucial parameter in topology optimization. In general, a higher contrast can complicate the computation of linear systems. In this study, we set $\mathtt{DoF}=512^3$ and $V^*=0.05$. We then varied the contrast $\mathtt{cr}$ from 4 to 6 to examine its effects on the resulting topological structures. This investigation aims to reveal how changes in contrast impact the characteristics of the final designs. The results are shown in \cref{fig:differentcr}.

The average temperature for these 3 cases are respectively 167.3$^\circ \mathrm{C}$, 138.1$^\circ \mathrm{C}$ and 137.9 $^\circ \mathrm{C}$. While the overall shapes of the resulting topological structures remain relatively similar across the different contrast values, there are noticeable variations in their finer details. Such insights underline the importance of carefully considering the contrast in topology optimization to achieve desired performance outcomes.

\subsection{Scalability tests}

\begin{figure}[!ht]
  \centering
  \input{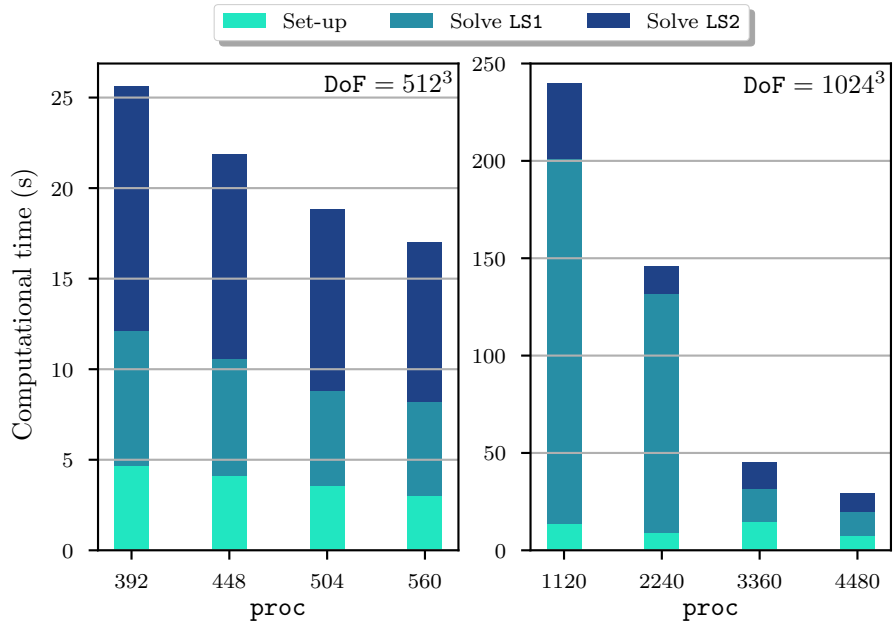}
  \caption{The results of strong scaling test} \label{fig:strong scaling}
\end{figure}

In this subsection, we fix $V^*=0.05$ and $\mathtt{cr}=6$. For our strong scaling analysis, we maintain the problem size at $512^3$ and $1024^3$ while varying the number of computational resources employed. Specifically, we set $(\mathtt{sd}, L_\mathup{c}, L_\mathup{cc})$ as $(7,4,4)$ for $\mathtt{DoF}=512^3$ and $(14,4,17)$ for $\mathtt{DoF}=1024^3$. This experiment allows us to evaluate how effectively the computational workload is distributed across multiple processors, providing insights into the efficiency gains achieved as we increase the number of available resources.  The results, as shown in \cref{fig:strong scaling}, reveal a consistent decrease in computational time as the number of processors increases, particularly evident when $DoF=512^3$. Notably, the right part of the figure indicates an accelerated decline in computational time, especially when the number of processors increases from 2240 to 3360 when $\mathtt{DoF}=1024^3$.

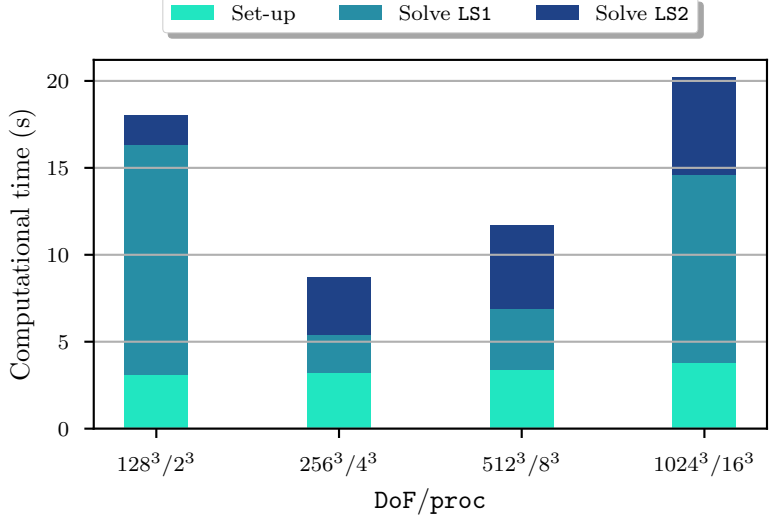
\begin{figure}[!ht]
  \centering
  \input{Weak_Scaling.pgf}
  \caption{The result of weak scaling test} \label{fig:weak scaling}
\end{figure}

In contrast, the weak scaling analysis will focus on increasing the problem size proportionally with the number of processors. Here we adjust the resolution of $\Omega$ while simultaneously increasing the number of computational resources. The coefficient field $\kappa$ for different $\mathtt{DoF}$ is obtained by interpolating $\kappa$ of $\mathtt{DoF}=128^3$ after 20 MMA iterations to the corresponding resolution, thus ensure a rather fair comparison. To ensure that each process in MPI is responsible for a coarse-coarse element $\mathcal{K}_{cc}$ with the the same size, we choose the problem size from $\{128^3, 256^3, 512^3, 1024^3\}$ and $\mathtt{proc}$ from $\{2^3, 4^3, 8^3, 16^3\}$ such that each coarse-coarse element is a $64^3$ cube. Moreover, we set $7^3$ sub domains in each MPI process and fix $L_\mathup{c}=4$, $L_\mathup{cc}=17$. The results in \cref{fig:weak scaling} show that $256^3/4^3$ achieves the best efficiency. This observation may be attributed to the fact that when $\mathtt{proc}\leq 56$, the computing potential of a single node has not yet been fully utilized while communication between MPI processes significantly reduces the efficiency of parallelism when $\mathtt{proc}>56$.

\subsection{Robustness tests for different parameters}
In this subsection, we focus on examining the robustness and efficiency of the proposed preconditioner (MMG). We will frequently write the computational time $\mathtt{time}$ and iteration count $\mathtt{iter}$ of an experiment in a format of $\mathtt{time}(\mathtt{iter})$, where one effective digit is reserved in the unit of seconds for $\mathtt{time}$.
\begin{table}[!ht]
  \footnotesize
  \caption{Records of elapsed wall time and iteration numbers of AMG and the proposed preconditioner w.r.t.\ different parameters ($\mathtt{sd}$, $L_\mathup{c}$, $L_\mathup{cc}$) and fixed contrast $\mathtt{cr}=1$ and $\mathtt{DoF}=512^3$.}\label{tab:comparison cr1}
  \begin{center}
    \makegapedcells
    \begin{tabular}{c c c c c c c c } \toprule
      Preconditioner              &  AMG    & $(7,4,8)$ & $(7,4,11)$  & $(7,4,17)$ & $(7,8,17)$& $(6,4,17)$     & $(8,4,17)$     \\ \midrule
      Set-up                               & $2.7$      & $3.1$      & $3.0$      & $2.9$      & $4.2$     &$3.2$     & $2.9$          \\
      Solve $\mathtt{LS1}$                              & $3.6(9)$      & $2.6(16)$      & $2.3(15)$      & $2.9(13)$     & $2.6(12)$   &$2.5(14)$ & $2.5(13)$ \\
      Solve $\mathtt{LS2}$  & $1.8(15)$      & $3.3(27)$      & $2.6(25)$      & $2.3(22)$     & $3.0(20)$   &$2.8(26)$      & $2.5(23)$         \\
 \midrule
      Relative time in total  & $100\%$      & $111.1\%$      & $97.5\%$           & $100.0\%$ & $121.0\%$  &$104.9\%$      & $97.5\%$         \\
  \bottomrule
    \end{tabular}
  \end{center}
\end{table}

\begin{table}[!ht]
  \footnotesize
  \caption{Records of elapsed wall time and iteration numbers of AMG and the proposed preconditioner w.r.t.\ different parameters ($\mathtt{sd}$, $L_\mathup{c}$, $L_\mathup{cc}$) and fixed contrast $\mathtt{cr}=2$ and $\mathtt{DoF}=512^3$.}\label{tab:comparison cr2}
  \begin{center}
    \makegapedcells
    \begin{tabular}{c c c c c c c c } \toprule
      Preconditioner              &  AMG    & $(7,4,8)$ & $(7,4,11)$  & $(7,4,17)$ & $(7,8,17)$& $(6,4,17)$     & $(8,4,17)$     \\ \midrule
      Set-up                               & $2.8$      & $3.1$      & $3.2$      & $3.2$      & $4.2$     &$3.2$     & $2.9$          \\
      Solve $\mathtt{LS1}$                              & $4.0(11)$      & $2.6(18)$      & $2.5(16)$      & $3.0(14)$     & $3.5(12)$   &$2.9(14)$ & $3.3(13)$ \\
      Solve $\mathtt{LS2}$  & $2.2(19)$      & $3.5(33)$      & $3.1(29)$      & $2.7(25)$     & $2.6(22)$   &$3.4(26)$      & $2.5(24)$         \\
 \midrule
      Relative time in total  & $100\%$      & $102.2\%$      & $97.8\%$           & $98.9\%$ & $111.1\%$  &$105.6\%$      & $96.7\%$         \\
  \bottomrule
    \end{tabular}
  \end{center}
\end{table}

\begin{table}[!ht]
  \footnotesize
  \caption{Records of elapsed wall time and iteration numbers of AMG and the proposed preconditioner w.r.t.\ different parameters ($\mathtt{sd}$, $L_\mathup{c}$, $L_\mathup{cc}$) and fixed contrast $\mathtt{cr}=3$ and $\mathtt{DoF}=512^3$.}\label{tab:comparison cr3}
  \begin{center}
    \makegapedcells
    \begin{tabular}{c c c c c c c c } \toprule
      Preconditioner              &  AMG    & $(7,4,8)$ & $(7,4,11)$  & $(7,4,17)$ & $(7,8,17)$& $(6,4,17)$     & $(8,4,17)$     \\ \midrule
      Set-up                               & $2.9$      & $3.1$      & $2.9$      & $2.9$      & $4.4$     &$3.2$     & $3.0$          \\
      Solve $\mathtt{LS1}$                              & $4.4(14)$      & $3.3(22)$      & $2.8(19)$      & $3.2(16)$     & $3.2(14)$   &$2.9(17)$ & $3.6(17)$ \\
      Solve $\mathtt{LS2}$  & $3.0(27)$      & $4.5(40)$      & $4.8(37)$      & $3.7(29)$     & $3.1(26)$   &$3.6(31)$      & $3.3(31)$         \\
 \midrule
      Relative time in total  & $100\%$      & $105.8\%$      & $101.9\%$           & $95.1\%$ & $103.9\%$  &$92.2\%$      & $93.2\%$         \\
  \bottomrule
    \end{tabular}
  \end{center}
\end{table}

First we fix the resolution level at $\mathtt{DoF}=512^3$ and test the influence of different parameters $(\mathtt{sd}, L_\mathup{c}, L_\mathup{cc})$ to our preconditioner. Here only 1 iteration in MMA is implemented for testing and $\mathtt{cr}$ are chosen from $\{1,2,3,4,5,6\}$ to get a comprehensive analysis. The results are shown in \cref{tab:comparison cr1,tab:comparison cr2,tab:comparison cr3,tab:comparison cr4,tab:comparison cr5,tab:comparison cr6}.

From \cref{tab:comparison cr1,tab:comparison cr2,tab:comparison cr3}, we observe that our preconditioner performs comparably to AMG in terms of both $\mathtt{time}$ and $\mathtt{iter}$ in the low-contrast case. Among all the parameters we tested, $(\mathtt{sd}, L_\mathup{c}, L_\mathup{cc})=(7,4,17)$ demonstrates the best performance, slightly outperforming AMG regardless of $\mathtt{cr}$. Increasing $\mathtt{sd}$ and $L_\mathup{cc}$ proves effective in reducing the number of iterations and can lead to a slight decrease in computational time. Although increasing $L_\mathup{c}$ may also reduce the iteration count, it significantly prolongs the set-up phase, making it less advantageous.

\begin{table}[!ht]
  \footnotesize
  \caption{Records of elapsed wall time and iteration numbers of AMG and the proposed preconditioner w.r.t.\ different parameters ($\mathtt{sd}$, $L_\mathup{c}$, $L_\mathup{cc}$) and fixed contrast $\mathtt{cr}=4$ and $\mathtt{DoF}=512^3$.}\label{tab:comparison cr4}
  \begin{center}
    \makegapedcells
    \begin{tabular}{c c c c c c c c } \toprule
      Preconditioner              &  AMG    & $(7,4,8)$ & $(7,4,11)$  & $(7,4,17)$& $(7,8,17)$ & $(6,4,17)$     & $(8,4,17)$     \\ \midrule
      Set-up                               & $2.7$      & $3.0$      & $2.9$      & $2.9$      & $4.1$     &$4.0$     & $2.9$          \\
      Solve $\mathtt{LS1}$                              & $5.4(24)$      & $3.5(26)$      & $3.3(23)$      & $3.5(20)$     & $3.7(18)$ &$3.2(18)$ & $3.5(20)$ \\
      Solve $\mathtt{LS2}$  & $6.8(65)$      & $6.3(50)$      & $5.0(44)$      & $4.8(39)$     & $4.0(35)$   &$3.6(33)$      & $4.8(38)$         \\
      \midrule
      Relative time in total  & $100\%$      & $85.9\%$      & $75.2\%$           & $75.2\%$ & $79.2\%$  &$72.5\%$      & $75.2\%$         \\
 \bottomrule
    \end{tabular}
  \end{center}
\end{table}

\begin{table}[!ht]
  \footnotesize
  \caption{Records of elapsed wall time and iteration numbers of AMG and the proposed preconditioner w.r.t.\ different parameters ($\mathtt{sd}$, $L_\mathup{c}$, $L_\mathup{cc}$) and fixed contrast $\mathtt{cr}=5$ and $\mathtt{DoF}=512^3$.}\label{tab:comparison cr5}
  \begin{center}
    \makegapedcells
    \begin{tabular}{c c c c c c c c } \toprule
      Preconditioner              &  AMG    & $(7,4,8)$ & $(7,4,11)$  & $(7,4,17)$ & $(7,8,17)$& $(6,4,17)$     & $(8,4,17)$     \\ \midrule
      Set-up                               & $2.8$      & $3.1$      & $2.9$        & $3.1$     & $4.2$    &$3.2$     & $3.1$          \\
      Solve $\mathtt{LS1}$                              & $8.6(60)$      & $3.8(24)$      & $4.0(28)$        & $3.7(24)$  & $4.0(24)$  &$3.2(29)$ & $4.3(25)$ \\
      Solve $\mathtt{LS2}$  & $16(164)$      & $6.1(48)$      & $6.0(54)$          & $5.4(48)$ & $5.9(45)$   &$6.0(54)$      & $6.1(48)$         \\
      \midrule
      Relative time in total  & $100\%$      & $47.1\%$      & $46.7\%$           & $44.2\%$ & $51.1\%$  &$44.9\%$      & $48.9\%$         \\
 \bottomrule
    \end{tabular}
  \end{center}
\end{table}

\begin{table}[!ht]
  \footnotesize
  \caption{Records of elapsed wall time and iteration numbers of AMG and the proposed preconditioner w.r.t.\ different parameters ($\mathtt{sd}$, $L_\mathup{c}$, $L_\mathup{cc}$) and fixed contrast $\mathtt{cr}=6$ and $\mathtt{DoF}=512^3$.}\label{tab:comparison cr6}
  \begin{center}
    \makegapedcells
    \begin{tabular}{c c c c c c c c } \toprule
      Preconditioner             &  AMG    & $(7,4,8)$ & $(7,4,11)$  & $(7,4,17)$& $(7,8,17)$ & $(6,4,17)$     & $(8,4,17)$     \\ \midrule
      Set-up                      & $2.7$      & $3.1$      & $3.0$            & $2.9$ & $4.2$    &$2.9$     & $3.1$          \\
      Solve $\mathtt{LS1}$                              & $14.6(121)$      & $9.5(69)$      & $6.4(50)$          & $6.7(50)$ & $5.6(35)$  &$4.9(36)$ & $5.2(38)$ \\
      Solve $\mathtt{LS2}$  & $32.0(314)$      & $17.7(155)$      & $13.0(111)$          & $12.3(111)$ & $10.0(74)$   &$8.2(74)$      & $9.9(84)$         \\
      \midrule
      Relative time in total  & $100\%$      & $61.5\%$      & $45.4\%$         & $44.4\%$  & $40.2\%$   &$32.5\%$      & $36.9\%$         \\
 \bottomrule
    \end{tabular}
  \end{center}
\end{table}
The results are different in high-contrast scenarios (see \cref{tab:comparison cr4,tab:comparison cr5,tab:comparison cr6}). The combination $(\mathtt{sd}, L_\mathup{c}, L_\mathup{cc})=(6,4,17)$ yields the best performance when $\mathtt{cr}=4$, resulting in substantially fewer iterations and a time savings of 27.5\%. As $\mathtt{cr}$ increases further, the advantage of our method become even more pronounced. Almost all parameters in MMG exhibit a 2$\times$ speedup compared to AMG when $\mathtt{cr}=5$ and some parameters in MMG achieve a 3$\times$ speedup when $\mathtt{cr}=6$. While the iteration numbers $\mathtt{iter}$ are not entirely related to $\mathtt{time}$, they remain well-controlled and show significant improvement over AMG.

\subsection{Applications of our solver to topology optimization}

\begin{figure}[!ht]
  \centering
  \resizebox{\linewidth}{!}{\input{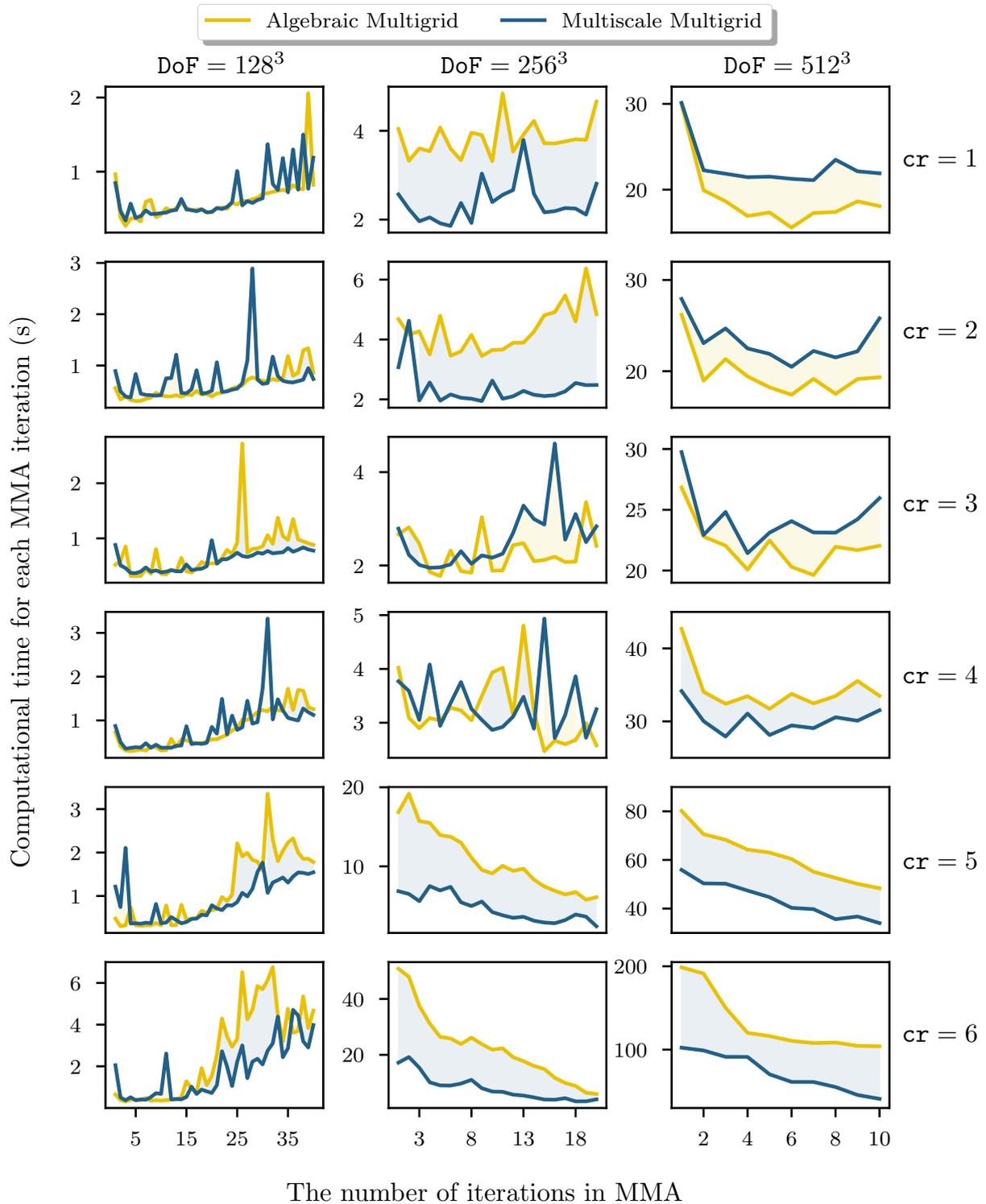}}
  \caption{The computational time of each MMA iterations for a whole optimization process for different $\mathtt{cr}$} \label{fig:MMA test with different contrast}
\end{figure}

Next, we evaluate our solver through a comprehensive optimization process, implementing $40$ MMA iterations for $\mathtt{DoF}=128^3$, $20$ MMA iterations for $\mathtt{DoF}=256^3$ and $10$ MMA iterations for $\mathtt{DoF}=512^3$. The results presented in \cref{fig:MMA test with different contrast} indicates that there are only slight differences between the two approaches, and both demonstrate advantages in different contexts in low-contrast cases. It is noteworthy to mention that the default AMG in PETSc employs a technique of reusing interpolation when rebuilding the precondtioner, which significantly reduces the time in setting the preconditioner. However, even with this advantage, AMG is still far less effective than our preconditioner when $\mathtt{cr}>3$ since last $10$ MMA iterations in $\mathtt{DoF}=512^3$ actually occupy most of our computational resources.

\begin{figure}[!ht]
  \centering
  \resizebox{\linewidth}{!}{\input{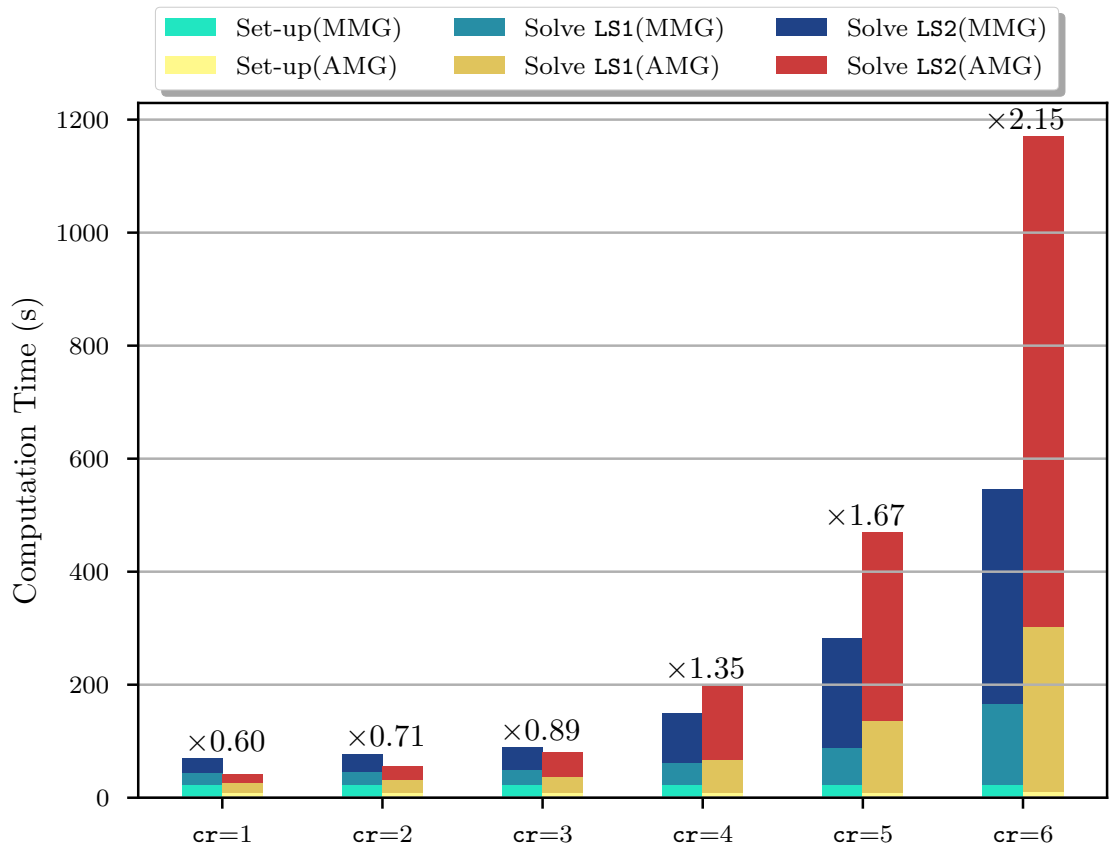}}
  \caption{The computational time for MMA iterations} \label{fig:Total time for different contrast}
\end{figure}

The specific time required to solve the linear system during the $10$ MMA iterations when $\mathtt{DoF}=512^3$ is illustrated in \cref{fig:Total time for different contrast}. The results indicate that as the contrast ratio $\mathtt{cr}$ increases, the set-up phase constitutes a smaller portion of the total computational time. This reduction allows the MMG method to demonstrate a significant advantage over AMG. In extreme cases where $\mathtt{cr}=6$, our method is $2.15$ times faster than AMG. Specifically, the efficiency of our preconditioner becomes increasingly apparent with higher contrast levels, leading to faster convergence and reduced overall computational time. Consequently, this trend highlights the robustness of our approach against the contrast ratio of the coefficient field.

\section{Conclusions}\label{sec:conclusions}
In this paper, we developed a robust parallel solver for heat transfer topology optimization using the PETSc numerical library. Our solver significantly enhances both the efficiency and robustness of the optimization process, particularly in scenarios characterized by high-contrast coefficients and high-resolution domains.

To ensure comprehensive performance, all components of the solver are effectively parallelized. Performance analysis demonstrates that our proposed solver not only exhibits excellent scalability, accommodating up to 4480 MPI processes, but also achieves up to $2\times$ speedup compared to the default algebraic multigrid method available in PETSc for high-contrast cases. Moreover, the proposed preconditioner is significantly robust against the contrast ratios in the coefficient field. 

Several approaches can be employed to further improve the performance of the proposed solver. One of the most promising avenue is to utilize GPUs to efficiently manage the computationally intensive tasks, such as matrix-vector multiplications within our solver \cite{ye2024fast}.




\section*{Acknowledgments}
EC's research is partially supported by the Hong Kong RGC General Research Fund (Project numbers: 14305222 and 14304021).

\input{main.bbl}\end{document}

%% file: output.pgf
\begingroup%
\makeatletter%
\begin{pgfpicture}%
\pgfpathrectangle{\pgfpointorigin}{\pgfqpoint{4.802013in}{2.952747in}}%
\pgfusepath{use as bounding box, clip}%
\begin{pgfscope}%
\pgfsetbuttcap%
\pgfsetmiterjoin%
\definecolor{currentfill}{rgb}{1.000000,1.000000,1.000000}%
\pgfsetfillcolor{currentfill}%
\pgfsetlinewidth{0.000000pt}%
\definecolor{currentstroke}{rgb}{1.000000,1.000000,1.000000}%
\pgfsetstrokecolor{currentstroke}%
\pgfsetdash{}{0pt}%
\pgfpathmoveto{\pgfqpoint{0.000000in}{0.000000in}}%
\pgfpathlineto{\pgfqpoint{4.802013in}{0.000000in}}%
\pgfpathlineto{\pgfqpoint{4.802013in}{2.952747in}}%
\pgfpathlineto{\pgfqpoint{0.000000in}{2.952747in}}%
\pgfpathlineto{\pgfqpoint{0.000000in}{0.000000in}}%
\pgfpathclose%
\pgfusepath{fill}%
\end{pgfscope}%
\begin{pgfscope}%
\pgfsetbuttcap%
\pgfsetmiterjoin%
\definecolor{currentfill}{rgb}{1.000000,1.000000,1.000000}%
\pgfsetfillcolor{currentfill}%
\pgfsetlinewidth{0.000000pt}%
\definecolor{currentstroke}{rgb}{0.000000,0.000000,0.000000}%
\pgfsetstrokecolor{currentstroke}%
\pgfsetstrokeopacity{0.000000}%
\pgfsetdash{}{0pt}%
\pgfpathmoveto{\pgfqpoint{0.656948in}{0.504167in}}%
\pgfpathlineto{\pgfqpoint{4.531948in}{0.504167in}}%
\pgfpathlineto{\pgfqpoint{4.531948in}{2.814167in}}%
\pgfpathlineto{\pgfqpoint{0.656948in}{2.814167in}}%
\pgfpathlineto{\pgfqpoint{0.656948in}{0.504167in}}%
\pgfpathclose%
\pgfusepath{fill}%
\end{pgfscope}%
\begin{pgfscope}%
\pgfpathrectangle{\pgfqpoint{0.656948in}{0.504167in}}{\pgfqpoint{3.875000in}{2.310000in}}%
\pgfusepath{clip}%
\pgfsetbuttcap%
\pgfsetmiterjoin%
\definecolor{currentfill}{rgb}{0.941176,0.501961,0.501961}%
\pgfsetfillcolor{currentfill}%
\pgfsetfillopacity{0.200000}%
\pgfsetlinewidth{0.240900pt}%
\definecolor{currentstroke}{rgb}{0.941176,0.501961,0.501961}%
\pgfsetstrokecolor{currentstroke}%
\pgfsetstrokeopacity{0.200000}%
\pgfsetdash{}{0pt}%
\pgfpathmoveto{\pgfqpoint{0.833084in}{0.504167in}}%
\pgfpathlineto{\pgfqpoint{2.846071in}{0.504167in}}%
\pgfpathlineto{\pgfqpoint{2.846071in}{2.814167in}}%
\pgfpathlineto{\pgfqpoint{0.833084in}{2.814167in}}%
\pgfpathlineto{\pgfqpoint{0.833084in}{0.504167in}}%
\pgfpathclose%
\pgfusepath{stroke,fill}%
\end{pgfscope}%
\begin{pgfscope}%
\pgfpathrectangle{\pgfqpoint{0.656948in}{0.504167in}}{\pgfqpoint{3.875000in}{2.310000in}}%
\pgfusepath{clip}%
\pgfsetbuttcap%
\pgfsetmiterjoin%
\definecolor{currentfill}{rgb}{0.564706,0.933333,0.564706}%
\pgfsetfillcolor{currentfill}%
\pgfsetfillopacity{0.200000}%
\pgfsetlinewidth{0.240900pt}%
\definecolor{currentstroke}{rgb}{0.564706,0.933333,0.564706}%
\pgfsetstrokecolor{currentstroke}%
\pgfsetstrokeopacity{0.200000}%
\pgfsetdash{}{0pt}%
\pgfpathmoveto{\pgfqpoint{2.846071in}{0.504167in}}%
\pgfpathlineto{\pgfqpoint{3.852565in}{0.504167in}}%
\pgfpathlineto{\pgfqpoint{3.852565in}{2.814167in}}%
\pgfpathlineto{\pgfqpoint{2.846071in}{2.814167in}}%
\pgfpathlineto{\pgfqpoint{2.846071in}{0.504167in}}%
\pgfpathclose%
\pgfusepath{stroke,fill}%
\end{pgfscope}%
\begin{pgfscope}%
\pgfpathrectangle{\pgfqpoint{0.656948in}{0.504167in}}{\pgfqpoint{3.875000in}{2.310000in}}%
\pgfusepath{clip}%
\pgfsetbuttcap%
\pgfsetmiterjoin%
\definecolor{currentfill}{rgb}{0.678431,0.847059,0.901961}%
\pgfsetfillcolor{currentfill}%
\pgfsetfillopacity{0.200000}%
\pgfsetlinewidth{0.240900pt}%
\definecolor{currentstroke}{rgb}{0.678431,0.847059,0.901961}%
\pgfsetstrokecolor{currentstroke}%
\pgfsetstrokeopacity{0.200000}%
\pgfsetdash{}{0pt}%
\pgfpathmoveto{\pgfqpoint{3.852565in}{0.504167in}}%
\pgfpathlineto{\pgfqpoint{4.355811in}{0.504167in}}%
\pgfpathlineto{\pgfqpoint{4.355811in}{2.814167in}}%
\pgfpathlineto{\pgfqpoint{3.852565in}{2.814167in}}%
\pgfpathlineto{\pgfqpoint{3.852565in}{0.504167in}}%
\pgfpathclose%
\pgfusepath{stroke,fill}%
\end{pgfscope}%
\begin{pgfscope}%
\pgfpathrectangle{\pgfqpoint{0.656948in}{0.504167in}}{\pgfqpoint{3.875000in}{2.310000in}}%
\pgfusepath{clip}%
\pgfsetbuttcap%
\pgfsetroundjoin%
\definecolor{currentfill}{rgb}{0.341176,0.317647,0.317647}%
\pgfsetfillcolor{currentfill}%
\pgfsetlinewidth{0.240900pt}%
\definecolor{currentstroke}{rgb}{0.341176,0.317647,0.317647}%
\pgfsetstrokecolor{currentstroke}%
\pgfsetdash{}{0pt}%
\pgfsys@defobject{currentmarker}{\pgfqpoint{-0.038036in}{-0.038036in}}{\pgfqpoint{0.038036in}{0.038036in}}{%
\pgfpathmoveto{\pgfqpoint{0.000000in}{-0.038036in}}%
\pgfpathcurveto{\pgfqpoint{0.010087in}{-0.038036in}}{\pgfqpoint{0.019763in}{-0.034029in}}{\pgfqpoint{0.026896in}{-0.026896in}}%
\pgfpathcurveto{\pgfqpoint{0.034029in}{-0.019763in}}{\pgfqpoint{0.038036in}{-0.010087in}}{\pgfqpoint{0.038036in}{0.000000in}}%
\pgfpathcurveto{\pgfqpoint{0.038036in}{0.010087in}}{\pgfqpoint{0.034029in}{0.019763in}}{\pgfqpoint{0.026896in}{0.026896in}}%
\pgfpathcurveto{\pgfqpoint{0.019763in}{0.034029in}}{\pgfqpoint{0.010087in}{0.038036in}}{\pgfqpoint{0.000000in}{0.038036in}}%
\pgfpathcurveto{\pgfqpoint{-0.010087in}{0.038036in}}{\pgfqpoint{-0.019763in}{0.034029in}}{\pgfqpoint{-0.026896in}{0.026896in}}%
\pgfpathcurveto{\pgfqpoint{-0.034029in}{0.019763in}}{\pgfqpoint{-0.038036in}{0.010087in}}{\pgfqpoint{-0.038036in}{0.000000in}}%
\pgfpathcurveto{\pgfqpoint{-0.038036in}{-0.010087in}}{\pgfqpoint{-0.034029in}{-0.019763in}}{\pgfqpoint{-0.026896in}{-0.026896in}}%
\pgfpathcurveto{\pgfqpoint{-0.019763in}{-0.034029in}}{\pgfqpoint{-0.010087in}{-0.038036in}}{\pgfqpoint{0.000000in}{-0.038036in}}%
\pgfpathlineto{\pgfqpoint{0.000000in}{-0.038036in}}%
\pgfpathclose%
\pgfusepath{stroke,fill}%
}%
\begin{pgfscope}%
\pgfsys@transformshift{2.846071in}{0.802364in}%
\pgfsys@useobject{currentmarker}{}%
\end{pgfscope}%
\begin{pgfscope}%
\pgfsys@transformshift{3.852565in}{0.704959in}%
\pgfsys@useobject{currentmarker}{}%
\end{pgfscope}%
\begin{pgfscope}%
\pgfsys@transformshift{4.355811in}{0.697387in}%
\pgfsys@useobject{currentmarker}{}%
\end{pgfscope}%
\end{pgfscope}%
\begin{pgfscope}%
\pgfpathrectangle{\pgfqpoint{0.656948in}{0.504167in}}{\pgfqpoint{3.875000in}{2.310000in}}%
\pgfusepath{clip}%
\pgfsetrectcap%
\pgfsetroundjoin%
\pgfsetlinewidth{0.803000pt}%
\definecolor{currentstroke}{rgb}{0.690196,0.690196,0.690196}%
\pgfsetstrokecolor{currentstroke}%
\pgfsetdash{}{0pt}%
\pgfpathmoveto{\pgfqpoint{0.833084in}{0.504167in}}%
\pgfpathlineto{\pgfqpoint{0.833084in}{2.814167in}}%
\pgfusepath{stroke}%
\end{pgfscope}%
\begin{pgfscope}%
\pgfsetbuttcap%
\pgfsetroundjoin%
\definecolor{currentfill}{rgb}{0.000000,0.000000,0.000000}%
\pgfsetfillcolor{currentfill}%
\pgfsetlinewidth{0.803000pt}%
\definecolor{currentstroke}{rgb}{0.000000,0.000000,0.000000}%
\pgfsetstrokecolor{currentstroke}%
\pgfsetdash{}{0pt}%
\pgfsys@defobject{currentmarker}{\pgfqpoint{0.000000in}{-0.048611in}}{\pgfqpoint{0.000000in}{0.000000in}}{%
\pgfpathmoveto{\pgfqpoint{0.000000in}{0.000000in}}%
\pgfpathlineto{\pgfqpoint{0.000000in}{-0.048611in}}%
\pgfusepath{stroke,fill}%
}%
\begin{pgfscope}%
\pgfsys@transformshift{0.833084in}{0.504167in}%
\pgfsys@useobject{currentmarker}{}%
\end{pgfscope}%
\end{pgfscope}%
\begin{pgfscope}%
\definecolor{textcolor}{rgb}{0.000000,0.000000,0.000000}%
\pgfsetstrokecolor{textcolor}%
\pgfsetfillcolor{textcolor}%
\pgftext[x=0.833084in,y=0.377778in,,top]{\color{textcolor}{\rmfamily\fontsize{8.000000}{9.600000}\selectfont\catcode`\^=\active\def^{\ifmmode\sp\else\^{}\fi}\catcode`\%=\active\def
\end{pgfscope}%
\begin{pgfscope}%
\pgfpathrectangle{\pgfqpoint{0.656948in}{0.504167in}}{\pgfqpoint{3.875000in}{2.310000in}}%
\pgfusepath{clip}%
\pgfsetrectcap%
\pgfsetroundjoin%
\pgfsetlinewidth{0.803000pt}%
\definecolor{currentstroke}{rgb}{0.690196,0.690196,0.690196}%
\pgfsetstrokecolor{currentstroke}%
\pgfsetdash{}{0pt}%
\pgfpathmoveto{\pgfqpoint{1.336331in}{0.504167in}}%
\pgfpathlineto{\pgfqpoint{1.336331in}{2.814167in}}%
\pgfusepath{stroke}%
\end{pgfscope}%
\begin{pgfscope}%
\pgfsetbuttcap%
\pgfsetroundjoin%
\definecolor{currentfill}{rgb}{0.000000,0.000000,0.000000}%
\pgfsetfillcolor{currentfill}%
\pgfsetlinewidth{0.803000pt}%
\definecolor{currentstroke}{rgb}{0.000000,0.000000,0.000000}%
\pgfsetstrokecolor{currentstroke}%
\pgfsetdash{}{0pt}%
\pgfsys@defobject{currentmarker}{\pgfqpoint{0.000000in}{-0.048611in}}{\pgfqpoint{0.000000in}{0.000000in}}{%
\pgfpathmoveto{\pgfqpoint{0.000000in}{0.000000in}}%
\pgfpathlineto{\pgfqpoint{0.000000in}{-0.048611in}}%
\pgfusepath{stroke,fill}%
}%
\begin{pgfscope}%
\pgfsys@transformshift{1.336331in}{0.504167in}%
\pgfsys@useobject{currentmarker}{}%
\end{pgfscope}%
\end{pgfscope}%
\begin{pgfscope}%
\definecolor{textcolor}{rgb}{0.000000,0.000000,0.000000}%
\pgfsetstrokecolor{textcolor}%
\pgfsetfillcolor{textcolor}%
\pgftext[x=1.336331in,y=0.377778in,,top]{\color{textcolor}{\rmfamily\fontsize{8.000000}{9.600000}\selectfont\catcode`\^=\active\def^{\ifmmode\sp\else\^{}\fi}\catcode`\%=\active\def
\end{pgfscope}%
\begin{pgfscope}%
\pgfpathrectangle{\pgfqpoint{0.656948in}{0.504167in}}{\pgfqpoint{3.875000in}{2.310000in}}%
\pgfusepath{clip}%
\pgfsetrectcap%
\pgfsetroundjoin%
\pgfsetlinewidth{0.803000pt}%
\definecolor{currentstroke}{rgb}{0.690196,0.690196,0.690196}%
\pgfsetstrokecolor{currentstroke}%
\pgfsetdash{}{0pt}%
\pgfpathmoveto{\pgfqpoint{1.839578in}{0.504167in}}%
\pgfpathlineto{\pgfqpoint{1.839578in}{2.814167in}}%
\pgfusepath{stroke}%
\end{pgfscope}%
\begin{pgfscope}%
\pgfsetbuttcap%
\pgfsetroundjoin%
\definecolor{currentfill}{rgb}{0.000000,0.000000,0.000000}%
\pgfsetfillcolor{currentfill}%
\pgfsetlinewidth{0.803000pt}%
\definecolor{currentstroke}{rgb}{0.000000,0.000000,0.000000}%
\pgfsetstrokecolor{currentstroke}%
\pgfsetdash{}{0pt}%
\pgfsys@defobject{currentmarker}{\pgfqpoint{0.000000in}{-0.048611in}}{\pgfqpoint{0.000000in}{0.000000in}}{%
\pgfpathmoveto{\pgfqpoint{0.000000in}{0.000000in}}%
\pgfpathlineto{\pgfqpoint{0.000000in}{-0.048611in}}%
\pgfusepath{stroke,fill}%
}%
\begin{pgfscope}%
\pgfsys@transformshift{1.839578in}{0.504167in}%
\pgfsys@useobject{currentmarker}{}%
\end{pgfscope}%
\end{pgfscope}%
\begin{pgfscope}%
\definecolor{textcolor}{rgb}{0.000000,0.000000,0.000000}%
\pgfsetstrokecolor{textcolor}%
\pgfsetfillcolor{textcolor}%
\pgftext[x=1.839578in,y=0.377778in,,top]{\color{textcolor}{\rmfamily\fontsize{8.000000}{9.600000}\selectfont\catcode`\^=\active\def^{\ifmmode\sp\else\^{}\fi}\catcode`\%=\active\def
\end{pgfscope}%
\begin{pgfscope}%
\pgfpathrectangle{\pgfqpoint{0.656948in}{0.504167in}}{\pgfqpoint{3.875000in}{2.310000in}}%
\pgfusepath{clip}%
\pgfsetrectcap%
\pgfsetroundjoin%
\pgfsetlinewidth{0.803000pt}%
\definecolor{currentstroke}{rgb}{0.690196,0.690196,0.690196}%
\pgfsetstrokecolor{currentstroke}%
\pgfsetdash{}{0pt}%
\pgfpathmoveto{\pgfqpoint{2.342824in}{0.504167in}}%
\pgfpathlineto{\pgfqpoint{2.342824in}{2.814167in}}%
\pgfusepath{stroke}%
\end{pgfscope}%
\begin{pgfscope}%
\pgfsetbuttcap%
\pgfsetroundjoin%
\definecolor{currentfill}{rgb}{0.000000,0.000000,0.000000}%
\pgfsetfillcolor{currentfill}%
\pgfsetlinewidth{0.803000pt}%
\definecolor{currentstroke}{rgb}{0.000000,0.000000,0.000000}%
\pgfsetstrokecolor{currentstroke}%
\pgfsetdash{}{0pt}%
\pgfsys@defobject{currentmarker}{\pgfqpoint{0.000000in}{-0.048611in}}{\pgfqpoint{0.000000in}{0.000000in}}{%
\pgfpathmoveto{\pgfqpoint{0.000000in}{0.000000in}}%
\pgfpathlineto{\pgfqpoint{0.000000in}{-0.048611in}}%
\pgfusepath{stroke,fill}%
}%
\begin{pgfscope}%
\pgfsys@transformshift{2.342824in}{0.504167in}%
\pgfsys@useobject{currentmarker}{}%
\end{pgfscope}%
\end{pgfscope}%
\begin{pgfscope}%
\definecolor{textcolor}{rgb}{0.000000,0.000000,0.000000}%
\pgfsetstrokecolor{textcolor}%
\pgfsetfillcolor{textcolor}%
\pgftext[x=2.342824in,y=0.377778in,,top]{\color{textcolor}{\rmfamily\fontsize{8.000000}{9.600000}\selectfont\catcode`\^=\active\def^{\ifmmode\sp\else\^{}\fi}\catcode`\%=\active\def
\end{pgfscope}%
\begin{pgfscope}%
\pgfpathrectangle{\pgfqpoint{0.656948in}{0.504167in}}{\pgfqpoint{3.875000in}{2.310000in}}%
\pgfusepath{clip}%
\pgfsetrectcap%
\pgfsetroundjoin%
\pgfsetlinewidth{0.803000pt}%
\definecolor{currentstroke}{rgb}{0.690196,0.690196,0.690196}%
\pgfsetstrokecolor{currentstroke}%
\pgfsetdash{}{0pt}%
\pgfpathmoveto{\pgfqpoint{2.846071in}{0.504167in}}%
\pgfpathlineto{\pgfqpoint{2.846071in}{2.814167in}}%
\pgfusepath{stroke}%
\end{pgfscope}%
\begin{pgfscope}%
\pgfsetbuttcap%
\pgfsetroundjoin%
\definecolor{currentfill}{rgb}{0.000000,0.000000,0.000000}%
\pgfsetfillcolor{currentfill}%
\pgfsetlinewidth{0.803000pt}%
\definecolor{currentstroke}{rgb}{0.000000,0.000000,0.000000}%
\pgfsetstrokecolor{currentstroke}%
\pgfsetdash{}{0pt}%
\pgfsys@defobject{currentmarker}{\pgfqpoint{0.000000in}{-0.048611in}}{\pgfqpoint{0.000000in}{0.000000in}}{%
\pgfpathmoveto{\pgfqpoint{0.000000in}{0.000000in}}%
\pgfpathlineto{\pgfqpoint{0.000000in}{-0.048611in}}%
\pgfusepath{stroke,fill}%
}%
\begin{pgfscope}%
\pgfsys@transformshift{2.846071in}{0.504167in}%
\pgfsys@useobject{currentmarker}{}%
\end{pgfscope}%
\end{pgfscope}%
\begin{pgfscope}%
\definecolor{textcolor}{rgb}{0.000000,0.000000,0.000000}%
\pgfsetstrokecolor{textcolor}%
\pgfsetfillcolor{textcolor}%
\pgftext[x=2.846071in,y=0.377778in,,top]{\color{textcolor}{\rmfamily\fontsize{8.000000}{9.600000}\selectfont\catcode`\^=\active\def^{\ifmmode\sp\else\^{}\fi}\catcode`\%=\active\def
\end{pgfscope}%
\begin{pgfscope}%
\pgfpathrectangle{\pgfqpoint{0.656948in}{0.504167in}}{\pgfqpoint{3.875000in}{2.310000in}}%
\pgfusepath{clip}%
\pgfsetrectcap%
\pgfsetroundjoin%
\pgfsetlinewidth{0.803000pt}%
\definecolor{currentstroke}{rgb}{0.690196,0.690196,0.690196}%
\pgfsetstrokecolor{currentstroke}%
\pgfsetdash{}{0pt}%
\pgfpathmoveto{\pgfqpoint{3.349318in}{0.504167in}}%
\pgfpathlineto{\pgfqpoint{3.349318in}{2.814167in}}%
\pgfusepath{stroke}%
\end{pgfscope}%
\begin{pgfscope}%
\pgfsetbuttcap%
\pgfsetroundjoin%
\definecolor{currentfill}{rgb}{0.000000,0.000000,0.000000}%
\pgfsetfillcolor{currentfill}%
\pgfsetlinewidth{0.803000pt}%
\definecolor{currentstroke}{rgb}{0.000000,0.000000,0.000000}%
\pgfsetstrokecolor{currentstroke}%
\pgfsetdash{}{0pt}%
\pgfsys@defobject{currentmarker}{\pgfqpoint{0.000000in}{-0.048611in}}{\pgfqpoint{0.000000in}{0.000000in}}{%
\pgfpathmoveto{\pgfqpoint{0.000000in}{0.000000in}}%
\pgfpathlineto{\pgfqpoint{0.000000in}{-0.048611in}}%
\pgfusepath{stroke,fill}%
}%
\begin{pgfscope}%
\pgfsys@transformshift{3.349318in}{0.504167in}%
\pgfsys@useobject{currentmarker}{}%
\end{pgfscope}%
\end{pgfscope}%
\begin{pgfscope}%
\definecolor{textcolor}{rgb}{0.000000,0.000000,0.000000}%
\pgfsetstrokecolor{textcolor}%
\pgfsetfillcolor{textcolor}%
\pgftext[x=3.349318in,y=0.377778in,,top]{\color{textcolor}{\rmfamily\fontsize{8.000000}{9.600000}\selectfont\catcode`\^=\active\def^{\ifmmode\sp\else\^{}\fi}\catcode`\%=\active\def
\end{pgfscope}%
\begin{pgfscope}%
\pgfpathrectangle{\pgfqpoint{0.656948in}{0.504167in}}{\pgfqpoint{3.875000in}{2.310000in}}%
\pgfusepath{clip}%
\pgfsetrectcap%
\pgfsetroundjoin%
\pgfsetlinewidth{0.803000pt}%
\definecolor{currentstroke}{rgb}{0.690196,0.690196,0.690196}%
\pgfsetstrokecolor{currentstroke}%
\pgfsetdash{}{0pt}%
\pgfpathmoveto{\pgfqpoint{3.852565in}{0.504167in}}%
\pgfpathlineto{\pgfqpoint{3.852565in}{2.814167in}}%
\pgfusepath{stroke}%
\end{pgfscope}%
\begin{pgfscope}%
\pgfsetbuttcap%
\pgfsetroundjoin%
\definecolor{currentfill}{rgb}{0.000000,0.000000,0.000000}%
\pgfsetfillcolor{currentfill}%
\pgfsetlinewidth{0.803000pt}%
\definecolor{currentstroke}{rgb}{0.000000,0.000000,0.000000}%
\pgfsetstrokecolor{currentstroke}%
\pgfsetdash{}{0pt}%
\pgfsys@defobject{currentmarker}{\pgfqpoint{0.000000in}{-0.048611in}}{\pgfqpoint{0.000000in}{0.000000in}}{%
\pgfpathmoveto{\pgfqpoint{0.000000in}{0.000000in}}%
\pgfpathlineto{\pgfqpoint{0.000000in}{-0.048611in}}%
\pgfusepath{stroke,fill}%
}%
\begin{pgfscope}%
\pgfsys@transformshift{3.852565in}{0.504167in}%
\pgfsys@useobject{currentmarker}{}%
\end{pgfscope}%
\end{pgfscope}%
\begin{pgfscope}%
\definecolor{textcolor}{rgb}{0.000000,0.000000,0.000000}%
\pgfsetstrokecolor{textcolor}%
\pgfsetfillcolor{textcolor}%
\pgftext[x=3.852565in,y=0.377778in,,top]{\color{textcolor}{\rmfamily\fontsize{8.000000}{9.600000}\selectfont\catcode`\^=\active\def^{\ifmmode\sp\else\^{}\fi}\catcode`\%=\active\def
\end{pgfscope}%
\begin{pgfscope}%
\pgfpathrectangle{\pgfqpoint{0.656948in}{0.504167in}}{\pgfqpoint{3.875000in}{2.310000in}}%
\pgfusepath{clip}%
\pgfsetrectcap%
\pgfsetroundjoin%
\pgfsetlinewidth{0.803000pt}%
\definecolor{currentstroke}{rgb}{0.690196,0.690196,0.690196}%
\pgfsetstrokecolor{currentstroke}%
\pgfsetdash{}{0pt}%
\pgfpathmoveto{\pgfqpoint{4.355811in}{0.504167in}}%
\pgfpathlineto{\pgfqpoint{4.355811in}{2.814167in}}%
\pgfusepath{stroke}%
\end{pgfscope}%
\begin{pgfscope}%
\pgfsetbuttcap%
\pgfsetroundjoin%
\definecolor{currentfill}{rgb}{0.000000,0.000000,0.000000}%
\pgfsetfillcolor{currentfill}%
\pgfsetlinewidth{0.803000pt}%
\definecolor{currentstroke}{rgb}{0.000000,0.000000,0.000000}%
\pgfsetstrokecolor{currentstroke}%
\pgfsetdash{}{0pt}%
\pgfsys@defobject{currentmarker}{\pgfqpoint{0.000000in}{-0.048611in}}{\pgfqpoint{0.000000in}{0.000000in}}{%
\pgfpathmoveto{\pgfqpoint{0.000000in}{0.000000in}}%
\pgfpathlineto{\pgfqpoint{0.000000in}{-0.048611in}}%
\pgfusepath{stroke,fill}%
}%
\begin{pgfscope}%
\pgfsys@transformshift{4.355811in}{0.504167in}%
\pgfsys@useobject{currentmarker}{}%
\end{pgfscope}%
\end{pgfscope}%
\begin{pgfscope}%
\definecolor{textcolor}{rgb}{0.000000,0.000000,0.000000}%
\pgfsetstrokecolor{textcolor}%
\pgfsetfillcolor{textcolor}%
\pgftext[x=4.355811in,y=0.377778in,,top]{\color{textcolor}{\rmfamily\fontsize{8.000000}{9.600000}\selectfont\catcode`\^=\active\def^{\ifmmode\sp\else\^{}\fi}\catcode`\%=\active\def
\end{pgfscope}%
\begin{pgfscope}%
\definecolor{textcolor}{rgb}{0.000000,0.000000,0.000000}%
\pgfsetstrokecolor{textcolor}%
\pgfsetfillcolor{textcolor}%
\pgftext[x=2.594448in,y=0.223457in,,top]{\color{textcolor}{\rmfamily\fontsize{10.000000}{12.000000}\selectfont\catcode`\^=\active\def^{\ifmmode\sp\else\^{}\fi}\catcode`\%=\active\def
\end{pgfscope}%
\begin{pgfscope}%
\pgfpathrectangle{\pgfqpoint{0.656948in}{0.504167in}}{\pgfqpoint{3.875000in}{2.310000in}}%
\pgfusepath{clip}%
\pgfsetrectcap%
\pgfsetroundjoin%
\pgfsetlinewidth{0.803000pt}%
\definecolor{currentstroke}{rgb}{0.690196,0.690196,0.690196}%
\pgfsetstrokecolor{currentstroke}%
\pgfsetdash{}{0pt}%
\pgfpathmoveto{\pgfqpoint{0.656948in}{0.504167in}}%
\pgfpathlineto{\pgfqpoint{4.531948in}{0.504167in}}%
\pgfusepath{stroke}%
\end{pgfscope}%
\begin{pgfscope}%
\pgfsetbuttcap%
\pgfsetroundjoin%
\definecolor{currentfill}{rgb}{0.000000,0.000000,0.000000}%
\pgfsetfillcolor{currentfill}%
\pgfsetlinewidth{0.803000pt}%
\definecolor{currentstroke}{rgb}{0.000000,0.000000,0.000000}%
\pgfsetstrokecolor{currentstroke}%
\pgfsetdash{}{0pt}%
\pgfsys@defobject{currentmarker}{\pgfqpoint{-0.048611in}{0.000000in}}{\pgfqpoint{-0.000000in}{0.000000in}}{%
\pgfpathmoveto{\pgfqpoint{-0.000000in}{0.000000in}}%
\pgfpathlineto{\pgfqpoint{-0.048611in}{0.000000in}}%
\pgfusepath{stroke,fill}%
}%
\begin{pgfscope}%
\pgfsys@transformshift{0.656948in}{0.504167in}%
\pgfsys@useobject{currentmarker}{}%
\end{pgfscope}%
\end{pgfscope}%
\begin{pgfscope}%
\definecolor{textcolor}{rgb}{0.000000,0.000000,0.000000}%
\pgfsetstrokecolor{textcolor}%
\pgfsetfillcolor{textcolor}%
\pgftext[x=0.471530in, y=0.465586in, left, base]{\color{textcolor}{\rmfamily\fontsize{8.000000}{9.600000}\selectfont\catcode`\^=\active\def^{\ifmmode\sp\else\^{}\fi}\catcode`\%=\active\def
\end{pgfscope}%
\begin{pgfscope}%
\pgfpathrectangle{\pgfqpoint{0.656948in}{0.504167in}}{\pgfqpoint{3.875000in}{2.310000in}}%
\pgfusepath{clip}%
\pgfsetrectcap%
\pgfsetroundjoin%
\pgfsetlinewidth{0.803000pt}%
\definecolor{currentstroke}{rgb}{0.690196,0.690196,0.690196}%
\pgfsetstrokecolor{currentstroke}%
\pgfsetdash{}{0pt}%
\pgfpathmoveto{\pgfqpoint{0.656948in}{0.792917in}}%
\pgfpathlineto{\pgfqpoint{4.531948in}{0.792917in}}%
\pgfusepath{stroke}%
\end{pgfscope}%
\begin{pgfscope}%
\pgfsetbuttcap%
\pgfsetroundjoin%
\definecolor{currentfill}{rgb}{0.000000,0.000000,0.000000}%
\pgfsetfillcolor{currentfill}%
\pgfsetlinewidth{0.803000pt}%
\definecolor{currentstroke}{rgb}{0.000000,0.000000,0.000000}%
\pgfsetstrokecolor{currentstroke}%
\pgfsetdash{}{0pt}%
\pgfsys@defobject{currentmarker}{\pgfqpoint{-0.048611in}{0.000000in}}{\pgfqpoint{-0.000000in}{0.000000in}}{%
\pgfpathmoveto{\pgfqpoint{-0.000000in}{0.000000in}}%
\pgfpathlineto{\pgfqpoint{-0.048611in}{0.000000in}}%
\pgfusepath{stroke,fill}%
}%
\begin{pgfscope}%
\pgfsys@transformshift{0.656948in}{0.792917in}%
\pgfsys@useobject{currentmarker}{}%
\end{pgfscope}%
\end{pgfscope}%
\begin{pgfscope}%
\definecolor{textcolor}{rgb}{0.000000,0.000000,0.000000}%
\pgfsetstrokecolor{textcolor}%
\pgfsetfillcolor{textcolor}%
\pgftext[x=0.353473in, y=0.754336in, left, base]{\color{textcolor}{\rmfamily\fontsize{8.000000}{9.600000}\selectfont\catcode`\^=\active\def^{\ifmmode\sp\else\^{}\fi}\catcode`\%=\active\def
\end{pgfscope}%
\begin{pgfscope}%
\pgfpathrectangle{\pgfqpoint{0.656948in}{0.504167in}}{\pgfqpoint{3.875000in}{2.310000in}}%
\pgfusepath{clip}%
\pgfsetrectcap%
\pgfsetroundjoin%
\pgfsetlinewidth{0.803000pt}%
\definecolor{currentstroke}{rgb}{0.690196,0.690196,0.690196}%
\pgfsetstrokecolor{currentstroke}%
\pgfsetdash{}{0pt}%
\pgfpathmoveto{\pgfqpoint{0.656948in}{1.081667in}}%
\pgfpathlineto{\pgfqpoint{4.531948in}{1.081667in}}%
\pgfusepath{stroke}%
\end{pgfscope}%
\begin{pgfscope}%
\pgfsetbuttcap%
\pgfsetroundjoin%
\definecolor{currentfill}{rgb}{0.000000,0.000000,0.000000}%
\pgfsetfillcolor{currentfill}%
\pgfsetlinewidth{0.803000pt}%
\definecolor{currentstroke}{rgb}{0.000000,0.000000,0.000000}%
\pgfsetstrokecolor{currentstroke}%
\pgfsetdash{}{0pt}%
\pgfsys@defobject{currentmarker}{\pgfqpoint{-0.048611in}{0.000000in}}{\pgfqpoint{-0.000000in}{0.000000in}}{%
\pgfpathmoveto{\pgfqpoint{-0.000000in}{0.000000in}}%
\pgfpathlineto{\pgfqpoint{-0.048611in}{0.000000in}}%
\pgfusepath{stroke,fill}%
}%
\begin{pgfscope}%
\pgfsys@transformshift{0.656948in}{1.081667in}%
\pgfsys@useobject{currentmarker}{}%
\end{pgfscope}%
\end{pgfscope}%
\begin{pgfscope}%
\definecolor{textcolor}{rgb}{0.000000,0.000000,0.000000}%
\pgfsetstrokecolor{textcolor}%
\pgfsetfillcolor{textcolor}%
\pgftext[x=0.353473in, y=1.043086in, left, base]{\color{textcolor}{\rmfamily\fontsize{8.000000}{9.600000}\selectfont\catcode`\^=\active\def^{\ifmmode\sp\else\^{}\fi}\catcode`\%=\active\def
\end{pgfscope}%
\begin{pgfscope}%
\pgfpathrectangle{\pgfqpoint{0.656948in}{0.504167in}}{\pgfqpoint{3.875000in}{2.310000in}}%
\pgfusepath{clip}%
\pgfsetrectcap%
\pgfsetroundjoin%
\pgfsetlinewidth{0.803000pt}%
\definecolor{currentstroke}{rgb}{0.690196,0.690196,0.690196}%
\pgfsetstrokecolor{currentstroke}%
\pgfsetdash{}{0pt}%
\pgfpathmoveto{\pgfqpoint{0.656948in}{1.370417in}}%
\pgfpathlineto{\pgfqpoint{4.531948in}{1.370417in}}%
\pgfusepath{stroke}%
\end{pgfscope}%
\begin{pgfscope}%
\pgfsetbuttcap%
\pgfsetroundjoin%
\definecolor{currentfill}{rgb}{0.000000,0.000000,0.000000}%
\pgfsetfillcolor{currentfill}%
\pgfsetlinewidth{0.803000pt}%
\definecolor{currentstroke}{rgb}{0.000000,0.000000,0.000000}%
\pgfsetstrokecolor{currentstroke}%
\pgfsetdash{}{0pt}%
\pgfsys@defobject{currentmarker}{\pgfqpoint{-0.048611in}{0.000000in}}{\pgfqpoint{-0.000000in}{0.000000in}}{%
\pgfpathmoveto{\pgfqpoint{-0.000000in}{0.000000in}}%
\pgfpathlineto{\pgfqpoint{-0.048611in}{0.000000in}}%
\pgfusepath{stroke,fill}%
}%
\begin{pgfscope}%
\pgfsys@transformshift{0.656948in}{1.370417in}%
\pgfsys@useobject{currentmarker}{}%
\end{pgfscope}%
\end{pgfscope}%
\begin{pgfscope}%
\definecolor{textcolor}{rgb}{0.000000,0.000000,0.000000}%
\pgfsetstrokecolor{textcolor}%
\pgfsetfillcolor{textcolor}%
\pgftext[x=0.353473in, y=1.331836in, left, base]{\color{textcolor}{\rmfamily\fontsize{8.000000}{9.600000}\selectfont\catcode`\^=\active\def^{\ifmmode\sp\else\^{}\fi}\catcode`\%=\active\def
\end{pgfscope}%
\begin{pgfscope}%
\pgfpathrectangle{\pgfqpoint{0.656948in}{0.504167in}}{\pgfqpoint{3.875000in}{2.310000in}}%
\pgfusepath{clip}%
\pgfsetrectcap%
\pgfsetroundjoin%
\pgfsetlinewidth{0.803000pt}%
\definecolor{currentstroke}{rgb}{0.690196,0.690196,0.690196}%
\pgfsetstrokecolor{currentstroke}%
\pgfsetdash{}{0pt}%
\pgfpathmoveto{\pgfqpoint{0.656948in}{1.659167in}}%
\pgfpathlineto{\pgfqpoint{4.531948in}{1.659167in}}%
\pgfusepath{stroke}%
\end{pgfscope}%
\begin{pgfscope}%
\pgfsetbuttcap%
\pgfsetroundjoin%
\definecolor{currentfill}{rgb}{0.000000,0.000000,0.000000}%
\pgfsetfillcolor{currentfill}%
\pgfsetlinewidth{0.803000pt}%
\definecolor{currentstroke}{rgb}{0.000000,0.000000,0.000000}%
\pgfsetstrokecolor{currentstroke}%
\pgfsetdash{}{0pt}%
\pgfsys@defobject{currentmarker}{\pgfqpoint{-0.048611in}{0.000000in}}{\pgfqpoint{-0.000000in}{0.000000in}}{%
\pgfpathmoveto{\pgfqpoint{-0.000000in}{0.000000in}}%
\pgfpathlineto{\pgfqpoint{-0.048611in}{0.000000in}}%
\pgfusepath{stroke,fill}%
}%
\begin{pgfscope}%
\pgfsys@transformshift{0.656948in}{1.659167in}%
\pgfsys@useobject{currentmarker}{}%
\end{pgfscope}%
\end{pgfscope}%
\begin{pgfscope}%
\definecolor{textcolor}{rgb}{0.000000,0.000000,0.000000}%
\pgfsetstrokecolor{textcolor}%
\pgfsetfillcolor{textcolor}%
\pgftext[x=0.294444in, y=1.620586in, left, base]{\color{textcolor}{\rmfamily\fontsize{8.000000}{9.600000}\selectfont\catcode`\^=\active\def^{\ifmmode\sp\else\^{}\fi}\catcode`\%=\active\def
\end{pgfscope}%
\begin{pgfscope}%
\pgfpathrectangle{\pgfqpoint{0.656948in}{0.504167in}}{\pgfqpoint{3.875000in}{2.310000in}}%
\pgfusepath{clip}%
\pgfsetrectcap%
\pgfsetroundjoin%
\pgfsetlinewidth{0.803000pt}%
\definecolor{currentstroke}{rgb}{0.690196,0.690196,0.690196}%
\pgfsetstrokecolor{currentstroke}%
\pgfsetdash{}{0pt}%
\pgfpathmoveto{\pgfqpoint{0.656948in}{1.947917in}}%
\pgfpathlineto{\pgfqpoint{4.531948in}{1.947917in}}%
\pgfusepath{stroke}%
\end{pgfscope}%
\begin{pgfscope}%
\pgfsetbuttcap%
\pgfsetroundjoin%
\definecolor{currentfill}{rgb}{0.000000,0.000000,0.000000}%
\pgfsetfillcolor{currentfill}%
\pgfsetlinewidth{0.803000pt}%
\definecolor{currentstroke}{rgb}{0.000000,0.000000,0.000000}%
\pgfsetstrokecolor{currentstroke}%
\pgfsetdash{}{0pt}%
\pgfsys@defobject{currentmarker}{\pgfqpoint{-0.048611in}{0.000000in}}{\pgfqpoint{-0.000000in}{0.000000in}}{%
\pgfpathmoveto{\pgfqpoint{-0.000000in}{0.000000in}}%
\pgfpathlineto{\pgfqpoint{-0.048611in}{0.000000in}}%
\pgfusepath{stroke,fill}%
}%
\begin{pgfscope}%
\pgfsys@transformshift{0.656948in}{1.947917in}%
\pgfsys@useobject{currentmarker}{}%
\end{pgfscope}%
\end{pgfscope}%
\begin{pgfscope}%
\definecolor{textcolor}{rgb}{0.000000,0.000000,0.000000}%
\pgfsetstrokecolor{textcolor}%
\pgfsetfillcolor{textcolor}%
\pgftext[x=0.294444in, y=1.909336in, left, base]{\color{textcolor}{\rmfamily\fontsize{8.000000}{9.600000}\selectfont\catcode`\^=\active\def^{\ifmmode\sp\else\^{}\fi}\catcode`\%=\active\def
\end{pgfscope}%
\begin{pgfscope}%
\pgfpathrectangle{\pgfqpoint{0.656948in}{0.504167in}}{\pgfqpoint{3.875000in}{2.310000in}}%
\pgfusepath{clip}%
\pgfsetrectcap%
\pgfsetroundjoin%
\pgfsetlinewidth{0.803000pt}%
\definecolor{currentstroke}{rgb}{0.690196,0.690196,0.690196}%
\pgfsetstrokecolor{currentstroke}%
\pgfsetdash{}{0pt}%
\pgfpathmoveto{\pgfqpoint{0.656948in}{2.236667in}}%
\pgfpathlineto{\pgfqpoint{4.531948in}{2.236667in}}%
\pgfusepath{stroke}%
\end{pgfscope}%
\begin{pgfscope}%
\pgfsetbuttcap%
\pgfsetroundjoin%
\definecolor{currentfill}{rgb}{0.000000,0.000000,0.000000}%
\pgfsetfillcolor{currentfill}%
\pgfsetlinewidth{0.803000pt}%
\definecolor{currentstroke}{rgb}{0.000000,0.000000,0.000000}%
\pgfsetstrokecolor{currentstroke}%
\pgfsetdash{}{0pt}%
\pgfsys@defobject{currentmarker}{\pgfqpoint{-0.048611in}{0.000000in}}{\pgfqpoint{-0.000000in}{0.000000in}}{%
\pgfpathmoveto{\pgfqpoint{-0.000000in}{0.000000in}}%
\pgfpathlineto{\pgfqpoint{-0.048611in}{0.000000in}}%
\pgfusepath{stroke,fill}%
}%
\begin{pgfscope}%
\pgfsys@transformshift{0.656948in}{2.236667in}%
\pgfsys@useobject{currentmarker}{}%
\end{pgfscope}%
\end{pgfscope}%
\begin{pgfscope}%
\definecolor{textcolor}{rgb}{0.000000,0.000000,0.000000}%
\pgfsetstrokecolor{textcolor}%
\pgfsetfillcolor{textcolor}%
\pgftext[x=0.294444in, y=2.198086in, left, base]{\color{textcolor}{\rmfamily\fontsize{8.000000}{9.600000}\selectfont\catcode`\^=\active\def^{\ifmmode\sp\else\^{}\fi}\catcode`\%=\active\def
\end{pgfscope}%
\begin{pgfscope}%
\pgfpathrectangle{\pgfqpoint{0.656948in}{0.504167in}}{\pgfqpoint{3.875000in}{2.310000in}}%
\pgfusepath{clip}%
\pgfsetrectcap%
\pgfsetroundjoin%
\pgfsetlinewidth{0.803000pt}%
\definecolor{currentstroke}{rgb}{0.690196,0.690196,0.690196}%
\pgfsetstrokecolor{currentstroke}%
\pgfsetdash{}{0pt}%
\pgfpathmoveto{\pgfqpoint{0.656948in}{2.525417in}}%
\pgfpathlineto{\pgfqpoint{4.531948in}{2.525417in}}%
\pgfusepath{stroke}%
\end{pgfscope}%
\begin{pgfscope}%
\pgfsetbuttcap%
\pgfsetroundjoin%
\definecolor{currentfill}{rgb}{0.000000,0.000000,0.000000}%
\pgfsetfillcolor{currentfill}%
\pgfsetlinewidth{0.803000pt}%
\definecolor{currentstroke}{rgb}{0.000000,0.000000,0.000000}%
\pgfsetstrokecolor{currentstroke}%
\pgfsetdash{}{0pt}%
\pgfsys@defobject{currentmarker}{\pgfqpoint{-0.048611in}{0.000000in}}{\pgfqpoint{-0.000000in}{0.000000in}}{%
\pgfpathmoveto{\pgfqpoint{-0.000000in}{0.000000in}}%
\pgfpathlineto{\pgfqpoint{-0.048611in}{0.000000in}}%
\pgfusepath{stroke,fill}%
}%
\begin{pgfscope}%
\pgfsys@transformshift{0.656948in}{2.525417in}%
\pgfsys@useobject{currentmarker}{}%
\end{pgfscope}%
\end{pgfscope}%
\begin{pgfscope}%
\definecolor{textcolor}{rgb}{0.000000,0.000000,0.000000}%
\pgfsetstrokecolor{textcolor}%
\pgfsetfillcolor{textcolor}%
\pgftext[x=0.294444in, y=2.486836in, left, base]{\color{textcolor}{\rmfamily\fontsize{8.000000}{9.600000}\selectfont\catcode`\^=\active\def^{\ifmmode\sp\else\^{}\fi}\catcode`\%=\active\def
\end{pgfscope}%
\begin{pgfscope}%
\pgfpathrectangle{\pgfqpoint{0.656948in}{0.504167in}}{\pgfqpoint{3.875000in}{2.310000in}}%
\pgfusepath{clip}%
\pgfsetrectcap%
\pgfsetroundjoin%
\pgfsetlinewidth{0.803000pt}%
\definecolor{currentstroke}{rgb}{0.690196,0.690196,0.690196}%
\pgfsetstrokecolor{currentstroke}%
\pgfsetdash{}{0pt}%
\pgfpathmoveto{\pgfqpoint{0.656948in}{2.814167in}}%
\pgfpathlineto{\pgfqpoint{4.531948in}{2.814167in}}%
\pgfusepath{stroke}%
\end{pgfscope}%
\begin{pgfscope}%
\pgfsetbuttcap%
\pgfsetroundjoin%
\definecolor{currentfill}{rgb}{0.000000,0.000000,0.000000}%
\pgfsetfillcolor{currentfill}%
\pgfsetlinewidth{0.803000pt}%
\definecolor{currentstroke}{rgb}{0.000000,0.000000,0.000000}%
\pgfsetstrokecolor{currentstroke}%
\pgfsetdash{}{0pt}%
\pgfsys@defobject{currentmarker}{\pgfqpoint{-0.048611in}{0.000000in}}{\pgfqpoint{-0.000000in}{0.000000in}}{%
\pgfpathmoveto{\pgfqpoint{-0.000000in}{0.000000in}}%
\pgfpathlineto{\pgfqpoint{-0.048611in}{0.000000in}}%
\pgfusepath{stroke,fill}%
}%
\begin{pgfscope}%
\pgfsys@transformshift{0.656948in}{2.814167in}%
\pgfsys@useobject{currentmarker}{}%
\end{pgfscope}%
\end{pgfscope}%
\begin{pgfscope}%
\definecolor{textcolor}{rgb}{0.000000,0.000000,0.000000}%
\pgfsetstrokecolor{textcolor}%
\pgfsetfillcolor{textcolor}%
\pgftext[x=0.294444in, y=2.775586in, left, base]{\color{textcolor}{\rmfamily\fontsize{8.000000}{9.600000}\selectfont\catcode`\^=\active\def^{\ifmmode\sp\else\^{}\fi}\catcode`\%=\active\def
\end{pgfscope}%
\begin{pgfscope}%
\definecolor{textcolor}{rgb}{0.000000,0.000000,0.000000}%
\pgfsetstrokecolor{textcolor}%
\pgfsetfillcolor{textcolor}%
\pgftext[x=0.238889in,y=1.659167in,,bottom,rotate=90.000000]{\color{textcolor}{\rmfamily\fontsize{10.000000}{12.000000}\selectfont\catcode`\^=\active\def^{\ifmmode\sp\else\^{}\fi}\catcode`\%=\active\def
\end{pgfscope}%
\begin{pgfscope}%
\pgfpathrectangle{\pgfqpoint{0.656948in}{0.504167in}}{\pgfqpoint{3.875000in}{2.310000in}}%
\pgfusepath{clip}%
\pgfsetrectcap%
\pgfsetroundjoin%
\pgfsetlinewidth{1.405250pt}%
\definecolor{currentstroke}{rgb}{0.341176,0.317647,0.317647}%
\pgfsetstrokecolor{currentstroke}%
\pgfsetdash{}{0pt}%
\pgfpathmoveto{\pgfqpoint{0.969236in}{2.824167in}}%
\pgfpathlineto{\pgfqpoint{0.984058in}{2.550304in}}%
\pgfpathlineto{\pgfqpoint{1.034383in}{2.131383in}}%
\pgfpathlineto{\pgfqpoint{1.084708in}{1.874904in}}%
\pgfpathlineto{\pgfqpoint{1.135032in}{1.686530in}}%
\pgfpathlineto{\pgfqpoint{1.185357in}{1.537875in}}%
\pgfpathlineto{\pgfqpoint{1.235682in}{1.419416in}}%
\pgfpathlineto{\pgfqpoint{1.286006in}{1.325433in}}%
\pgfpathlineto{\pgfqpoint{1.336331in}{1.251054in}}%
\pgfpathlineto{\pgfqpoint{1.386656in}{1.191607in}}%
\pgfpathlineto{\pgfqpoint{1.436980in}{1.143306in}}%
\pgfpathlineto{\pgfqpoint{1.487305in}{1.102878in}}%
\pgfpathlineto{\pgfqpoint{1.537630in}{1.068135in}}%
\pgfpathlineto{\pgfqpoint{1.587954in}{1.039192in}}%
\pgfpathlineto{\pgfqpoint{1.638279in}{1.016359in}}%
\pgfpathlineto{\pgfqpoint{1.688604in}{0.998758in}}%
\pgfpathlineto{\pgfqpoint{1.738928in}{0.986175in}}%
\pgfpathlineto{\pgfqpoint{1.789253in}{0.978001in}}%
\pgfpathlineto{\pgfqpoint{1.839578in}{0.972950in}}%
\pgfpathlineto{\pgfqpoint{1.889902in}{0.970255in}}%
\pgfpathlineto{\pgfqpoint{1.940227in}{0.968752in}}%
\pgfpathlineto{\pgfqpoint{1.990552in}{0.967891in}}%
\pgfpathlineto{\pgfqpoint{2.040876in}{0.967411in}}%
\pgfpathlineto{\pgfqpoint{2.091201in}{0.967117in}}%
\pgfpathlineto{\pgfqpoint{2.141526in}{0.961314in}}%
\pgfpathlineto{\pgfqpoint{2.191850in}{0.947874in}}%
\pgfpathlineto{\pgfqpoint{2.242175in}{0.927153in}}%
\pgfpathlineto{\pgfqpoint{2.292500in}{0.918776in}}%
\pgfpathlineto{\pgfqpoint{2.342824in}{0.915228in}}%
\pgfpathlineto{\pgfqpoint{2.393149in}{0.903408in}}%
\pgfpathlineto{\pgfqpoint{2.443474in}{0.896210in}}%
\pgfpathlineto{\pgfqpoint{2.493798in}{0.882339in}}%
\pgfpathlineto{\pgfqpoint{2.544123in}{0.870418in}}%
\pgfpathlineto{\pgfqpoint{2.594448in}{0.857021in}}%
\pgfpathlineto{\pgfqpoint{2.644772in}{0.845113in}}%
\pgfpathlineto{\pgfqpoint{2.695097in}{0.833316in}}%
\pgfpathlineto{\pgfqpoint{2.745422in}{0.822571in}}%
\pgfpathlineto{\pgfqpoint{2.795746in}{0.811697in}}%
\pgfpathlineto{\pgfqpoint{2.846071in}{0.802364in}}%
\pgfpathlineto{\pgfqpoint{2.896396in}{0.792707in}}%
\pgfpathlineto{\pgfqpoint{2.946721in}{0.792722in}}%
\pgfpathlineto{\pgfqpoint{2.997045in}{0.792624in}}%
\pgfpathlineto{\pgfqpoint{3.047370in}{0.792395in}}%
\pgfpathlineto{\pgfqpoint{3.097695in}{0.791849in}}%
\pgfpathlineto{\pgfqpoint{3.148019in}{0.790554in}}%
\pgfpathlineto{\pgfqpoint{3.198344in}{0.787755in}}%
\pgfpathlineto{\pgfqpoint{3.248669in}{0.782462in}}%
\pgfpathlineto{\pgfqpoint{3.298993in}{0.774244in}}%
\pgfpathlineto{\pgfqpoint{3.349318in}{0.763969in}}%
\pgfpathlineto{\pgfqpoint{3.399643in}{0.753745in}}%
\pgfpathlineto{\pgfqpoint{3.449967in}{0.745033in}}%
\pgfpathlineto{\pgfqpoint{3.500292in}{0.737098in}}%
\pgfpathlineto{\pgfqpoint{3.550617in}{0.730132in}}%
\pgfpathlineto{\pgfqpoint{3.600941in}{0.724254in}}%
\pgfpathlineto{\pgfqpoint{3.651266in}{0.719371in}}%
\pgfpathlineto{\pgfqpoint{3.701591in}{0.716562in}}%
\pgfpathlineto{\pgfqpoint{3.751915in}{0.713729in}}%
\pgfpathlineto{\pgfqpoint{3.802240in}{0.709633in}}%
\pgfpathlineto{\pgfqpoint{3.852565in}{0.704959in}}%
\pgfpathlineto{\pgfqpoint{3.902889in}{0.700498in}}%
\pgfpathlineto{\pgfqpoint{3.953214in}{0.700483in}}%
\pgfpathlineto{\pgfqpoint{4.003539in}{0.700454in}}%
\pgfpathlineto{\pgfqpoint{4.053863in}{0.700369in}}%
\pgfpathlineto{\pgfqpoint{4.104188in}{0.700210in}}%
\pgfpathlineto{\pgfqpoint{4.154513in}{0.699967in}}%
\pgfpathlineto{\pgfqpoint{4.204837in}{0.699631in}}%
\pgfpathlineto{\pgfqpoint{4.255162in}{0.699172in}}%
\pgfpathlineto{\pgfqpoint{4.305487in}{0.698496in}}%
\pgfpathlineto{\pgfqpoint{4.355811in}{0.697387in}}%
\pgfusepath{stroke}%
\end{pgfscope}%
\begin{pgfscope}%
\pgfsetrectcap%
\pgfsetmiterjoin%
\pgfsetlinewidth{0.803000pt}%
\definecolor{currentstroke}{rgb}{0.000000,0.000000,0.000000}%
\pgfsetstrokecolor{currentstroke}%
\pgfsetdash{}{0pt}%
\pgfpathmoveto{\pgfqpoint{0.656948in}{0.504167in}}%
\pgfpathlineto{\pgfqpoint{4.531948in}{0.504167in}}%
\pgfusepath{stroke}%
\end{pgfscope}%
\begin{pgfscope}%
\definecolor{textcolor}{rgb}{0.341176,0.317647,0.317647}%
\pgfsetstrokecolor{textcolor}%
\pgfsetfillcolor{textcolor}%
\pgftext[x=2.544123in,y=0.917864in,left,base]{\color{textcolor}{\rmfamily\fontsize{10.000000}{12.000000}\selectfont\catcode`\^=\active\def^{\ifmmode\sp\else\^{}\fi}\catcode`\%=\active\def
\end{pgfscope}%
\begin{pgfscope}%
\definecolor{textcolor}{rgb}{0.341176,0.317647,0.317647}%
\pgfsetstrokecolor{textcolor}%
\pgfsetfillcolor{textcolor}%
\pgftext[x=3.550617in,y=0.589459in,left,base]{\color{textcolor}{\rmfamily\fontsize{10.000000}{12.000000}\selectfont\catcode`\^=\active\def^{\ifmmode\sp\else\^{}\fi}\catcode`\%=\active\def
\end{pgfscope}%
\begin{pgfscope}%
\definecolor{textcolor}{rgb}{0.341176,0.317647,0.317647}%
\pgfsetstrokecolor{textcolor}%
\pgfsetfillcolor{textcolor}%
\pgftext[x=4.053863in,y=0.755137in,left,base]{\color{textcolor}{\rmfamily\fontsize{10.000000}{12.000000}\selectfont\catcode`\^=\active\def^{\ifmmode\sp\else\^{}\fi}\catcode`\%=\active\def
\end{pgfscope}%
\begin{pgfscope}%
\definecolor{textcolor}{rgb}{1.000000,0.121569,0.352941}%
\pgfsetstrokecolor{textcolor}%
\pgfsetfillcolor{textcolor}%
\pgftext[x=1.575373in,y=2.698667in,left,base]{\color{textcolor}{\rmfamily\fontsize{10.000000}{12.000000}\selectfont\catcode`\^=\active\def^{\ifmmode\sp\else\^{}\fi}\catcode`\%=\active\def
\end{pgfscope}%
\begin{pgfscope}%
\definecolor{textcolor}{rgb}{0.090196,0.725490,0.470588}%
\pgfsetstrokecolor{textcolor}%
\pgfsetfillcolor{textcolor}%
\pgftext[x=3.085113in,y=2.698667in,left,base]{\color{textcolor}{\rmfamily\fontsize{10.000000}{12.000000}\selectfont\catcode`\^=\active\def^{\ifmmode\sp\else\^{}\fi}\catcode`\%=\active\def
\end{pgfscope}%
\begin{pgfscope}%
\definecolor{textcolor}{rgb}{0.282353,0.184314,0.968627}%
\pgfsetstrokecolor{textcolor}%
\pgfsetfillcolor{textcolor}%
\pgftext[x=3.839983in,y=2.698667in,left,base]{\color{textcolor}{\rmfamily\fontsize{10.000000}{12.000000}\selectfont\catcode`\^=\active\def^{\ifmmode\sp\else\^{}\fi}\catcode`\%=\active\def
\end{pgfscope}%
\begin{pgfscope}%
\pgfsetbuttcap%
\pgfsetmiterjoin%
\definecolor{currentfill}{rgb}{1.000000,1.000000,1.000000}%
\pgfsetfillcolor{currentfill}%
\pgfsetfillopacity{0.800000}%
\pgfsetlinewidth{0.240900pt}%
\definecolor{currentstroke}{rgb}{0.800000,0.800000,0.800000}%
\pgfsetstrokecolor{currentstroke}%
\pgfsetstrokeopacity{0.800000}%
\pgfsetdash{}{0pt}%
\pgfpathmoveto{\pgfqpoint{0.734726in}{0.559722in}}%
\pgfpathlineto{\pgfqpoint{1.327052in}{0.559722in}}%
\pgfpathquadraticcurveto{\pgfqpoint{1.349274in}{0.559722in}}{\pgfqpoint{1.349274in}{0.581944in}}%
\pgfpathlineto{\pgfqpoint{1.349274in}{0.725772in}}%
\pgfpathquadraticcurveto{\pgfqpoint{1.349274in}{0.747994in}}{\pgfqpoint{1.327052in}{0.747994in}}%
\pgfpathlineto{\pgfqpoint{0.734726in}{0.747994in}}%
\pgfpathquadraticcurveto{\pgfqpoint{0.712503in}{0.747994in}}{\pgfqpoint{0.712503in}{0.725772in}}%
\pgfpathlineto{\pgfqpoint{0.712503in}{0.581944in}}%
\pgfpathquadraticcurveto{\pgfqpoint{0.712503in}{0.559722in}}{\pgfqpoint{0.734726in}{0.559722in}}%
\pgfpathlineto{\pgfqpoint{0.734726in}{0.559722in}}%
\pgfpathclose%
\pgfusepath{stroke,fill}%
\end{pgfscope}%
\begin{pgfscope}%
\pgfsetrectcap%
\pgfsetroundjoin%
\pgfsetlinewidth{1.405250pt}%
\definecolor{currentstroke}{rgb}{0.341176,0.317647,0.317647}%
\pgfsetstrokecolor{currentstroke}%
\pgfsetdash{}{0pt}%
\pgfpathmoveto{\pgfqpoint{0.756948in}{0.664661in}}%
\pgfpathlineto{\pgfqpoint{0.868059in}{0.664661in}}%
\pgfpathlineto{\pgfqpoint{0.979170in}{0.664661in}}%
\pgfusepath{stroke}%
\end{pgfscope}%
\begin{pgfscope}%
\definecolor{textcolor}{rgb}{0.000000,0.000000,0.000000}%
\pgfsetstrokecolor{textcolor}%
\pgfsetfillcolor{textcolor}%
\pgftext[x=1.068059in,y=0.625772in,left,base]{\color{textcolor}{\rmfamily\fontsize{8.000000}{9.600000}\selectfont\catcode`\^=\active\def^{\ifmmode\sp\else\^{}\fi}\catcode`\%=\active\def
\end{pgfscope}%
\end{pgfpicture}%
\makeatother%
\endgroup%

%% file: Weak_Scaling.pgf
\begingroup%
\makeatletter%
\begin{pgfpicture}%
\pgfpathrectangle{\pgfpointorigin}{\pgfqpoint{4.126391in}{2.902917in}}%
\pgfusepath{use as bounding box, clip}%
\begin{pgfscope}%
\pgfsetbuttcap%
\pgfsetmiterjoin%
\definecolor{currentfill}{rgb}{1.000000,1.000000,1.000000}%
\pgfsetfillcolor{currentfill}%
\pgfsetlinewidth{0.000000pt}%
\definecolor{currentstroke}{rgb}{1.000000,1.000000,1.000000}%
\pgfsetstrokecolor{currentstroke}%
\pgfsetdash{}{0pt}%
\pgfpathmoveto{\pgfqpoint{0.000000in}{-0.000000in}}%
\pgfpathlineto{\pgfqpoint{4.126391in}{-0.000000in}}%
\pgfpathlineto{\pgfqpoint{4.126391in}{2.902917in}}%
\pgfpathlineto{\pgfqpoint{0.000000in}{2.902917in}}%
\pgfpathlineto{\pgfqpoint{0.000000in}{-0.000000in}}%
\pgfpathclose%
\pgfusepath{fill}%
\end{pgfscope}%
\begin{pgfscope}%
\pgfsetbuttcap%
\pgfsetmiterjoin%
\definecolor{currentfill}{rgb}{1.000000,1.000000,1.000000}%
\pgfsetfillcolor{currentfill}%
\pgfsetlinewidth{0.000000pt}%
\definecolor{currentstroke}{rgb}{0.000000,0.000000,0.000000}%
\pgfsetstrokecolor{currentstroke}%
\pgfsetstrokeopacity{0.000000}%
\pgfsetdash{}{0pt}%
\pgfpathmoveto{\pgfqpoint{0.538891in}{0.548472in}}%
\pgfpathlineto{\pgfqpoint{4.026391in}{0.548472in}}%
\pgfpathlineto{\pgfqpoint{4.026391in}{2.473472in}}%
\pgfpathlineto{\pgfqpoint{0.538891in}{2.473472in}}%
\pgfpathlineto{\pgfqpoint{0.538891in}{0.548472in}}%
\pgfpathclose%
\pgfusepath{fill}%
\end{pgfscope}%
\begin{pgfscope}%
\pgfpathrectangle{\pgfqpoint{0.538891in}{0.548472in}}{\pgfqpoint{3.487500in}{1.925000in}}%
\pgfusepath{clip}%
\pgfsetbuttcap%
\pgfsetmiterjoin%
\definecolor{currentfill}{rgb}{0.129412,0.901961,0.756863}%
\pgfsetfillcolor{currentfill}%
\pgfsetlinewidth{0.000000pt}%
\definecolor{currentstroke}{rgb}{0.000000,0.000000,0.000000}%
\pgfsetstrokecolor{currentstroke}%
\pgfsetstrokeopacity{0.000000}%
\pgfsetdash{}{0pt}%
\pgfpathmoveto{\pgfqpoint{0.697413in}{0.548472in}}%
\pgfpathlineto{\pgfqpoint{1.028655in}{0.548472in}}%
\pgfpathlineto{\pgfqpoint{1.028655in}{0.829825in}}%
\pgfpathlineto{\pgfqpoint{0.697413in}{0.829825in}}%
\pgfpathlineto{\pgfqpoint{0.697413in}{0.548472in}}%
\pgfpathclose%
\pgfusepath{fill}%
\end{pgfscope}%
\begin{pgfscope}%
\pgfpathrectangle{\pgfqpoint{0.538891in}{0.548472in}}{\pgfqpoint{3.487500in}{1.925000in}}%
\pgfusepath{clip}%
\pgfsetbuttcap%
\pgfsetmiterjoin%
\definecolor{currentfill}{rgb}{0.129412,0.901961,0.756863}%
\pgfsetfillcolor{currentfill}%
\pgfsetlinewidth{0.000000pt}%
\definecolor{currentstroke}{rgb}{0.000000,0.000000,0.000000}%
\pgfsetstrokecolor{currentstroke}%
\pgfsetstrokeopacity{0.000000}%
\pgfsetdash{}{0pt}%
\pgfpathmoveto{\pgfqpoint{1.643818in}{0.548472in}}%
\pgfpathlineto{\pgfqpoint{1.975059in}{0.548472in}}%
\pgfpathlineto{\pgfqpoint{1.975059in}{0.838901in}}%
\pgfpathlineto{\pgfqpoint{1.643818in}{0.838901in}}%
\pgfpathlineto{\pgfqpoint{1.643818in}{0.548472in}}%
\pgfpathclose%
\pgfusepath{fill}%
\end{pgfscope}%
\begin{pgfscope}%
\pgfpathrectangle{\pgfqpoint{0.538891in}{0.548472in}}{\pgfqpoint{3.487500in}{1.925000in}}%
\pgfusepath{clip}%
\pgfsetbuttcap%
\pgfsetmiterjoin%
\definecolor{currentfill}{rgb}{0.129412,0.901961,0.756863}%
\pgfsetfillcolor{currentfill}%
\pgfsetlinewidth{0.000000pt}%
\definecolor{currentstroke}{rgb}{0.000000,0.000000,0.000000}%
\pgfsetstrokecolor{currentstroke}%
\pgfsetstrokeopacity{0.000000}%
\pgfsetdash{}{0pt}%
\pgfpathmoveto{\pgfqpoint{2.590222in}{0.548472in}}%
\pgfpathlineto{\pgfqpoint{2.921463in}{0.548472in}}%
\pgfpathlineto{\pgfqpoint{2.921463in}{0.857053in}}%
\pgfpathlineto{\pgfqpoint{2.590222in}{0.857053in}}%
\pgfpathlineto{\pgfqpoint{2.590222in}{0.548472in}}%
\pgfpathclose%
\pgfusepath{fill}%
\end{pgfscope}%
\begin{pgfscope}%
\pgfpathrectangle{\pgfqpoint{0.538891in}{0.548472in}}{\pgfqpoint{3.487500in}{1.925000in}}%
\pgfusepath{clip}%
\pgfsetbuttcap%
\pgfsetmiterjoin%
\definecolor{currentfill}{rgb}{0.129412,0.901961,0.756863}%
\pgfsetfillcolor{currentfill}%
\pgfsetlinewidth{0.000000pt}%
\definecolor{currentstroke}{rgb}{0.000000,0.000000,0.000000}%
\pgfsetstrokecolor{currentstroke}%
\pgfsetstrokeopacity{0.000000}%
\pgfsetdash{}{0pt}%
\pgfpathmoveto{\pgfqpoint{3.536626in}{0.548472in}}%
\pgfpathlineto{\pgfqpoint{3.867868in}{0.548472in}}%
\pgfpathlineto{\pgfqpoint{3.867868in}{0.893357in}}%
\pgfpathlineto{\pgfqpoint{3.536626in}{0.893357in}}%
\pgfpathlineto{\pgfqpoint{3.536626in}{0.548472in}}%
\pgfpathclose%
\pgfusepath{fill}%
\end{pgfscope}%
\begin{pgfscope}%
\pgfpathrectangle{\pgfqpoint{0.538891in}{0.548472in}}{\pgfqpoint{3.487500in}{1.925000in}}%
\pgfusepath{clip}%
\pgfsetbuttcap%
\pgfsetmiterjoin%
\definecolor{currentfill}{rgb}{0.152941,0.556863,0.647059}%
\pgfsetfillcolor{currentfill}%
\pgfsetlinewidth{0.000000pt}%
\definecolor{currentstroke}{rgb}{0.000000,0.000000,0.000000}%
\pgfsetstrokecolor{currentstroke}%
\pgfsetstrokeopacity{0.000000}%
\pgfsetdash{}{0pt}%
\pgfpathmoveto{\pgfqpoint{0.697413in}{0.829825in}}%
\pgfpathlineto{\pgfqpoint{1.028655in}{0.829825in}}%
\pgfpathlineto{\pgfqpoint{1.028655in}{2.027845in}}%
\pgfpathlineto{\pgfqpoint{0.697413in}{2.027845in}}%
\pgfpathlineto{\pgfqpoint{0.697413in}{0.829825in}}%
\pgfpathclose%
\pgfusepath{fill}%
\end{pgfscope}%
\begin{pgfscope}%
\pgfpathrectangle{\pgfqpoint{0.538891in}{0.548472in}}{\pgfqpoint{3.487500in}{1.925000in}}%
\pgfusepath{clip}%
\pgfsetbuttcap%
\pgfsetmiterjoin%
\definecolor{currentfill}{rgb}{0.152941,0.556863,0.647059}%
\pgfsetfillcolor{currentfill}%
\pgfsetlinewidth{0.000000pt}%
\definecolor{currentstroke}{rgb}{0.000000,0.000000,0.000000}%
\pgfsetstrokecolor{currentstroke}%
\pgfsetstrokeopacity{0.000000}%
\pgfsetdash{}{0pt}%
\pgfpathmoveto{\pgfqpoint{1.643818in}{0.838901in}}%
\pgfpathlineto{\pgfqpoint{1.975059in}{0.838901in}}%
\pgfpathlineto{\pgfqpoint{1.975059in}{1.038571in}}%
\pgfpathlineto{\pgfqpoint{1.643818in}{1.038571in}}%
\pgfpathlineto{\pgfqpoint{1.643818in}{0.838901in}}%
\pgfpathclose%
\pgfusepath{fill}%
\end{pgfscope}%
\begin{pgfscope}%
\pgfpathrectangle{\pgfqpoint{0.538891in}{0.548472in}}{\pgfqpoint{3.487500in}{1.925000in}}%
\pgfusepath{clip}%
\pgfsetbuttcap%
\pgfsetmiterjoin%
\definecolor{currentfill}{rgb}{0.152941,0.556863,0.647059}%
\pgfsetfillcolor{currentfill}%
\pgfsetlinewidth{0.000000pt}%
\definecolor{currentstroke}{rgb}{0.000000,0.000000,0.000000}%
\pgfsetstrokecolor{currentstroke}%
\pgfsetstrokeopacity{0.000000}%
\pgfsetdash{}{0pt}%
\pgfpathmoveto{\pgfqpoint{2.590222in}{0.857053in}}%
\pgfpathlineto{\pgfqpoint{2.921463in}{0.857053in}}%
\pgfpathlineto{\pgfqpoint{2.921463in}{1.174710in}}%
\pgfpathlineto{\pgfqpoint{2.590222in}{1.174710in}}%
\pgfpathlineto{\pgfqpoint{2.590222in}{0.857053in}}%
\pgfpathclose%
\pgfusepath{fill}%
\end{pgfscope}%
\begin{pgfscope}%
\pgfpathrectangle{\pgfqpoint{0.538891in}{0.548472in}}{\pgfqpoint{3.487500in}{1.925000in}}%
\pgfusepath{clip}%
\pgfsetbuttcap%
\pgfsetmiterjoin%
\definecolor{currentfill}{rgb}{0.152941,0.556863,0.647059}%
\pgfsetfillcolor{currentfill}%
\pgfsetlinewidth{0.000000pt}%
\definecolor{currentstroke}{rgb}{0.000000,0.000000,0.000000}%
\pgfsetstrokecolor{currentstroke}%
\pgfsetstrokeopacity{0.000000}%
\pgfsetdash{}{0pt}%
\pgfpathmoveto{\pgfqpoint{3.536626in}{0.893357in}}%
\pgfpathlineto{\pgfqpoint{3.867868in}{0.893357in}}%
\pgfpathlineto{\pgfqpoint{3.867868in}{1.873555in}}%
\pgfpathlineto{\pgfqpoint{3.536626in}{1.873555in}}%
\pgfpathlineto{\pgfqpoint{3.536626in}{0.893357in}}%
\pgfpathclose%
\pgfusepath{fill}%
\end{pgfscope}%
\begin{pgfscope}%
\pgfpathrectangle{\pgfqpoint{0.538891in}{0.548472in}}{\pgfqpoint{3.487500in}{1.925000in}}%
\pgfusepath{clip}%
\pgfsetbuttcap%
\pgfsetmiterjoin%
\definecolor{currentfill}{rgb}{0.121569,0.258824,0.529412}%
\pgfsetfillcolor{currentfill}%
\pgfsetlinewidth{0.000000pt}%
\definecolor{currentstroke}{rgb}{0.000000,0.000000,0.000000}%
\pgfsetstrokecolor{currentstroke}%
\pgfsetstrokeopacity{0.000000}%
\pgfsetdash{}{0pt}%
\pgfpathmoveto{\pgfqpoint{0.697413in}{2.027845in}}%
\pgfpathlineto{\pgfqpoint{1.028655in}{2.027845in}}%
\pgfpathlineto{\pgfqpoint{1.028655in}{2.182135in}}%
\pgfpathlineto{\pgfqpoint{0.697413in}{2.182135in}}%
\pgfpathlineto{\pgfqpoint{0.697413in}{2.027845in}}%
\pgfpathclose%
\pgfusepath{fill}%
\end{pgfscope}%
\begin{pgfscope}%
\pgfpathrectangle{\pgfqpoint{0.538891in}{0.548472in}}{\pgfqpoint{3.487500in}{1.925000in}}%
\pgfusepath{clip}%
\pgfsetbuttcap%
\pgfsetmiterjoin%
\definecolor{currentfill}{rgb}{0.121569,0.258824,0.529412}%
\pgfsetfillcolor{currentfill}%
\pgfsetlinewidth{0.000000pt}%
\definecolor{currentstroke}{rgb}{0.000000,0.000000,0.000000}%
\pgfsetstrokecolor{currentstroke}%
\pgfsetstrokeopacity{0.000000}%
\pgfsetdash{}{0pt}%
\pgfpathmoveto{\pgfqpoint{1.643818in}{1.038571in}}%
\pgfpathlineto{\pgfqpoint{1.975059in}{1.038571in}}%
\pgfpathlineto{\pgfqpoint{1.975059in}{1.338076in}}%
\pgfpathlineto{\pgfqpoint{1.643818in}{1.338076in}}%
\pgfpathlineto{\pgfqpoint{1.643818in}{1.038571in}}%
\pgfpathclose%
\pgfusepath{fill}%
\end{pgfscope}%
\begin{pgfscope}%
\pgfpathrectangle{\pgfqpoint{0.538891in}{0.548472in}}{\pgfqpoint{3.487500in}{1.925000in}}%
\pgfusepath{clip}%
\pgfsetbuttcap%
\pgfsetmiterjoin%
\definecolor{currentfill}{rgb}{0.121569,0.258824,0.529412}%
\pgfsetfillcolor{currentfill}%
\pgfsetlinewidth{0.000000pt}%
\definecolor{currentstroke}{rgb}{0.000000,0.000000,0.000000}%
\pgfsetstrokecolor{currentstroke}%
\pgfsetstrokeopacity{0.000000}%
\pgfsetdash{}{0pt}%
\pgfpathmoveto{\pgfqpoint{2.590222in}{1.174710in}}%
\pgfpathlineto{\pgfqpoint{2.921463in}{1.174710in}}%
\pgfpathlineto{\pgfqpoint{2.921463in}{1.610353in}}%
\pgfpathlineto{\pgfqpoint{2.590222in}{1.610353in}}%
\pgfpathlineto{\pgfqpoint{2.590222in}{1.174710in}}%
\pgfpathclose%
\pgfusepath{fill}%
\end{pgfscope}%
\begin{pgfscope}%
\pgfpathrectangle{\pgfqpoint{0.538891in}{0.548472in}}{\pgfqpoint{3.487500in}{1.925000in}}%
\pgfusepath{clip}%
\pgfsetbuttcap%
\pgfsetmiterjoin%
\definecolor{currentfill}{rgb}{0.121569,0.258824,0.529412}%
\pgfsetfillcolor{currentfill}%
\pgfsetlinewidth{0.000000pt}%
\definecolor{currentstroke}{rgb}{0.000000,0.000000,0.000000}%
\pgfsetstrokecolor{currentstroke}%
\pgfsetstrokeopacity{0.000000}%
\pgfsetdash{}{0pt}%
\pgfpathmoveto{\pgfqpoint{3.536626in}{1.873555in}}%
\pgfpathlineto{\pgfqpoint{3.867868in}{1.873555in}}%
\pgfpathlineto{\pgfqpoint{3.867868in}{2.381805in}}%
\pgfpathlineto{\pgfqpoint{3.536626in}{2.381805in}}%
\pgfpathlineto{\pgfqpoint{3.536626in}{1.873555in}}%
\pgfpathclose%
\pgfusepath{fill}%
\end{pgfscope}%
\begin{pgfscope}%
\pgfsetbuttcap%
\pgfsetroundjoin%
\definecolor{currentfill}{rgb}{0.000000,0.000000,0.000000}%
\pgfsetfillcolor{currentfill}%
\pgfsetlinewidth{0.803000pt}%
\definecolor{currentstroke}{rgb}{0.000000,0.000000,0.000000}%
\pgfsetstrokecolor{currentstroke}%
\pgfsetdash{}{0pt}%
\pgfsys@defobject{currentmarker}{\pgfqpoint{0.000000in}{-0.048611in}}{\pgfqpoint{0.000000in}{0.000000in}}{%
\pgfpathmoveto{\pgfqpoint{0.000000in}{0.000000in}}%
\pgfpathlineto{\pgfqpoint{0.000000in}{-0.048611in}}%
\pgfusepath{stroke,fill}%
}%
\begin{pgfscope}%
\pgfsys@transformshift{0.863034in}{0.548472in}%
\pgfsys@useobject{currentmarker}{}%
\end{pgfscope}%
\end{pgfscope}%
\begin{pgfscope}%
\definecolor{textcolor}{rgb}{0.000000,0.000000,0.000000}%
\pgfsetstrokecolor{textcolor}%
\pgfsetfillcolor{textcolor}%
\pgftext[x=0.863034in,y=0.422083in,,top]{\color{textcolor}{\rmfamily\fontsize{8.000000}{9.600000}\selectfont\catcode`\^=\active\def^{\ifmmode\sp\else\^{}\fi}\catcode`\%=\active\def
\end{pgfscope}%
\begin{pgfscope}%
\pgfsetbuttcap%
\pgfsetroundjoin%
\definecolor{currentfill}{rgb}{0.000000,0.000000,0.000000}%
\pgfsetfillcolor{currentfill}%
\pgfsetlinewidth{0.803000pt}%
\definecolor{currentstroke}{rgb}{0.000000,0.000000,0.000000}%
\pgfsetstrokecolor{currentstroke}%
\pgfsetdash{}{0pt}%
\pgfsys@defobject{currentmarker}{\pgfqpoint{0.000000in}{-0.048611in}}{\pgfqpoint{0.000000in}{0.000000in}}{%
\pgfpathmoveto{\pgfqpoint{0.000000in}{0.000000in}}%
\pgfpathlineto{\pgfqpoint{0.000000in}{-0.048611in}}%
\pgfusepath{stroke,fill}%
}%
\begin{pgfscope}%
\pgfsys@transformshift{1.809438in}{0.548472in}%
\pgfsys@useobject{currentmarker}{}%
\end{pgfscope}%
\end{pgfscope}%
\begin{pgfscope}%
\definecolor{textcolor}{rgb}{0.000000,0.000000,0.000000}%
\pgfsetstrokecolor{textcolor}%
\pgfsetfillcolor{textcolor}%
\pgftext[x=1.809438in,y=0.422083in,,top]{\color{textcolor}{\rmfamily\fontsize{8.000000}{9.600000}\selectfont\catcode`\^=\active\def^{\ifmmode\sp\else\^{}\fi}\catcode`\%=\active\def
\end{pgfscope}%
\begin{pgfscope}%
\pgfsetbuttcap%
\pgfsetroundjoin%
\definecolor{currentfill}{rgb}{0.000000,0.000000,0.000000}%
\pgfsetfillcolor{currentfill}%
\pgfsetlinewidth{0.803000pt}%
\definecolor{currentstroke}{rgb}{0.000000,0.000000,0.000000}%
\pgfsetstrokecolor{currentstroke}%
\pgfsetdash{}{0pt}%
\pgfsys@defobject{currentmarker}{\pgfqpoint{0.000000in}{-0.048611in}}{\pgfqpoint{0.000000in}{0.000000in}}{%
\pgfpathmoveto{\pgfqpoint{0.000000in}{0.000000in}}%
\pgfpathlineto{\pgfqpoint{0.000000in}{-0.048611in}}%
\pgfusepath{stroke,fill}%
}%
\begin{pgfscope}%
\pgfsys@transformshift{2.755843in}{0.548472in}%
\pgfsys@useobject{currentmarker}{}%
\end{pgfscope}%
\end{pgfscope}%
\begin{pgfscope}%
\definecolor{textcolor}{rgb}{0.000000,0.000000,0.000000}%
\pgfsetstrokecolor{textcolor}%
\pgfsetfillcolor{textcolor}%
\pgftext[x=2.755843in,y=0.422083in,,top]{\color{textcolor}{\rmfamily\fontsize{8.000000}{9.600000}\selectfont\catcode`\^=\active\def^{\ifmmode\sp\else\^{}\fi}\catcode`\%=\active\def
\end{pgfscope}%
\begin{pgfscope}%
\pgfsetbuttcap%
\pgfsetroundjoin%
\definecolor{currentfill}{rgb}{0.000000,0.000000,0.000000}%
\pgfsetfillcolor{currentfill}%
\pgfsetlinewidth{0.803000pt}%
\definecolor{currentstroke}{rgb}{0.000000,0.000000,0.000000}%
\pgfsetstrokecolor{currentstroke}%
\pgfsetdash{}{0pt}%
\pgfsys@defobject{currentmarker}{\pgfqpoint{0.000000in}{-0.048611in}}{\pgfqpoint{0.000000in}{0.000000in}}{%
\pgfpathmoveto{\pgfqpoint{0.000000in}{0.000000in}}%
\pgfpathlineto{\pgfqpoint{0.000000in}{-0.048611in}}%
\pgfusepath{stroke,fill}%
}%
\begin{pgfscope}%
\pgfsys@transformshift{3.702247in}{0.548472in}%
\pgfsys@useobject{currentmarker}{}%
\end{pgfscope}%
\end{pgfscope}%
\begin{pgfscope}%
\definecolor{textcolor}{rgb}{0.000000,0.000000,0.000000}%
\pgfsetstrokecolor{textcolor}%
\pgfsetfillcolor{textcolor}%
\pgftext[x=3.702247in,y=0.422083in,,top]{\color{textcolor}{\rmfamily\fontsize{8.000000}{9.600000}\selectfont\catcode`\^=\active\def^{\ifmmode\sp\else\^{}\fi}\catcode`\%=\active\def
\end{pgfscope}%
\begin{pgfscope}%
\definecolor{textcolor}{rgb}{0.000000,0.000000,0.000000}%
\pgfsetstrokecolor{textcolor}%
\pgfsetfillcolor{textcolor}%
\pgftext[x=2.282641in,y=0.238889in,,top]{\color{textcolor}{\rmfamily\fontsize{10.000000}{12.000000}\selectfont\catcode`\^=\active\def^{\ifmmode\sp\else\^{}\fi}\catcode`\%=\active\def
\end{pgfscope}%
\begin{pgfscope}%
\pgfpathrectangle{\pgfqpoint{0.538891in}{0.548472in}}{\pgfqpoint{3.487500in}{1.925000in}}%
\pgfusepath{clip}%
\pgfsetrectcap%
\pgfsetroundjoin%
\pgfsetlinewidth{0.803000pt}%
\definecolor{currentstroke}{rgb}{0.690196,0.690196,0.690196}%
\pgfsetstrokecolor{currentstroke}%
\pgfsetdash{}{0pt}%
\pgfpathmoveto{\pgfqpoint{0.538891in}{0.548472in}}%
\pgfpathlineto{\pgfqpoint{4.026391in}{0.548472in}}%
\pgfusepath{stroke}%
\end{pgfscope}%
\begin{pgfscope}%
\pgfsetbuttcap%
\pgfsetroundjoin%
\definecolor{currentfill}{rgb}{0.000000,0.000000,0.000000}%
\pgfsetfillcolor{currentfill}%
\pgfsetlinewidth{0.803000pt}%
\definecolor{currentstroke}{rgb}{0.000000,0.000000,0.000000}%
\pgfsetstrokecolor{currentstroke}%
\pgfsetdash{}{0pt}%
\pgfsys@defobject{currentmarker}{\pgfqpoint{-0.048611in}{0.000000in}}{\pgfqpoint{-0.000000in}{0.000000in}}{%
\pgfpathmoveto{\pgfqpoint{-0.000000in}{0.000000in}}%
\pgfpathlineto{\pgfqpoint{-0.048611in}{0.000000in}}%
\pgfusepath{stroke,fill}%
}%
\begin{pgfscope}%
\pgfsys@transformshift{0.538891in}{0.548472in}%
\pgfsys@useobject{currentmarker}{}%
\end{pgfscope}%
\end{pgfscope}%
\begin{pgfscope}%
\definecolor{textcolor}{rgb}{0.000000,0.000000,0.000000}%
\pgfsetstrokecolor{textcolor}%
\pgfsetfillcolor{textcolor}%
\pgftext[x=0.353473in, y=0.509892in, left, base]{\color{textcolor}{\rmfamily\fontsize{8.000000}{9.600000}\selectfont\catcode`\^=\active\def^{\ifmmode\sp\else\^{}\fi}\catcode`\%=\active\def
\end{pgfscope}%
\begin{pgfscope}%
\pgfpathrectangle{\pgfqpoint{0.538891in}{0.548472in}}{\pgfqpoint{3.487500in}{1.925000in}}%
\pgfusepath{clip}%
\pgfsetrectcap%
\pgfsetroundjoin%
\pgfsetlinewidth{0.803000pt}%
\definecolor{currentstroke}{rgb}{0.690196,0.690196,0.690196}%
\pgfsetstrokecolor{currentstroke}%
\pgfsetdash{}{0pt}%
\pgfpathmoveto{\pgfqpoint{0.538891in}{1.002267in}}%
\pgfpathlineto{\pgfqpoint{4.026391in}{1.002267in}}%
\pgfusepath{stroke}%
\end{pgfscope}%
\begin{pgfscope}%
\pgfsetbuttcap%
\pgfsetroundjoin%
\definecolor{currentfill}{rgb}{0.000000,0.000000,0.000000}%
\pgfsetfillcolor{currentfill}%
\pgfsetlinewidth{0.803000pt}%
\definecolor{currentstroke}{rgb}{0.000000,0.000000,0.000000}%
\pgfsetstrokecolor{currentstroke}%
\pgfsetdash{}{0pt}%
\pgfsys@defobject{currentmarker}{\pgfqpoint{-0.048611in}{0.000000in}}{\pgfqpoint{-0.000000in}{0.000000in}}{%
\pgfpathmoveto{\pgfqpoint{-0.000000in}{0.000000in}}%
\pgfpathlineto{\pgfqpoint{-0.048611in}{0.000000in}}%
\pgfusepath{stroke,fill}%
}%
\begin{pgfscope}%
\pgfsys@transformshift{0.538891in}{1.002267in}%
\pgfsys@useobject{currentmarker}{}%
\end{pgfscope}%
\end{pgfscope}%
\begin{pgfscope}%
\definecolor{textcolor}{rgb}{0.000000,0.000000,0.000000}%
\pgfsetstrokecolor{textcolor}%
\pgfsetfillcolor{textcolor}%
\pgftext[x=0.353473in, y=0.963687in, left, base]{\color{textcolor}{\rmfamily\fontsize{8.000000}{9.600000}\selectfont\catcode`\^=\active\def^{\ifmmode\sp\else\^{}\fi}\catcode`\%=\active\def
\end{pgfscope}%
\begin{pgfscope}%
\pgfpathrectangle{\pgfqpoint{0.538891in}{0.548472in}}{\pgfqpoint{3.487500in}{1.925000in}}%
\pgfusepath{clip}%
\pgfsetrectcap%
\pgfsetroundjoin%
\pgfsetlinewidth{0.803000pt}%
\definecolor{currentstroke}{rgb}{0.690196,0.690196,0.690196}%
\pgfsetstrokecolor{currentstroke}%
\pgfsetdash{}{0pt}%
\pgfpathmoveto{\pgfqpoint{0.538891in}{1.456063in}}%
\pgfpathlineto{\pgfqpoint{4.026391in}{1.456063in}}%
\pgfusepath{stroke}%
\end{pgfscope}%
\begin{pgfscope}%
\pgfsetbuttcap%
\pgfsetroundjoin%
\definecolor{currentfill}{rgb}{0.000000,0.000000,0.000000}%
\pgfsetfillcolor{currentfill}%
\pgfsetlinewidth{0.803000pt}%
\definecolor{currentstroke}{rgb}{0.000000,0.000000,0.000000}%
\pgfsetstrokecolor{currentstroke}%
\pgfsetdash{}{0pt}%
\pgfsys@defobject{currentmarker}{\pgfqpoint{-0.048611in}{0.000000in}}{\pgfqpoint{-0.000000in}{0.000000in}}{%
\pgfpathmoveto{\pgfqpoint{-0.000000in}{0.000000in}}%
\pgfpathlineto{\pgfqpoint{-0.048611in}{0.000000in}}%
\pgfusepath{stroke,fill}%
}%
\begin{pgfscope}%
\pgfsys@transformshift{0.538891in}{1.456063in}%
\pgfsys@useobject{currentmarker}{}%
\end{pgfscope}%
\end{pgfscope}%
\begin{pgfscope}%
\definecolor{textcolor}{rgb}{0.000000,0.000000,0.000000}%
\pgfsetstrokecolor{textcolor}%
\pgfsetfillcolor{textcolor}%
\pgftext[x=0.294444in, y=1.417483in, left, base]{\color{textcolor}{\rmfamily\fontsize{8.000000}{9.600000}\selectfont\catcode`\^=\active\def^{\ifmmode\sp\else\^{}\fi}\catcode`\%=\active\def
\end{pgfscope}%
\begin{pgfscope}%
\pgfpathrectangle{\pgfqpoint{0.538891in}{0.548472in}}{\pgfqpoint{3.487500in}{1.925000in}}%
\pgfusepath{clip}%
\pgfsetrectcap%
\pgfsetroundjoin%
\pgfsetlinewidth{0.803000pt}%
\definecolor{currentstroke}{rgb}{0.690196,0.690196,0.690196}%
\pgfsetstrokecolor{currentstroke}%
\pgfsetdash{}{0pt}%
\pgfpathmoveto{\pgfqpoint{0.538891in}{1.909858in}}%
\pgfpathlineto{\pgfqpoint{4.026391in}{1.909858in}}%
\pgfusepath{stroke}%
\end{pgfscope}%
\begin{pgfscope}%
\pgfsetbuttcap%
\pgfsetroundjoin%
\definecolor{currentfill}{rgb}{0.000000,0.000000,0.000000}%
\pgfsetfillcolor{currentfill}%
\pgfsetlinewidth{0.803000pt}%
\definecolor{currentstroke}{rgb}{0.000000,0.000000,0.000000}%
\pgfsetstrokecolor{currentstroke}%
\pgfsetdash{}{0pt}%
\pgfsys@defobject{currentmarker}{\pgfqpoint{-0.048611in}{0.000000in}}{\pgfqpoint{-0.000000in}{0.000000in}}{%
\pgfpathmoveto{\pgfqpoint{-0.000000in}{0.000000in}}%
\pgfpathlineto{\pgfqpoint{-0.048611in}{0.000000in}}%
\pgfusepath{stroke,fill}%
}%
\begin{pgfscope}%
\pgfsys@transformshift{0.538891in}{1.909858in}%
\pgfsys@useobject{currentmarker}{}%
\end{pgfscope}%
\end{pgfscope}%
\begin{pgfscope}%
\definecolor{textcolor}{rgb}{0.000000,0.000000,0.000000}%
\pgfsetstrokecolor{textcolor}%
\pgfsetfillcolor{textcolor}%
\pgftext[x=0.294444in, y=1.871278in, left, base]{\color{textcolor}{\rmfamily\fontsize{8.000000}{9.600000}\selectfont\catcode`\^=\active\def^{\ifmmode\sp\else\^{}\fi}\catcode`\%=\active\def
\end{pgfscope}%
\begin{pgfscope}%
\pgfpathrectangle{\pgfqpoint{0.538891in}{0.548472in}}{\pgfqpoint{3.487500in}{1.925000in}}%
\pgfusepath{clip}%
\pgfsetrectcap%
\pgfsetroundjoin%
\pgfsetlinewidth{0.803000pt}%
\definecolor{currentstroke}{rgb}{0.690196,0.690196,0.690196}%
\pgfsetstrokecolor{currentstroke}%
\pgfsetdash{}{0pt}%
\pgfpathmoveto{\pgfqpoint{0.538891in}{2.363654in}}%
\pgfpathlineto{\pgfqpoint{4.026391in}{2.363654in}}%
\pgfusepath{stroke}%
\end{pgfscope}%
\begin{pgfscope}%
\pgfsetbuttcap%
\pgfsetroundjoin%
\definecolor{currentfill}{rgb}{0.000000,0.000000,0.000000}%
\pgfsetfillcolor{currentfill}%
\pgfsetlinewidth{0.803000pt}%
\definecolor{currentstroke}{rgb}{0.000000,0.000000,0.000000}%
\pgfsetstrokecolor{currentstroke}%
\pgfsetdash{}{0pt}%
\pgfsys@defobject{currentmarker}{\pgfqpoint{-0.048611in}{0.000000in}}{\pgfqpoint{-0.000000in}{0.000000in}}{%
\pgfpathmoveto{\pgfqpoint{-0.000000in}{0.000000in}}%
\pgfpathlineto{\pgfqpoint{-0.048611in}{0.000000in}}%
\pgfusepath{stroke,fill}%
}%
\begin{pgfscope}%
\pgfsys@transformshift{0.538891in}{2.363654in}%
\pgfsys@useobject{currentmarker}{}%
\end{pgfscope}%
\end{pgfscope}%
\begin{pgfscope}%
\definecolor{textcolor}{rgb}{0.000000,0.000000,0.000000}%
\pgfsetstrokecolor{textcolor}%
\pgfsetfillcolor{textcolor}%
\pgftext[x=0.294444in, y=2.325073in, left, base]{\color{textcolor}{\rmfamily\fontsize{8.000000}{9.600000}\selectfont\catcode`\^=\active\def^{\ifmmode\sp\else\^{}\fi}\catcode`\%=\active\def
\end{pgfscope}%
\begin{pgfscope}%
\definecolor{textcolor}{rgb}{0.000000,0.000000,0.000000}%
\pgfsetstrokecolor{textcolor}%
\pgfsetfillcolor{textcolor}%
\pgftext[x=0.238889in,y=1.510972in,,bottom,rotate=90.000000]{\color{textcolor}{\rmfamily\fontsize{10.000000}{12.000000}\selectfont\catcode`\^=\active\def^{\ifmmode\sp\else\^{}\fi}\catcode`\%=\active\def
\end{pgfscope}%
\begin{pgfscope}%
\pgfsetrectcap%
\pgfsetmiterjoin%
\pgfsetlinewidth{0.803000pt}%
\definecolor{currentstroke}{rgb}{0.000000,0.000000,0.000000}%
\pgfsetstrokecolor{currentstroke}%
\pgfsetdash{}{0pt}%
\pgfpathmoveto{\pgfqpoint{0.538891in}{0.548472in}}%
\pgfpathlineto{\pgfqpoint{0.538891in}{2.473472in}}%
\pgfusepath{stroke}%
\end{pgfscope}%
\begin{pgfscope}%
\pgfsetrectcap%
\pgfsetmiterjoin%
\pgfsetlinewidth{0.803000pt}%
\definecolor{currentstroke}{rgb}{0.000000,0.000000,0.000000}%
\pgfsetstrokecolor{currentstroke}%
\pgfsetdash{}{0pt}%
\pgfpathmoveto{\pgfqpoint{4.026391in}{0.548472in}}%
\pgfpathlineto{\pgfqpoint{4.026391in}{2.473472in}}%
\pgfusepath{stroke}%
\end{pgfscope}%
\begin{pgfscope}%
\pgfsetrectcap%
\pgfsetmiterjoin%
\pgfsetlinewidth{0.803000pt}%
\definecolor{currentstroke}{rgb}{0.000000,0.000000,0.000000}%
\pgfsetstrokecolor{currentstroke}%
\pgfsetdash{}{0pt}%
\pgfpathmoveto{\pgfqpoint{0.538891in}{0.548472in}}%
\pgfpathlineto{\pgfqpoint{4.026391in}{0.548472in}}%
\pgfusepath{stroke}%
\end{pgfscope}%
\begin{pgfscope}%
\pgfsetrectcap%
\pgfsetmiterjoin%
\pgfsetlinewidth{0.803000pt}%
\definecolor{currentstroke}{rgb}{0.000000,0.000000,0.000000}%
\pgfsetstrokecolor{currentstroke}%
\pgfsetdash{}{0pt}%
\pgfpathmoveto{\pgfqpoint{0.538891in}{2.473472in}}%
\pgfpathlineto{\pgfqpoint{4.026391in}{2.473472in}}%
\pgfusepath{stroke}%
\end{pgfscope}%
\begin{pgfscope}%
\pgfsetbuttcap%
\pgfsetmiterjoin%
\definecolor{currentfill}{rgb}{0.300000,0.300000,0.300000}%
\pgfsetfillcolor{currentfill}%
\pgfsetfillopacity{0.500000}%
\pgfsetlinewidth{0.240900pt}%
\definecolor{currentstroke}{rgb}{0.300000,0.300000,0.300000}%
\pgfsetstrokecolor{currentstroke}%
\pgfsetstrokeopacity{0.500000}%
\pgfsetdash{}{0pt}%
\pgfpathmoveto{\pgfqpoint{0.946713in}{2.586867in}}%
\pgfpathlineto{\pgfqpoint{3.674124in}{2.586867in}}%
\pgfpathquadraticcurveto{\pgfqpoint{3.696346in}{2.586867in}}{\pgfqpoint{3.696346in}{2.609089in}}%
\pgfpathlineto{\pgfqpoint{3.696346in}{2.752917in}}%
\pgfpathquadraticcurveto{\pgfqpoint{3.696346in}{2.775139in}}{\pgfqpoint{3.674124in}{2.775139in}}%
\pgfpathlineto{\pgfqpoint{0.946713in}{2.775139in}}%
\pgfpathquadraticcurveto{\pgfqpoint{0.924491in}{2.775139in}}{\pgfqpoint{0.924491in}{2.752917in}}%
\pgfpathlineto{\pgfqpoint{0.924491in}{2.609089in}}%
\pgfpathquadraticcurveto{\pgfqpoint{0.924491in}{2.586867in}}{\pgfqpoint{0.946713in}{2.586867in}}%
\pgfpathlineto{\pgfqpoint{0.946713in}{2.586867in}}%
\pgfpathclose%
\pgfusepath{stroke,fill}%
\end{pgfscope}%
\begin{pgfscope}%
\pgfsetbuttcap%
\pgfsetmiterjoin%
\definecolor{currentfill}{rgb}{1.000000,1.000000,1.000000}%
\pgfsetfillcolor{currentfill}%
\pgfsetlinewidth{0.240900pt}%
\definecolor{currentstroke}{rgb}{0.800000,0.800000,0.800000}%
\pgfsetstrokecolor{currentstroke}%
\pgfsetdash{}{0pt}%
\pgfpathmoveto{\pgfqpoint{0.918935in}{2.614645in}}%
\pgfpathlineto{\pgfqpoint{3.646346in}{2.614645in}}%
\pgfpathquadraticcurveto{\pgfqpoint{3.668568in}{2.614645in}}{\pgfqpoint{3.668568in}{2.636867in}}%
\pgfpathlineto{\pgfqpoint{3.668568in}{2.780694in}}%
\pgfpathquadraticcurveto{\pgfqpoint{3.668568in}{2.802917in}}{\pgfqpoint{3.646346in}{2.802917in}}%
\pgfpathlineto{\pgfqpoint{0.918935in}{2.802917in}}%
\pgfpathquadraticcurveto{\pgfqpoint{0.896713in}{2.802917in}}{\pgfqpoint{0.896713in}{2.780694in}}%
\pgfpathlineto{\pgfqpoint{0.896713in}{2.636867in}}%
\pgfpathquadraticcurveto{\pgfqpoint{0.896713in}{2.614645in}}{\pgfqpoint{0.918935in}{2.614645in}}%
\pgfpathlineto{\pgfqpoint{0.918935in}{2.614645in}}%
\pgfpathclose%
\pgfusepath{stroke,fill}%
\end{pgfscope}%
\begin{pgfscope}%
\pgfsetbuttcap%
\pgfsetmiterjoin%
\definecolor{currentfill}{rgb}{0.129412,0.901961,0.756863}%
\pgfsetfillcolor{currentfill}%
\pgfsetlinewidth{0.000000pt}%
\definecolor{currentstroke}{rgb}{0.000000,0.000000,0.000000}%
\pgfsetstrokecolor{currentstroke}%
\pgfsetstrokeopacity{0.000000}%
\pgfsetdash{}{0pt}%
\pgfpathmoveto{\pgfqpoint{0.941157in}{2.680694in}}%
\pgfpathlineto{\pgfqpoint{1.163380in}{2.680694in}}%
\pgfpathlineto{\pgfqpoint{1.163380in}{2.758472in}}%
\pgfpathlineto{\pgfqpoint{0.941157in}{2.758472in}}%
\pgfpathlineto{\pgfqpoint{0.941157in}{2.680694in}}%
\pgfpathclose%
\pgfusepath{fill}%
\end{pgfscope}%
\begin{pgfscope}%
\definecolor{textcolor}{rgb}{0.000000,0.000000,0.000000}%
\pgfsetstrokecolor{textcolor}%
\pgfsetfillcolor{textcolor}%
\pgftext[x=1.252268in,y=2.680694in,left,base]{\color{textcolor}{\rmfamily\fontsize{8.000000}{9.600000}\selectfont\catcode`\^=\active\def^{\ifmmode\sp\else\^{}\fi}\catcode`\%=\active\def
\end{pgfscope}%
\begin{pgfscope}%
\pgfsetbuttcap%
\pgfsetmiterjoin%
\definecolor{currentfill}{rgb}{0.152941,0.556863,0.647059}%
\pgfsetfillcolor{currentfill}%
\pgfsetlinewidth{0.000000pt}%
\definecolor{currentstroke}{rgb}{0.000000,0.000000,0.000000}%
\pgfsetstrokecolor{currentstroke}%
\pgfsetstrokeopacity{0.000000}%
\pgfsetdash{}{0pt}%
\pgfpathmoveto{\pgfqpoint{1.808986in}{2.680694in}}%
\pgfpathlineto{\pgfqpoint{2.031208in}{2.680694in}}%
\pgfpathlineto{\pgfqpoint{2.031208in}{2.758472in}}%
\pgfpathlineto{\pgfqpoint{1.808986in}{2.758472in}}%
\pgfpathlineto{\pgfqpoint{1.808986in}{2.680694in}}%
\pgfpathclose%
\pgfusepath{fill}%
\end{pgfscope}%
\begin{pgfscope}%
\definecolor{textcolor}{rgb}{0.000000,0.000000,0.000000}%
\pgfsetstrokecolor{textcolor}%
\pgfsetfillcolor{textcolor}%
\pgftext[x=2.120097in,y=2.680694in,left,base]{\color{textcolor}{\rmfamily\fontsize{8.000000}{9.600000}\selectfont\catcode`\^=\active\def^{\ifmmode\sp\else\^{}\fi}\catcode`\%=\active\def
\end{pgfscope}%
\begin{pgfscope}%
\pgfsetbuttcap%
\pgfsetmiterjoin%
\definecolor{currentfill}{rgb}{0.121569,0.258824,0.529412}%
\pgfsetfillcolor{currentfill}%
\pgfsetlinewidth{0.000000pt}%
\definecolor{currentstroke}{rgb}{0.000000,0.000000,0.000000}%
\pgfsetstrokecolor{currentstroke}%
\pgfsetstrokeopacity{0.000000}%
\pgfsetdash{}{0pt}%
\pgfpathmoveto{\pgfqpoint{2.827666in}{2.680694in}}%
\pgfpathlineto{\pgfqpoint{3.049888in}{2.680694in}}%
\pgfpathlineto{\pgfqpoint{3.049888in}{2.758472in}}%
\pgfpathlineto{\pgfqpoint{2.827666in}{2.758472in}}%
\pgfpathlineto{\pgfqpoint{2.827666in}{2.680694in}}%
\pgfpathclose%
\pgfusepath{fill}%
\end{pgfscope}%
\begin{pgfscope}%
\definecolor{textcolor}{rgb}{0.000000,0.000000,0.000000}%
\pgfsetstrokecolor{textcolor}%
\pgfsetfillcolor{textcolor}%
\pgftext[x=3.138777in,y=2.680694in,left,base]{\color{textcolor}{\rmfamily\fontsize{8.000000}{9.600000}\selectfont\catcode`\^=\active\def^{\ifmmode\sp\else\^{}\fi}\catcode`\%=\active\def
\end{pgfscope}%
\end{pgfpicture}%
\makeatother%
\endgroup%